\documentclass[letterpaper]{IEEEtran}
\usepackage{csm16}
\usepackage[margin=1in]{geometry}
\usepackage[ascii]{inputenc}
\usepackage[T1]{fontenc}
\usepackage{lmodern}
\usepackage{amsmath, amssymb}
\usepackage{mdframed}
\usepackage{graphicx}
\usepackage{epstopdf}
\usepackage{overpic}
\usepackage{subcaption}
\usepackage{color}
\usepackage{rotating}
\usepackage{setspace}
\usepackage{palatino} 
\usepackage{mdframed}
\newcounter{frameNum}
\usepackage{booktabs}
\usepackage{array}
\usepackage{multirow,multicol}

\usepackage[innercaption]{sidecap}

\newcounter{seqn}
\newcounter{counteqn}
\newcounter{countfig}
\newcounter{sidefig}
\makeatletter
\newenvironment{Sidebar}[2]{ 
	\protected@edef\@currentlabelname{#2}
	\protected@edef\@currentlabel{#2}
	\setcounter{table}{0}
	\setcounter{seqn}{\value{equation}}
	\setcounter{equation}{\value{counteqn}}
	
	\setcounter{sidefig}{\value{figure}}
	\setcounter{figure}{\value{countfig}}
	\renewcommand{\thefigure}{F\arabic{figure}}
	\renewcommand{\theequation}{S\arabic{equation}}
	\renewcommand{\thetable}{T\arabic{table}}
	\begin{mdframed}[backgroundcolor=lightgray!20,shadow=true,roundcorner=8pt,frametitle={#1}]
	}{	 
\end{mdframed}
\setcounter{counteqn}{\value{equation}}
\setcounter{equation}{\value{seqn}}
\setcounter{countfig}{\value{figure}}
\setcounter{figure}{\value{sidefig}}
\addtocounter{frameNum}{1}
}
\makeatother

\usepackage{textcomp,xcolor}
\definecolor{MatlabCellColour}{RGB}{250,250,250}
\definecolor{MatPurp}{rgb}{.625,.1406,.9375}
\usepackage{listings}

\lstdefinestyle{customc}{
  belowcaptionskip=.25\baselineskip,
  breaklines=true,
  frame=L,
  xleftmargin=\parindent,
  language=Matlab,
  showstringspaces=false,
  basicstyle=\small\ttfamily,
  keywordstyle=\bfseries\color{white!30!black},
  identifierstyle=\color{blue},  
  commentstyle=\itshape\color{green!60!black},
  stringstyle=\color{MatPurp},
  backgroundcolor=\color{MatlabCellColour}
 }
 
\newcommand{\ba}{\mathbf{a}}
\newcommand{\bc}{\mathbf{c}}
\newcommand{\be}{\mathbf{e}}
\newcommand{\bx}{\mathbf{x}}
\newcommand{\by}{\mathbf{y}}

\newcommand{\bA}{\mathbf{A}}
\newcommand{\bB}{\mathbf{B}}
\newcommand{\bM}{\mathbf{M}}

\newcommand{\bs}{\mathbf{s}}
\newcommand{\bC}{\mathbf{C}}

\newcommand{\bR}{\mathbf{R}}
\newcommand{\bQ}{\mathbf{Q}}
\newcommand{\bX}{\mathbf{X}}
\newcommand{\bSigma}{\mathbf{\Sigma}}
\newcommand{\bV}{\mathbf{V}}
\newcommand{\bPsi}{\mathbf{\Psi}}

\newcommand{\bpsi}{\boldsymbol\psi}
\newcommand{\bTheta}{\mathbf{\Theta}}
\newcommand{\reals}{\mathbb{R}}
\newcommand{\ie}{\emph{i.e.,}}

\DeclareMathOperator*{\argmax}{arg\rm{}max}
\DeclareMathOperator*{\argmin}{arg\rm{}min}
\DeclareMathOperator*{\logdet}{log\rm{}det}
\DeclareMathOperator*{\subjto}{subject~to}

\newcommand\revision[3][]{%
\ifstrempty{#1}{}
{\textcolor{black}{}}
\ifstrempty{#3}{{\color{black}{#2}}}%
{{\color{black}{#2}}}
}
\newcommand{\btheta}{\boldsymbol{\theta}}
\newcommand{\ind}{\gamma}
\newcommand{\Ind}{\boldsymbol{\gamma}}
\usepackage{algorithm}
\usepackage[noend]{algpseudocode}

\IEEEoverridecommandlockouts
\graphicspath{{figures/}}

\title{\Large Data-Driven Sparse Sensor Placement for Reconstruction}
\author{Krithika Manohar$^*$, Bingni W. Brunton, J. Nathan Kutz, and Steven L. Brunton\\
	{\small $^*$Corresponding author: kmanohar@uw.edu}}
	
\newif\ifPDF \ifx\pdfoutput\undefined\PDFfalse \else\ifnum\pdfoutput > 0\PDFtrue \else\PDFfalse \fi \fi
\ifPDF 
\usepackage[pdftex, plainpages = false, colorlinks=true, linkcolor=black, citecolor = green!50!blue, urlcolor = blue, filecolor=black, pagebackref=false, hypertexnames=false,  pdfpagelabels ]{hyperref}
\fi

\begin{document}
\twocolumn[
\begin{@twocolumnfalse}
	\maketitle
	\CSMsetup	
	\vspace{-4em}
	\begin{abstract}
Optimal sensor placement is a central challenge in the design, prediction, estimation, and control of high-dimensional systems. High-dimensional states can often leverage a latent low-dimensional representation, and this inherent compressibility enables sparse sensing. This article explores optimized sensor placement for signal reconstruction based on a tailored library of features extracted from training data.  Sparse point sensors are discovered using the singular value decomposition and QR pivoting, which are two ubiquitous matrix computations that underpin modern linear dimensionality reduction. Sparse sensing in a tailored basis is contrasted with compressed sensing, a universal signal recovery method in which an unknown signal is reconstructed via a sparse representation in a universal basis.  Although compressed sensing can recover a wider class of signals, we demonstrate the benefits of exploiting known patterns in data with optimized sensing. In particular, drastic reductions in the required number of sensors and improved reconstruction are observed in examples ranging from facial images to fluid vorticity fields. Principled sensor placement may be critically enabling when sensors are costly and provides faster state estimation for low-latency, high-bandwidth control.
	\end{abstract}
\end{@twocolumnfalse}]


\begin{tabbing}
  XXXXXXXX \= \kill
   \textit{Scalars} \> \sc{Notation} \\
  $n$ \> State dimension \\
  $m$ \> Number of snapshots \\
    $p$ \> Number of sensors (measurements) \\  
  $r$ \> Intrinsic rank of tailored basis $\bPsi_r$\\ 
    $K$ \> Sparsity of state in universal basis $\bPsi$\\
    $\eta$ \> Variance of zero-mean sensor noise \\[5pt]
  \textit{Vectors} \> \\
  $\bx\in\reals^n$ \> High-dimensional state \\  
  $\by\in\reals^p$ \> Measurements of state \\
  $\ba\in\reals^r$ \> Tailored basis coefficients \\
  $\be_j\in\reals^n$ \> Canonical basis vectors for $\reals^n$ \\
  $\bs\in\reals^n$ \> $K$-sparse basis coefficients\\
  $\Ind\in\mathbb{N}^p$ \> Sensor placement indices \\
  $\bpsi\in\reals^n$ \> POD modes (columns of $\bPsi_r$) \\ 
  $\btheta\in\reals^{1\times r}$ \> Rows of $\bTheta$ \\ [5pt]
   \textit{Matrices} \> \\
  $\bC\in\reals^{p\times n}$ \> Measurement matrix \\
  $\bQ$ \> Unitary QR factor matrix \\
  $\bR$ \> Upper triangular QR factor matrix\\
  $\bPsi\in\reals^{n\times n}$ \> Universal basis \\
  $\bPsi_r\in\reals^{n\times r}$ \> Tailored basis of rank $r$ \\
  $\bTheta=\bC\bPsi$ \> Product of measurement and basis  \\
  $\bX\in\mathbb{R}^{n\times m}$ \> Data matrix with $m$ snapshots\\
 \end{tabbing}

\vspace{-2em}

\begin{figure*}
\begin{Sidebar}{Mathematical formulation of sensor selection}{Sidebar: Mathematical formulation of sensor selection}
\label{Sidebar1}

Many physical systems are described by a high-dimensional state $\bx \in\reals^n$, yet the dynamics evolve on a low-dimensional attractor that can be leveraged for prediction and control.  
Thus, a state $\bx$ that evolves according to nonlinear dynamics $\dot\bx(t) = \mathbf{f}(\bx(t))$ will often have a compact representation in a transform basis $\bPsi$.  
In a universal basis $\bPsi\in\mathbb{R}^{n\times n}$, such as Fourier or wavelet bases, $\bx$ may have a \emph{sparse} representation
\begin{align}
\label{eqn:stateCS}
\bx &= \bPsi\bs &\bs\in\reals^n,
\end{align}
where $\bs$ is a sparse vector indicating which few modes of $\bPsi$ are active.  In a tailored basis $\bPsi_r\in\mathbb{R}^{n\times r}$, such as a proper orthogonal decomposition (POD) basis, $\bx$ may have a \emph{low-rank} representation
\begin{align}
\label{eqn:state}
\bx &= \bPsi_r\ba & \ba\in\reals^r.
\end{align}

The central challenge in this work is to design a measurement matrix $\bC\in\mathbb{R}^{p\times n}$ consisting of a small number $(p\ll n)$ of optimized measurements
\begin{align}
\label{eqn:observations}
\by &= \bC\bx &\by\in\reals^p,
\end{align}
that facilitate accurate reconstruction of either $\bs$ or $\ba$, and hence $\bx$.  
Combining~\eqref{eqn:stateCS} and~\eqref{eqn:observations} yields
\begin{equation}
\by = (\bC\bPsi)\bs = \bTheta\bs,\label{Eq:LinSysCS}
\end{equation}
which is referred to as the \emph{compressed sensing} problem, while combining~\eqref{eqn:state} and~\eqref{eqn:observations} yields
\begin{equation}
\by = (\bC\bPsi_r)\ba = \bTheta\ba.\label{Eq:LinSysPOD}
\end{equation}
In both cases, effective measurements $\bC$ given a basis $\bPsi$ or $\bPsi_r$ are chosen so that the operator $\bTheta$ is well-conditioned for signal reconstruction. 
Thus, it is possible to solve for the sparse coefficients $\bs$ or the low-rank coefficients $\ba$ given the measurements $\by$, either by $\ell_1$ minimization in~\eqref{Eq:LinSysCS} or pseudoinverse of $\bTheta$ in~\eqref{Eq:LinSysPOD}, respectively.  
The goal of this work is to optimize the measurements in $\bC$. 
Moreover, in many physical applications, it is desired that $\bC$ consists of rows of the identify matrix, corresponding to individual point sensors of individual components of $\bx$.  
\end{Sidebar}
\end{figure*}


Optimal sensor and actuator placement is an important unsolved problem in control theory.  
Nearly every downstream control decision is affected by these sensor/actuator locations, but determining optimal locations amounts to an intractable brute force search among the combinatorial possibilities.
Indeed, there are ${n\choose p}=\frac{n!}{(n-p)!p!}$ possible choices of $p$ point sensors out of an $n$ dimensional state $\bx$.  
Determining optimal sensor and actuator placement in general, even for linear feedback control, is an open challenge.  
Instead, sensor and actuator locations are routinely chosen according to heuristics and intuition.  
For moderate sized search spaces, the sensor placement problem has well-known model-based solutions using optimal experiment design~\cite{Boyd2004convexbook,Joshi2009ieee}, information theoretic and Bayesian criteria~\cite{Caselton1984spl,Krause2008jmlr,Lindley1956ams,Sebastiani2000jrss,Paninski2005nc} .
We explore how to design optimal sensor locations for signal reconstruction in a framework that scales to arbitrarily large problems, leveraging modern techniques in machine learning and sparse sampling. 
Reducing the number of sensors through principled selection may be critically enabling when sensors are costly, and may also enable faster state estimation for low latency, high bandwidth control. 

\revision[]{This article explores optimized sensor placement for signal reconstruction based on a tailored library of features extracted from training data. 
In this paradigm, optimized sparse sensors are computed using a powerful sampling scheme based on the matrix QR factorization and singular value decomposition. 
Both procedures are natively implemented in modern scientific computing software, and Matlab code supplements are provided for all examples in this paper~\cite{manohar2017code}.
These data-driven computations are more efficient and easier to implement than the convex optimization methods used for sensor placement in classical design of experiments. 
In addition, data-driven sensing in a tailored basis is contrasted with compressed sensing, a universal signal recovery method in which an unknown signal is reconstructed using a sparse representation in a universal basis. 
Although compressed sensing can recover a wider class of signals, we demonstrate the benefits of exploiting known patterns in data with optimized sensing. 
In particular, drastic reductions in the required number of sensors and improved reconstruction are observed in examples ranging from facial images to fluid vorticity fields.
The overarching signal reconstruction problem is formulated in~``\nameref{Sidebar1}".
}{
This paper provides a tutorial overview of current sparse sampling methods for sensor placement and reconstruction of structured signals. 
We also connect and equate certain sampling strategies with analogues in the design of experiments literature.  
Near-optimal sensor locations are obtained using fast greedy procedures that scale well with large signal dimension. This work also generalizes and extends a powerful sampling scheme based on the matrix QR factorization and demonstrates its broad applicability to image and fluid flow reconstruction as well as polynomial interpolation.
The overarching sensor placement problem is summarized in~``\nameref{Sidebar1}", and Matlab code is provided for all examples.}

There are myriad complex systems that would benefit from principled, scaleable sensor and actuator placement, including fluid flow control~\cite{Brunton2015amr}, power grid optimization~\cite{susuki2012nonlinear}, epidemiological modeling and suppression~\cite{Proctor2015ih}, bio-regulatory network monitoring and control~\cite{Plenge2013}, and high-performance computing~\cite{Carlberg2011ijnme,Carlberg2013jcp}, to name only a few.  
{In applications where individual sensors are expensive, reducing the number of sensors through principled design may be critically enabling.  
In applications where fast decisions are required, such as in high performance computing or feedback control, computations may be accelerated by minimizing the number of sensors required.  
In other words, low-dimensional computations may be performed directly in the sensor space.}

Scaleable optimization of sensor and actuator placement is a grand challenge problem, with tremendous potential impact and considerable mathematical depth.  
With existing mathematical machinery, optimal placement can only be determined in general using a brute-force combinatorial search.  
Although this approach has been successful in small-scale problems~\cite{Chen:2011}, a combinatorial search does not scale well to larger problems.  Moore's law of exponentially increasing computer power cannot keep pace with this combinatorial growth in complexity.   

Despite the challenges of sensing and actuation in a high-dimensional, possibly nonlinear dynamical system, there are promising indicators that this problem may be tractable with modern techniques.  
High-dimensional systems, such as are found in fluids, epidemiology, neuroscience, and the power grid, typically exhibit dominant coherent structures that evolve on a low-dimensional attractor.  
{Indeed, much of the success of modern machine learning rests on the ability to identify and take advantage of patterns and features in high-dimensional data.  }
These low-dimensional patterns are often identified using dimensionality reduction techniques~\cite{Kutz:2013} such as the proper orthogonal decomposition (POD)~\cite{Sirovich:1987,berkooz1993proper,HLBR_turb}, which is a variant of principal component analysis (PCA), or more recently via dynamic mode decomposition (DMD)~\cite{Schmid2010jfm,Rowley2009jfm,Tu2014jcd,Kutz2016book}, diffusion maps~\cite{Coifman2005pnas,Coifman2006acha,Coifman2008mmas}, etc.   
In control theory, balanced truncation~\cite{Moore1981ieeetac}, balanced proper orthogonal decomposition (BPOD)~\cite{Willcox2002aiaaj,Rowley2005ijbc}, and the eigensystem realization algorithm (ERA)~\cite{ERA:1985}, have been successfully applied to obtain control-oriented reduced-order models for many high-dimensional systems.  

In addition to advances in dimensionality reduction, key developments in optimization, compression, and the geometry of sparse vectors in high-dimensional spaces are providing powerful new techniques to obtain approximate solutions to NP-hard, combinatorially difficult problems in scaleable convex optimization architectures.  
For example, compressed sensing~\cite{Candes2006cpam,Donoho2006ieeetit,Baraniuk2007ieeespm} provides convex algorithms to solve the combinatorial sparse signal reconstruction problem with high probability.  
Ideas from compressed sensing have been used to determine the optimal sensor locations for categorical decisions based on high-dimensional data~\cite{Brunton2016siap}. 
Recently, compressed sensing, sparsity-promoting algorithms such as the \emph{lasso} regression~\cite{Tibshirani1996lasso,Hastie2009book,James2013book}, and machine learning have been increasingly applied to characterize and control dynamical systems~\cite{Schaeffer2013pnas,Ozolicnvs2013pnas,Brunton2014siads,Proctor2014epj,Brunton2016pnas,Mangan2016ieee,Loiseau2016arxiv,Rudy2017sciadv,Duriez2016book,Schaeffer2017prsa}.  
These techniques have been effective in modeling high-dimensional fluid systems using POD~\cite{Bai2014aiaa} and DMD~\cite{Brunton2015jcd,Tu2014ef,Jovanovic2014pof,Gueniat2015pof}.  
Information criteria~\cite{Akaike1973,Bayesian1978} has also been leveraged for the sparse identification of nonlinear dynamics~\cite{Brunton2016pnas}, as in~~\cite{Giannakis:2013,Mangan2017arxiv}, and may also be useful for sensor placement. 

Thus, key advances in two fields are fundamentally changing our approach to the acquisition and analysis of data from complex dynamical systems:  1) \emph{machine learning}, which exploits patterns in data for low-dimensional representations, and 2) \emph{sparse sampling}, where a full signal can be reconstructed from a small subset of measurements.  The combination of machine learning and sparse sampling is synergistic, in that underlying low-rank representations facilitate sparse measurements.  
Exploiting coherent structures underlying a large state space allows us to estimate and control systems with few measurements and sparse actuation.  Low-dimensional, data-driven sensing and control is inspired in part by the high performance exhibited by biological organisms, such as an insect that performs robust, high-performance flight control in a turbulent fluid environment with minimal sensing and low-latency control~\cite{Manohar2016jfs}.  They provide proof-by-existence that it is possible to assimilate sparse measurements and perform low-dimensional computations to interact with coherent structures in a high-dimensional system (i.e., a turbulent fluid).

Here we explore two competing perspectives on high-dimensional signal reconstruction:  1)  the use of compressed sensing based on random measurements in a universal encoding basis, and 2) the use of highly specialized sensors for reconstruction in a tailored basis, such as POD or DMD.  
These choices are also discussed in the context of feedback control.  
Many competing factors impact control design, and a chief consideration is the latency in making a control decision, with large latency imposing limitations on robust performance~\cite{sp:book,dp:book}.  
Thus, for systems with fast dynamics and complex coherent structures, it is important to make control decisions quickly based on efficient low-order models, with sensors and actuators placed strategically to gather information and exploit sensitivities in the dynamics.

 \subsection{Extensions to dynamics, control, and multiscale physics}
Data-driven sensor selection is generally used for instantaneous full-state reconstruction, despite the fact that many signals are generated by a dynamical system~\cite{guckenheimer_holmes,HLBR_turb}.  
Even in reduced-order models, sensors are typically used to estimate nonlinear terms instantaneously without taking advantage of the underlying dynamics.  
However, it is well known that for linear control systems~\cite{dp:book,sp:book}, the high-dimensional state may be reconstructed with few sensors, if not a single sensor, by leveraging the time history in conjunction with a model of the dynamics, as exemplified by the Kalman filter~\cite{Kalman1960jfe,Welch1995book}.  
In dynamic estimation and control, prior placement of sensors and actuators is generally assumed. 
Extending the sensor placement optimization to the model reduction~\cite{Moore1981ieeetac,Willcox2002aiaaj,Rowley2005ijbc} and system identification~\cite{ERA:1985,ljung:book,Juang1991nasatm,Juang1994book} of linear control systems is an important avenue of ongoing work.  
In particular, sensors and actuators may be chosen to increase the volume of the controllability and observability Gramians, related to the original balanced truncation literature~\cite{Moore1981ieeetac}.  
More generally, sensor and actuator placement may be optimized for robustness~\cite{Doyle:1978,Doyle:1981}, or for network control and consensus problems~\cite{Leonard2001cdc,Olfati2004ieeetac,Doyle2005pnas,Leonard2007pieee,Rahmani:SIAMJCO09}.

The sensor placement algorithms discussed above are rooted firmly in linear algebra, making them readily extensible to linear control systems.  
Recent advances in dynamical systems are providing techniques to embed nonlinear systems in a linear framework through a suitable choice of \emph{measurement functions} of the state, opening up the possibility of optimized sensing for nonlinear systems.  
As early as the 1930s, Koopman demonstrated that a nonlinear system can be rewritten as an infinite-dimensional linear operator on the Hilbert space of measurement functions~\cite{Koopman1931pnas}.  
This perspective did not gain traction until modern computation and data collection capabilities enabled the analysis of large volumes of measurement data.

Modern Koopman theory may drive sensor placement and the selection of nonlinear measurement functions on the sensors to embed nonlinear dynamics in a linear framework for optimal nonlinear estimation and control.  
This approach is consistent with neural control systems, where biological sensor networks (e.g., strain sensors on an insect wing) are processed through nonlinear neural filters before being used for feedback control.  
Much of the modern Koopman operator theory has been recently developed~\cite{Mezic2005nd,Budivsic2009cdc,Budivsic2012chaos,Mezic2013arfm}, and it has been shown that under certain conditions DMD approximates the Koopman operator~\cite{Schmid2010jfm,Rowley2009jfm,Tu2014jcd,Kutz2016book}; sensor fusion is also possible in the Koopman framework~\cite{Williams2015epl}.  
Recently, Koopman analysis has been used to develop nonlinear estimators~\cite{Surana2016cdc,Surana2016nolcos} and controllers~\cite{Brunton2016plosone}, although establishing rigorous connections to control theory is an ongoing effort~\cite{Proctor:2016DMDc,Proctor2016arxiv,Korda2016control}.  
Koopman theory has also been used to analyze chaotic dynamical systems from time-series data~\cite{Giannakis2015arxiv,Brunton2017natcomm}, relying on the Takens embedding~\cite{Takens1981lnm}, which is  related to sensor selection.  

%

Beyond extending sensor selection to nonlinear systems and control, there is a significant opportunity to apply principled sensor selection to multiscale systems.  
Turbulence is an important high-dimensional system that exhibits multiscale phenomena~\cite{HLBR_turb,Majda2010dcds,Osth:2014}.  
Data-driven approaches have been used to characterize turbulent systems~\cite{Brunton2015amr}, including clustering~\cite{Kaiser2014jfm}, network theory~\cite{nair2015network,Taira2016jfm}, DMD-based model reduction~\cite{Noack2015arxiv,tissot2014model}, and local POD subspaces~\cite{Amsallem2012ijnme}, to name a few.  
Recently, a multiresolution DMD has been proposed~\cite{Kutz2016siads}, where a low-dimensional subspace may locally characterize the attractor, despite a high-dimensional global attractor.  
This approach may significantly reduce the number of sensors needed for multiscale problems.   
%

\section{Compressed sensing: Random measurements in a universal basis}\label{Sec:Random}

\begin{figure}
\begin{Sidebar}{Sidebar: Image compression}{Sidebar: Image compression}
	\label{Sidebar4}
Images and audio signals tend to be sparse in Fourier or wavelet bases, providing the foundation of JPEG and MP3 compression, respectively. This is shown schematically in Fig.~\ref{Fig:Compression} using the included Matlab code. 
\begin{lstlisting}[firstline=5,lastline=12]
X = imread('cappuccino','jpeg');  
[nx,ny] = size(X);
S = fft2(X); % X = Psi*s
Ssort = sort(abs(S(:)));
thresh = Ssort(ceil(.95*nx*ny));
bigind = abs(S)>thresh; % big coeffs
Strunc = S.*bigind; % keep big coeffs
Xrecon = uint8(ifft2(Strunc));
\end{lstlisting}
\begin{figure}[H]
		\begin{overpic}[width=\textwidth]{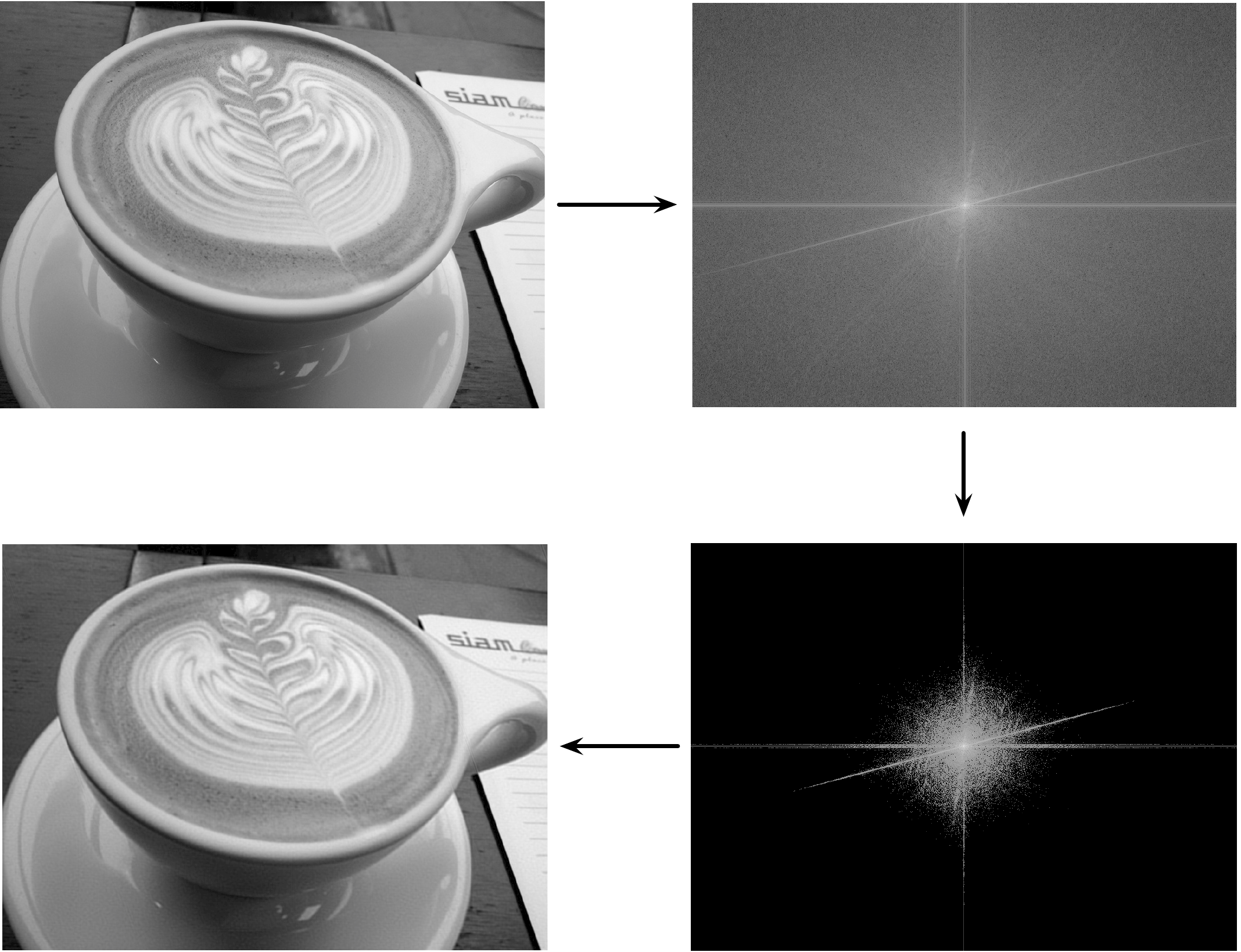}	
			\put(47,63){$\mathcal{F}$}
			\put(46,20){$\mathcal{F}^{-1}$}	
			\put(80,38){Keep $5\,\%$}
		\end{overpic}
	\caption{Fourier image compression. \label{Fig:Compression}} 
\end{figure}	 
\end{Sidebar}	
\end{figure}

The majority of natural signals, such as images and audio, are highly compressible, meaning that when the signal is written in an appropriate coordinate system, only a few basis modes are active.  
These few values corresponding to the large mode amplitudes must be stored for accurate reconstruction, providing a significant reduction compared to the original signal size.  
In other words, in the universal transform basis, the signal may be approximated by a \emph{sparse} vector containing mostly zeros.  
This inherent sparsity of natural signals is central to the mathematical framework of compressed sensing.
Signal compression in the Fourier domain is illustrated on an image example in~``\nameref{Sidebar4}". 
Further, sparse signal recovery using compressed sensing is demonstrated on a sinusoidal example in~``\nameref{sb:cs_example}".
The theory of compressed sensing~\cite{Donoho2006ieeetit,Candes2006cpam,Candes2006ieeetit,Candes2006bieeetit,Candes2006picm,Candes2008ieeespm,Baraniuk2007ieeespm,Baraniuk:2009} inverts this compression paradigm.  
Instead of collecting high-dimensional measurements just to compress and discard most of the information, it may be possible to collect  a low-dimensional subsample or compression of the data and then infer the sparse vector of coefficients in the transformed coordinate system. 
\subsection{Theory of compressed sensing}\label{Sec:CS}
\begin{figure*}
		\centering
		\begin{overpic}[width=\textwidth]{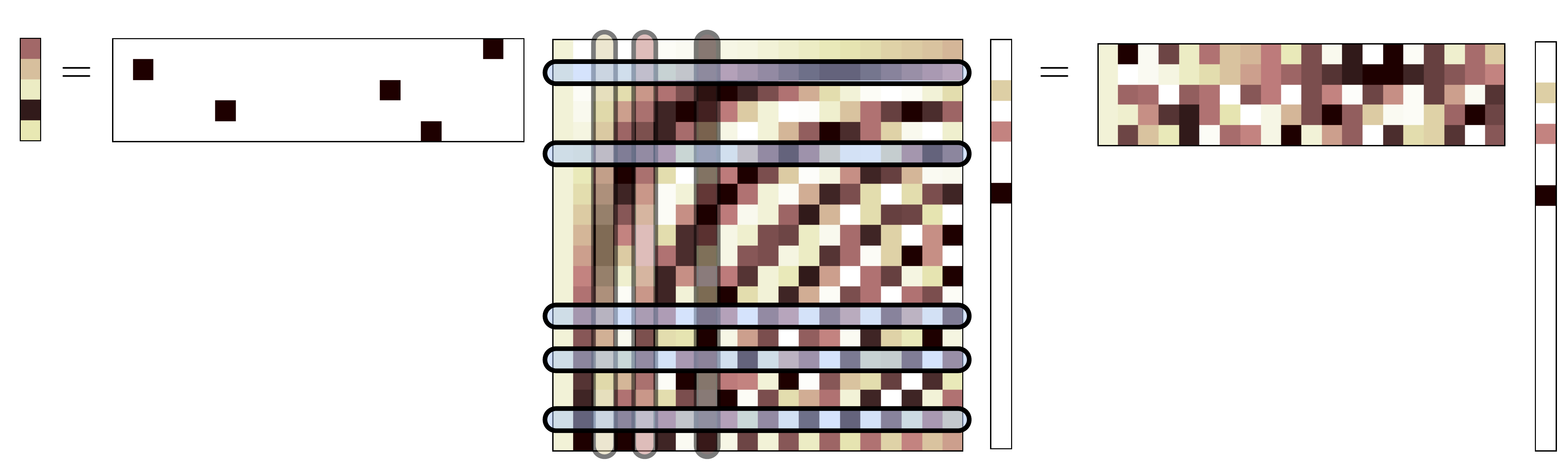}
			\put(1.5,28){$\by$}
			\put(19,28){$\bC$}
			\put(47.5,28){$\bPsi$}
			\put(63.5,28){$\bs$}	
			\put(82,28){$\bTheta$}
			\put(98.2,28){$\bs$}
		\end{overpic}
		\vspace{-.1in}
		\caption{Compressed sensing provides the sparsest solution to an underdetermined linear system.}\label{Fig:CSschematic}
	\end{figure*}
Mathematically, a compressible signal $\mathbf{x}\in\mathbb{R}^n$ may be written as a sparse vector $\mathbf{s}\in\mathbb{R}^n$ in a new basis $\boldsymbol{\Psi}\in\mathbb{R}^{n\times n}$ such that
\begin{equation}
\bx = \bPsi\bs.
\end{equation}
The vector $\bs$ is called $K$-sparse if there are exactly $K$ nonzero elements.  
To be able to represent \emph{any} natural signal, rather than just those from a tailored category, the basis $\bPsi$ must be complete.  

Consider a set of measurements $\by\in\mathbb{R}^p$, obtained via a measurement matrix $\bC\in\mathbb{R}^{p\times n}$, which satisfies
\begin{equation}
\by = \bC \bx = \bC\bPsi\bs = \bTheta\bs.\label{Eq:CS}
\end{equation}
In general, for $p<n$~\eqref{Eq:CS} is underdetermined, and there are infinitely many solutions.  The least least squaressquares (minimum $\|\bs\|_2$) solution is not sparse, and typically yields poor reconstruction.  
Instead, knowing that natural signals are sparse, we seek the sparsest $\bs$ consistent with the measurements $\by$,
\begin{align}
\bs = \argmin_{\bs'}  \|\bs'\|_0 \text{, such that }\by = \bC\bPsi\bs',\label{Eq:L0}
\end{align}
where $\|\bs\|_0$ is the $\ell_0$ pseudo-norm corresponding to the number of non-zero entries of $\bs$.  
Unfortunately, this optimization problem is intractable, requiring a combinatorial brute-force search across all sparse vectors $\bs$.  
A major innovation of compressed sensing is a set of conditions on the measurement matrix $\bC$ that allow the nonconvex $\ell_0$-minimization in~\eqref{Eq:L0} to be relaxed to the convex $\ell_1$-minimization~\cite{Candes2006bieeetit,Donoho:2006b}
\begin{equation}
\bs = \argmin_{\bs'}  \|\bs'\|_1 \text{, such that }\by = \bC\bPsi\bs',\label{Eq:L1}
\end{equation}
where $\|\bs\|_1 = \sum_{k=1}^n|s_k|$.  This formulation is shown schematically in Fig.~\ref{Fig:CSschematic}.
	
For the $\ell_1$-minimization in~\eqref{Eq:L1} to yield the sparsest solution in~\eqref{Eq:L0} with high probability, the measurements $\bC$ must be chosen so that ${\bTheta=\bC\bPsi}$ satisfies a \emph{restricted isometry property} (RIP)
\begin{equation}
(1-\delta_K)\|\bs\|_2^2 \leq \|\bC\bPsi\bs\|_2^2\leq (1+\delta_K)\|\bs\|_2^2,
\end{equation}
where $\delta_K$ is a small positive \emph{restricted isometry} constant~\cite{Candes:2005,Candes2008ieeespm}.  
In particular, there are two conditions on $\bC$ for a RIP to be satisfied for all $K$-sparse vectors $\bs$:
\begin{enumerate}
	\item The measurements $\bC$ must be \emph{incoherent} with respect to the basis $\bPsi$.  
	This incoherence means that the rows of $\bC$ are sufficiently uncorrelated with the columns of $\bPsi$, as quantified by $\mu$
	\begin{equation}
	\mu(\bC,\bPsi) = \sqrt{n} \max_{j,k}|\langle \mathbf{c}_k,\mathbf{\bpsi}_j\rangle|.
	\end{equation}
	Small $\mu$ indicates better incoherent measurements, with an optimal value of $\mu=1$.  
	Here, $\mathbf{c}_k$ denotes the $k$-th row of $\bC$ and $\mathbf{\bpsi}_j$ the $j$-th column of $\bPsi$, both of which are assumed to be normalized.  A more detailed discussion about incoherence and the RIP may be found in~\cite{Baraniuk2007ieeespm,Candes2008ieeespm}.  
	\item The number of measurements $p$ must satisfy~\cite{Candes2006picm,Candes2006bieeetit, Baraniuk2007ieeespm,Candes2008ieeespm,Candes:2010}
	\begin{equation}
	p\sim \mathcal{O}(K\log(n/K)).
	\end{equation}
	The $K\log(n/K)$ term above is generally multiplied by a small constant multiple of the incoherence.  
	Thus, fewer measurements are required if they are less coherent. 
\end{enumerate}

Intuitively, the existence of a RIP implies that the geometry of sparse vectors is preserved through the measurement matrix $\bC\bPsi$.  
Determining the exact constant $\delta_K$ may be extremely challenging in practice, and it tends to be more desirable to characterize the statistical properties of $\delta_K$, as the measurement matrix $\bC$ may be randomly chosen.  
``\nameref{Sidebar_6}" describes why it is not possible to use QR pivot locations as optimized sensors for compressed sensing, since they fail to identify the sparse structure of an unknown signal.

Often, a generic basis such as Fourier or wavelets may be used to represent the signal sparsely.  
Spatially localized measurements (i.e., single pixels in the case of an image) are optimally incoherent with respect to the Fourier basis, so that ${\mu(\bC,\bPsi)=1}$.  
Thus, single pixel measurements are ideal because they excite a broadband frequency response.  
In contrast, a measurement corresponding to a fixed Fourier mode would be uninformative; if the signal is not sparse in this particular frequency, this measurement provides no information about the other Fourier modes.  
For many engineering applications, spatially localized measurements are desirable, as they correspond to physically realizable sensors, such as buoys in the ocean.  

One of the major results of compressed sensing is that random projection measurements of the state (i.e., entries of $\bC$ that are Bernoulli or Gaussian random variables) are incoherent with respect to nearly any generic basis $\bPsi$~\cite{Candes2006picm,Candes2006ieeetit,Donoho2006ieeetit}.  
This result is truly remarkable; however, the incoherence of random projections is not optimal, and typically scales as $\mu\sim\sqrt{2\log(n)}$.  
Moreover, it may be difficult to obtain random projections of the full state $\bx$ in physical applications.  

There are many alternative strategies to solve for the sparsest solution to~\eqref{Eq:CS}.  
Greedy algorithms are often used~\cite{Tropp:2004,Tropp:2006,Tropp:2006b,Tropp2007ieeetit,Gilbert:2010}, including the compressed sampling matching pursuit (CoSaMP) algorithm~\cite{Needell:2010}.  
In addition, there is additional theory about how sparse the random projections may be for compressed sensing~\cite{Li:2006,Wang:2010}.  
\begin{figure*}[t]
	\centering
	\begin{overpic}[width=.75\textwidth]{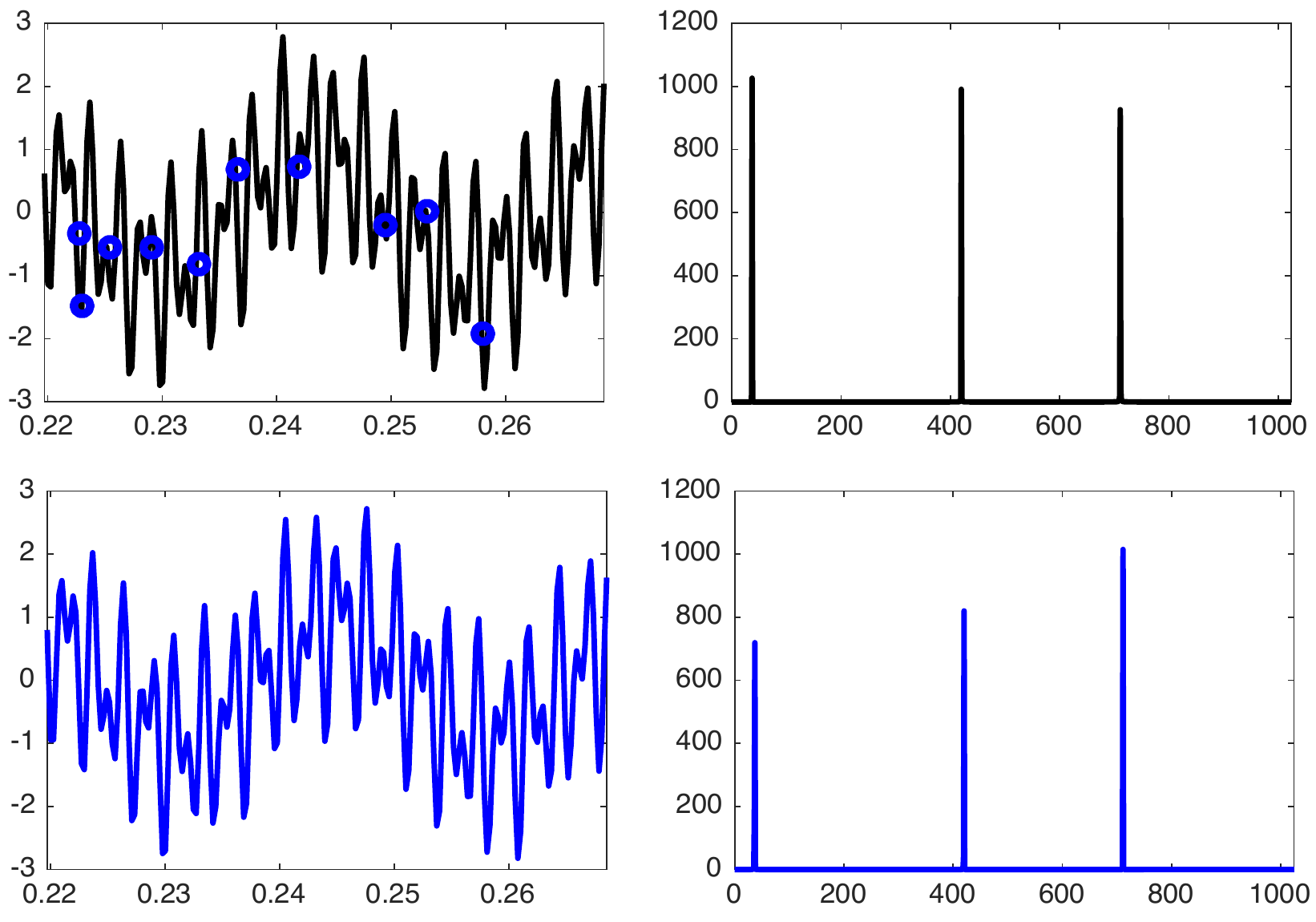}
	\small
		\put(19,-2){Time [s]}
		\put(68,-2){Frequency [Hz]}
		\put(-2,18.5){$\bx$}
		\put(-2,53.75){$\bx$}
		\put(62,30){\small Power Spectral Density}
		\put(62,65.5){\small Power Spectral Density}
	\end{overpic}
	\vspace{.15in}
	\caption{Compressed sensing, applied to three-tone signal.}\label{Fig:CSexample}
\end{figure*}

\subsection[Compressed sensing example]{Compressed sensing example}
\label{sb:cs_example}
As a simple example, we consider a sparse signal that is constructed as the sum of three distinct cosine waves,
\begin{equation}
x(t) = \cos(2\pi\times 37 t) +  \cos(2\pi\times 420 t) +  \cos(2\pi\times 711 t). 
\end{equation}
The Shannon-Nyquist sampling theorem~\cite{Nyquist1928taiee,Shannon1948bstj} states that for full signal reconstruction, we must sample at twice the highest frequency present, indicating a theoretical minimum sampling rate of $1422\,$Hz.  
However, since the signal is sparse, we may sample at considerably lower than the Nyquist rate, in this case at an average of $256\,$Hz, shown in Figure~\ref{Fig:CSexample}.
Note that for accurate sparse signal reconstruction, these measurements must be randomly spaced in time, so that the relative spacing of consecutive points may be quite close or quite far apart.  
Spacing points evenly with a sampling rate of $256\,$Hz would alias the signal, resulting in poor reconstruction. Matlab code for reproducing Figure~\ref{Fig:CSexample} is provided below.
\begin{lstlisting}[firstline=3,lastline=23]
N = 4096;
t = linspace(0, 1, N); 
x = cos(2*pi*37*t) + cos(2*pi*420*t) + cos(2*pi*711*t); 

%% Randomly sample signal
p = 256; % num. samples, p=N/16
perm = randperm(N,p);
y = x(perm); % compressed measurement

%% Solve compressed sensing problem
Psi = dct(eye(N, N));  % build Psi
Theta = Psi(perm, :);  % Measure Psi

cvx_begin; % L1-minimization with CVX
    variable s(N); 
    minimize( norm(s,1) ); 
    subject to 
        Theta*s == y';
cvx_end;

xrecon = idct(s); % reconstruct x
\end{lstlisting}

\section{Optimal sparse sensing in a tailored basis}\label{Sec:Tailored}

\begin{figure*}[t]
\begin{Sidebar}{Sidebar: Proper orthogonal decomposition and eigenfaces}{Sidebar: Proper orthogonal decomposition and eigenfaces}
	\label{Sidebar2}
	One of the most visually striking and intuitive applications of proper orthogonal decomposition (POD) is the feature extraction of facial images. These POD eigenmodes of these datasets are called {\em eigenfaces} due to their resemblance to generic human faces.  We demonstrate this application of POD on the extended Yale B dataset~\cite{YaleB2001ieee,YaleB2005ieee}, consisting of cropped and aligned images of several individuals in different lighting conditions. We obtain a resized version of the dataset in the form of Matlab data files from~\cite{yale}. Each image is a $32\times 32$ matrix of grayscale pixel values, reshaped into a column vector of length 1024 and assembled into a data matrix $\bX$.  This example, detailed in~``\nameref{Sec:Results:Eigenfaces}", is a benchmark problem for sensor selection.  
	
Matlab code for obtaining eigenfaces from training images is provided. First, training images are used to assemble a mean-subtracted data matrix.
\begin{lstlisting}[firstline=36,lastline=38]
meanface = mean(X,2);
X = X-repmat(meanface,1,size(X,2)); % mean centered data
\end{lstlisting}
Next, POD eigenfaces are obtained using the singular value decomposition of the data matrix. Outputs from both code snippets are visualized in Figure~\ref{fig:yale}.
\begin{lstlisting}[firstline=39,lastline=41]
% proper orthogonal decomposition
[Psi,S,V] = svd(X,'econ');
eigenfaces = Psi(:,1:10);
\end{lstlisting}
\begin{figure}[H]
\centering
		\begin{overpic}[width=.85\textwidth]{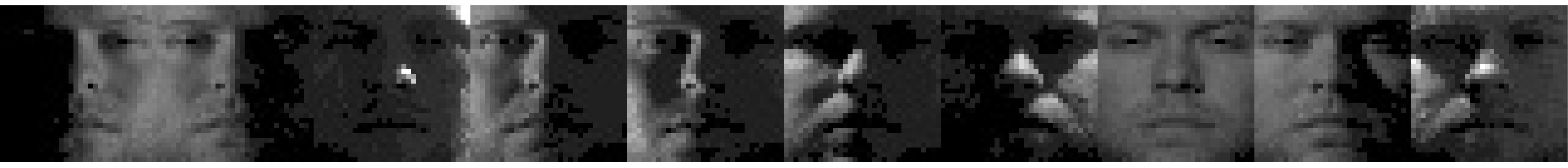}
		\put(17,8){\color{white} Training images from the Extended Yale B dataset\label{fig:yale_train}}
		\end{overpic}		
		\begin{overpic}[width=.85\textwidth]{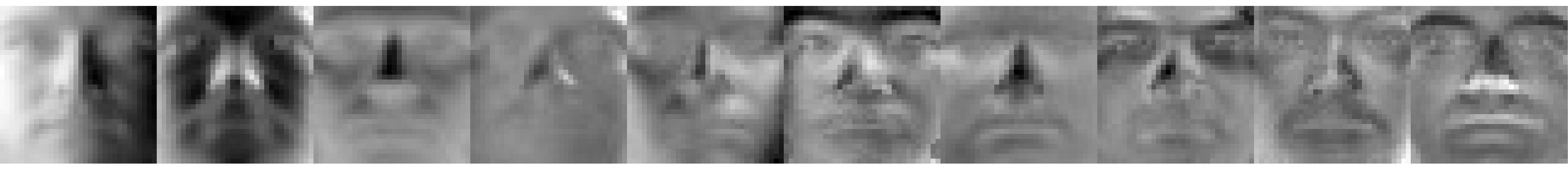}	
		\put(30,8){\color{white} First ten POD eigenfaces\label{fig:eigenfaces}}
		\end{overpic}
	\caption{Proper orthogonal decomposition (POD) modes successfully recover important facial information such as the main facial features (eyes, nose, mouth) followed by depth information (brows, ridges, chin).\label{fig:yale}}
\end{figure}

\end{Sidebar}
\end{figure*}

The compressed sensing strategy above is ideal for the recovery of a high-dimensional signal of \emph{unknown} content using random measurements in a universal basis.  
However, if information is available about the \emph{type} of signal (e.g., the signal is a turbulent velocity field or an image of a human face), it is possible to design optimized sensors that are tailored for the particular signals of interest.  
Dominant features are extracted from a training dataset consisting of representative exaples, for example using the proper orthogonal decomposition (POD).  
These low-rank features, mined from patterns in the data, facilitate the design of specialized sensors that are tailored to a specific problem.  

Low-rank embeddings, such as POD, have already been used in the ROM community to select measurements in the state space that are informative for feature space reconstruction.  
The so-called {\em empirical interpolation methods} seek the best interpolation points for a given basis of POD features.  
These methods have primarily been used to speed up the evaluation of nonlinear terms in high-dimensional, parameterized systems. However, the resulting interpolation points correspond to measurements in state space, and their use for data-driven sensor selection has largely been overlooked. We will focus on this formulation of sensor selection and explore sparse, convex, and greedy optimization methods for solving it. 

\begin{figure*}[t]
\begin{Sidebar}{Sidebar: Comparison -- Sparse sensing methods}{Sidebar: Comparison -- Sparse sensing methods}
	\label{sb:sensing_compare}
	\revision[]{
\begin{minipage}{\textwidth}
\centering
\captionof{table}{Comparison: compressed sensing and data-driven sensing} \label{tab:comparison} 	
		\noindent \begin{tabular}{|l|| l| l| } 
			\hline
			 & Compressed sensing & Data-driven QR sensing \tabularnewline \hline
			 No. samples & $\mathcal{O}(K\log\frac{n}{K})$ & $p \ge r, r\ll K$ \\
			Samples & Random & Optimized \\
			Basis & Universal & Tailored \\
			Training data & No & Yes \\			
			Solution procedure & Convex optimization & Least squares \\
			\hline
		\end{tabular}
\end{minipage}

\vspace{.5em}
	We summarize in Table~\ref{tab:comparison} the distinction between two competing perspectives for signal recovery: compressed sensing and data-driven sensing in a tailored basis. Compressed sensing recovers an unknown, underdetermined signal, by sparsity-promoting $\ell_1$ minimization. Figure~\ref{Fig:Lpnorms} further illustrates the $\ell_q$ norm constraints ($0\le q\le 1$), which promote sparsity of the solution. This approach mandates random measurements to maintain incoherence with the basis and requires $O(K\log\frac{n}{K})$ measurements, where $K$ is the signal's sparsity within this basis. In contrast, data-driven sensing requires only as many optimized samples as the data's intrinsic rank $r$, where $r$ is often much smaller than $K$. Sensor locations are thus selected to be informative based on system structure, yielding compressed measurements that streamline subsequent analysis, particularly when the original system is high-dimensional. 
Although compressed sensing can recover a wider range of signals, random sensing and convex optimization procedures may be impractical for high-dimensional, structured signals encountered in physical systems. In these cases, data-driven sensing is beneficial since it permits a drastic reduction in the required number of sensors and downstream computation.}{}

\begin{figure}[H]
\centering
\begin{overpic}[width=.425\textwidth]{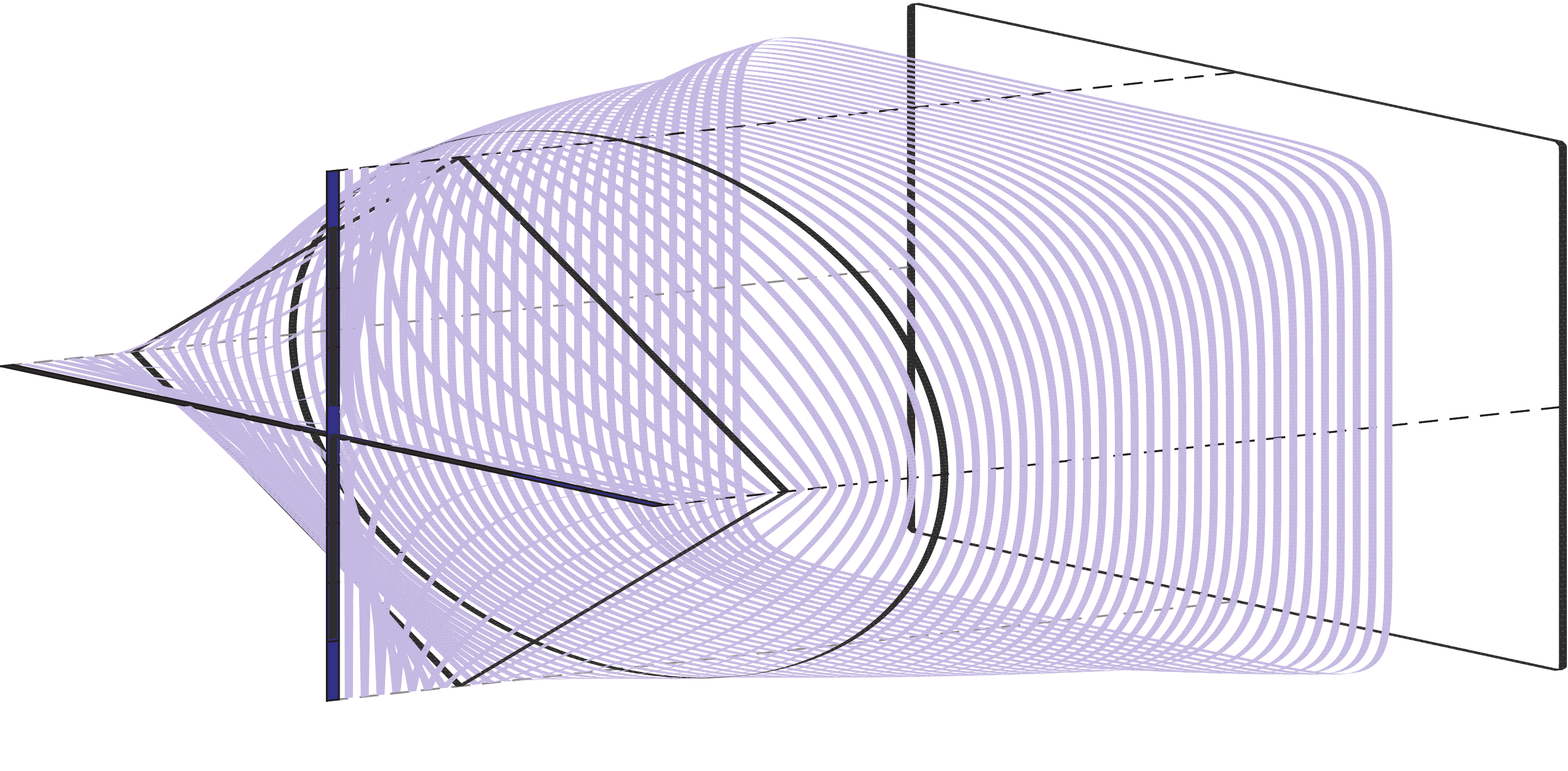}
		\put(-3,26){$a_1$}
		\put(20,40){$a_2$}
		\put(18,0){$\ell_0$}		
		\put(27,1){$\ell_1$}
		\put(41,2){$\ell_2$}
		\put(93,3){$\ell_\infty$}
	\end{overpic}
	\caption{Progression of unit balls $\|\ba\|_q=1$ of $\ell_q$ norm where $q\in[0,\infty]$. Unit balls of $\ell_0,\ell_1,\ell_2,\mbox{ and }\ell_\infty$ norms are outlined in black. Minimizing $\ell_q$ norms when $q\in[0,1]$ promotes sparsity since these unit balls are oriented toward the coordinate axes. }\label{Fig:Lpnorms}
\end{figure}
\end{Sidebar}
\end{figure*}

We begin with brief expositions on POD and our mathematical formulation of sensor placement, followed by an overview of related work in design of experiments and sparse sampling. We conclude this section with our generalized sensor selection method that connects empirical interpolation methods, such as QR pivoting to optimize condition number, with D-optimal experimental design~\cite{doptimal}. 
The QR pivoting method described in~``\nameref{Sec:QR}" is particularly favorable, as it is fast, simple to implement, and provides nearly optimal sensors tailored to a data-driven POD basis. \revision[]{Finally, the distinctions between compressed sensing and our data-driven sensing are summarized in~``\nameref{sb:sensing_compare}".}{}

\subsection{Proper orthogonal decomposition}

POD is a widespread data-driven dimensionality reduction technique~\cite{berkooz1993proper,HLBR_turb} used in many domains; it is also commonly known as the Karhunen-Lo\`eve expansion, principal component analysis (PCA)~\cite{Pearson:1901}, and empirical orthogonal functions~\cite{lorenzMITTR56}. 
POD expresses high-dimensional states $\bx\in\reals^n$ as linear combinations of a small number of orthonormal eigenmodes $\bpsi$ (i.e., POD modes) that define a low-dimensional embedding space. States are projected into this POD subspace, yielding a reduced representation that can be used to streamline tasks that would normally be expensive in the high-dimensional state space. 
This low-rank embedding does not come for free, but instead requires training data to \emph{tailor} the POD basis to a specific problem.  
POD is illustrated on a simple example of extracting coherent features in images of human faces in~``\nameref{Sidebar2}".

A low-dimensional representation of $\bx$ in terms of POD coefficients $\mathbf{a}$ can be lifted back to the full state with a linear combination of POD modes,
\begin{equation*}
\bx_i \approx \sum_{k=1}^r a_k(t_i)\bpsi_k(x).
\end{equation*}
For time-series data $\bx_i$, the coefficients $a_k(t_i)$ vary in time and $\bpsi_k(x)$ are purely spatial modes without time dependence, resulting in a space-time separation of variables.  
Thus, care should be taken applying POD to data from a traveling wave problem.   

The eigenmodes $\bpsi_k$ and POD coefficients $a_k$ are easily obtained from the singular value decomposition (SVD). Given a data matrix of state space observations $\bX=[\bx_1 ~\bx_2~\dots~\bx_m]$, the resulting eigenmodes are the orthonormal left singular vectors $\bPsi$ of $\bX$ obtained via the SVD,
\begin{equation}
\bX = \mathbf{\Psi\Sigma V}^T\approx \bPsi_r\bSigma_r\bV_r^T.
\end{equation}
The matrices $\bPsi_r$ and $\bV_r$ contain the first $r$ columns of $\bPsi$ and $\bV$ (left and right singular vectors, respectively), and the diagonal matrix $\bSigma_r$ contains the first $r\times r$ block of $\bSigma$ (singular values).  
 The SVD is the optimal least squares approximation to the data for a given rank $r$, as demonstrated by the Eckart-Young theorem~\cite{Eckart1936psych}
\begin{equation}
\bX_\star = \argmin_{\tilde{\bX}}\|\bX - \tilde{\bX}\|_F ~\text{ s. t. }~ \text{rank}(\tilde{\bX})=r, 
\end{equation}
where $\bX_\star = \bPsi_r\bSigma_r\bV_r^T$, and $\|\cdot\|_F$ is the Frobenius norm. 
The low-dimensional vector of POD coefficients for a state $\bx$ is given by the orthogonal projection ${\ba = \bPsi_r^T\bx}$. 
Thus, the POD is a widely used dimensionality reduction technique for high-dimensional systems. This reduction allows computational speedup of numerical time-stepping, parameter estimation, and control. 

Choosing the intrinsic target rank without magnifying noise in the data is a difficult task. In practice, $r$ is often chosen by thresholding the singular values to capture some percentage of the variance in the data. An optimal hard threshold is derived in~\cite{Gavish2014ieeetit} based on the singular value distribution and aspect ratio of the data matrix, assuming additive Gaussian white noise of unknown variance.  This threshold criterion has been effective in practice, even in cases where the noise is likely not Gaussian.

\begin{figure*}[t]
\begin{Sidebar}{Sidebar: Condition number}{Sidebar: Condition number}
	\label{Sidebar3}
		The condition number of a matrix $\mathbf{\Theta}$ gives a measure of how sensitive matrix multiplication or inversion is to errors in the input, with larger condition number indicating worse performance.  
		The condition number $\kappa(\mathbf{\Theta})$ is the ratio of the maximum and minimum singular values of $\mathbf{\Theta}$
		\begin{equation}
		\kappa(\mathbf{\Theta}) = \frac{\sigma_{\max}(\mathbf{\Theta})}{\sigma_{\min}(\mathbf{\Theta})}.
		\end{equation}
		To see the effect of the condition number on matrix multiplication, consider a square, invertible $\bTheta$ and an input signal $\mathbf{x}$ that is contaminated by noise $\mathbf{\epsilon_x}$.  
		Further, we consider the \emph{worst-case scenario} where $\mathbf{x}$ is aligned with the right singular vector of $\mathbf{\Theta}$ corresponding to the minimum singular value $\sigma_{\min}(\mathbf{\Theta})$ and where the error $\mathbf{\epsilon_x}$ is aligned with the right singular vector of $\mathbf{\Theta}$ corresponding to the maximum singular value $\sigma_{\max}(\mathbf{\Theta})$.  Thus, error is scaled by $\sigma_{\max}$ while the signal is scaled by $\sigma_{\min}$
		\begin{equation}
		\mathbf{\Theta}\left(\mathbf{x} + \mathbf{\epsilon_x}\right) = \sigma_{\min}\mathbf{x} + \sigma_{\max}\mathbf{\epsilon_x}.  
		\end{equation}
		Thus, we see that the signal to noise ratio (SNR) is reduced by a factor equal to the condition number,
		\begin{equation}
		\text{SNR}_{\text{in}} = \frac{\mathbf{x}}{\mathbf{\epsilon_x}} \quad\Longrightarrow\quad \text{SNR}_{\text{out}} = \frac{\sigma_{\min}}{\sigma_{\max}}\frac{\mathbf{x}}{\mathbf{\epsilon_x}} = {\text{SNR}_{\text{in}}}/{\kappa}.
		\end{equation}
		If the condition number is large, then error can be amplified relative to the signal.  
		The ideal condition number is $1$, where $\mathbf{\Theta}$ has all singular values equal to $1$, for example if $\mathbf{\Theta}$ is a unitary matrix. 
		
		Similarly, the sensitivity of the matrix inverse to error is also related to the condition number
		\begin{equation}
		\mathbf{x}+\mathbf{\epsilon_x} = \mathbf{\Theta}^{-1}\left(\mathbf{y} + \mathbf{\epsilon_y}\right).
		\end{equation}
		Again, worst-case scenario errors in the inversion are amplified by the condition number (this time with the error being aligned with the singular vector corresponding to minimum singular value).  
		
		In general, a random error $\mathbf{\epsilon_x}$ or $\mathbf{\epsilon_y}$ will have some component in this \emph{worst-case} direction, and will be amplified by the maximum singular value of $\mathbf{\Theta}$ or $\mathbf{\Theta}^{-1}$.  
		Thus, it is desirable to explicitly control the condition number of $\bTheta=\bC\bPsi_r$ by choice of the row selection operator $\bC$.  
		For invertible matrices, the condition numbers of $\bTheta$ and $\bTheta^{-1}$ are the same.  The discussion above generalizes for rectangular $\bTheta$.    
\end{Sidebar}
\end{figure*}

\subsection{Sensor placement for reconstruction}

We optimize sensor placement specifically to reconstruct high-dimensional states from point measurements, given data-driven or tailored bases. Recall that full states may be expressed as an unknown linear combination of basis vectors
\begin{equation}
\label{eqn:state_i} \color{black}
x_j = \sum_{k=1}^r \Psi_{jk} a_k,
\end{equation}
where $\Psi_{jk}$ is the coordinate form of $\bPsi_r$ from~\eqref{eqn:state}.  
Effective sensor placement results in a point measurement matrix $\bC$ that is optimized to recover the modal mixture $\ba$ from sensor outputs $\by$.  Point measurements require that the sampling matrix $\bC\in\reals^{p\times n}$ be structured in the following way
\begin{equation}
\bC = \begin{bmatrix} 
\be_{\ind_1} &
\be_{\ind_2} &
\dots &
\be_{\ind_p}
\end{bmatrix}^T,
\end{equation}
where $\be_j$ are the canonical basis vectors for $\reals^n$ with a unit entry at index $j$ and zeros elsewhere.  Note that point measurements are fundamentally different than the suggested random projections of compressive sensing.  
The measurement matrix results in the linear system
\begin{equation}
\label{eqn:observations_i} \color{black}
y_i = \sum_{j=1}^n C_{ij} x_j = \sum_{j=1}^n C_{ij} \sum_{k=1}^r \Psi_{jk} a_k,
\end{equation}
where $C_{ij}$ is the coordinate form of $\bC$ from~\eqref{eqn:observations}.
The observations in $\by$ consist of $p$ elements selected from $\bx$
\begin{equation}
\by = \bC\bx = [x_{\ind_1} ~x_{\ind_2}~ \dots~ x_{\ind_p}]^T,
\end{equation}
where $\Ind = \{\ind_1,\dots,\ind_p\} \subset \{1,\dots,n\}$ denotes the index set of sensor locations with cardinality $|\Ind|=p$. 

\begin{figure}
	\centering
	\includegraphics[width=.5\textwidth]{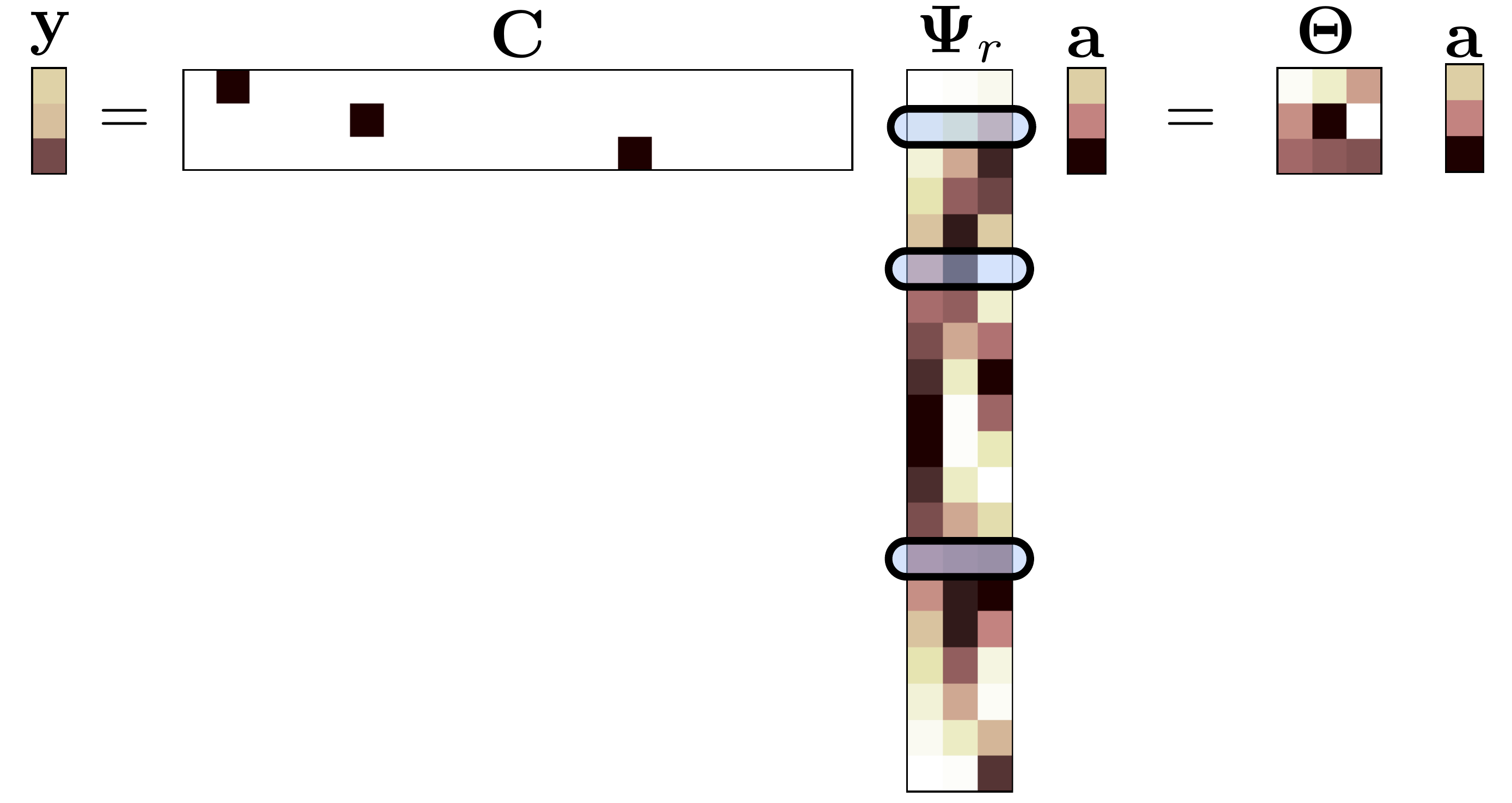}
	\vspace{-.1in}
	\caption{Full state reconstruction of $\bx$ from point  observations ($\by$) is accomplished using least squares estimation of POD coefficients ($\ba=\bTheta^\dagger\by$). }\label{Fig:QRscheme}
\end{figure}

When $\bx$ is unknown, it can be reconstructed by approximating the unknown basis coefficients $\ba$ with the Moore-Penrose pseudoinverse, $\ba= \bTheta^{\dagger}\by=(\bC\bPsi_r)^\dagger \by$.  Equivalently, the reconstruction is obtained using
\begin{equation}
	\hat{\bx} = \bPsi_r\hat{\ba}, \mbox{ where } 
	\hat{\ba} = \begin{cases} \bTheta^{-1}\by =  (\bC\bPsi_r)^{-1}\by, & p=r, \\
	\bTheta^{\dagger}\by = (\bC\bPsi_r)^{\dagger}\by, & p>r
	\end{cases}.
\end{equation}
A schematic of sparse sampling in a tailored basis $\bPsi_r$ is shown in Fig.~\ref{Fig:QRscheme}.  
The optimal sensor locations are those that permit the best possible reconstruction $\hat\bx$. 
Thus, the sensor placement problem seeks rows of $\bPsi_r$, corresponding to point sensor locations in state space, that optimally condition inversion of the matrix $\bTheta$. 
For brevity in the following discussion we denote the matrix to be inverted by $\bM_{\Ind}=\bTheta^T\bTheta$ ($\bM_{\Ind}=\bTheta$ if $p=r$). 
Recall that $\Ind$ determines the structure of $\bC$, i.e. the sensor locations, and hence affects the condition numbers of $\bTheta$ and $\bM_{\Ind}$.  
The condition number of the system may be indirectly bounded by optimizing the spectral content of $\bM_{\Ind}$ using its determinant, trace, or spectral radius. 
For example, the spectral radius criterion for $\bM_{\Ind}^{-1}$ maximizes the smallest singular value of $\bM_{\Ind}$
\begin{equation}
\label{eqn:e_opt}
	\Ind_{\star} = \argmin_{\Ind,|\Ind|=p} \|\bM_{\Ind}^{-1}\|_2 = \argmax_{\Ind,|\Ind|=p} \sigma_{\min}(\bM_{\Ind}) .
\end{equation}
Likewise, the sum (trace) or product of magnitudes (determinant) of its eigenvalue or singular value spectrum may be optimized
\begin{equation}
\label{eqn:a_opt}
	\Ind_{\star} = \argmax_{\Ind,|\Ind|=p} \mbox{tr}(\bM_{\Ind,|\Ind|=p}) = \argmax_{\Ind} \sum_{i} \lambda_i(\bM_{\Ind}), 
\end{equation}
\begin{align}
\label{eqn:d_opt}
	\Ind_{\star} = \argmax_{\Ind,|\Ind|=p} |\det \bM_{\Ind}| &= \argmax_{\Ind,|\Ind|=p} \prod_{i} \left|\lambda_i(\bM_{\Ind})\right| \nonumber \\ &=\argmax_{\Ind,|\Ind|=p} \prod_i \sigma_i(\bM_{\Ind}).
\end{align}
Direct optimization of the above criteria requires a combinatorial search over $n\choose p$ possible sensor configurations and is hence computationally intractable even for moderate $n$. Several heuristic greedy sampling methods have emerged for state reconstruction specifically with POD bases. These {\em gappy POD}~\cite{Everson1995gappy} methods originally relied on random sub-sampling.  However, significant performance advances where demonstrated by using principled sampling strategies for reduced-order models (ROMs) in fluid dynamics~\cite{Willcox2006compfl}, ocean modeling~\cite{Yildirim2009oceanmod} and aerodynamics~\cite{Chen2011aiaa}.  More recently, variants of the so-called {\em empirical interpolation method} (EIM, DEIM and Q-DEIM)~\cite{Barrault2004crm,Chaturantabut2010siamjsc,drmac2016siam} have provided near optimal sampling for interpolative reconstruction of nonlinear terms in ROMs. This work examines an approximate greedy solution given by the matrix QR factorization with column pivoting of $\bPsi_r^T$, which builds upon the Q-DEIM method~\cite{drmac2016siam}.

\revision[] 
{
	
}{ 
The goal of sparse measurement selection is to choose a measurement matrix $\bC$ representing $p$ distinct point measurements in state space, so that $\bC$ consists of $p$ rows of the identity matrix.  
Although more general linear measurements of $\bx$ may be admissible in some problems, point measurements are physically appealing as spatially localized sensors.  
The row-wise sum $\mathbf{c}^+=\sum_j \mathbf{c}_k$ satisfies $\| \mathbf{c}^+ \|_0 = p$,
where $\mathbf{c}_k$ is the $k$-th row of $\bC$. 
The sensor selection problem can be formulated as choosing $\bC$ to make the pseudoinverse of $\bTheta=\bC\bPsi_r$ as well-conditioned as possible, thus making the estimation of the coefficients in $\ba$ robust to measurement noise on $\by$.  The condition number is discussed in more detail in~``\nameref{Sidebar3}".  

Condition number minimization is the viewpoint taken by the reduced-order modeling community for selecting point measurements within a low-rank POD basis. These empirical interpolation methods (EIM)~\cite{Barrault2004crm} assume $p=r$ sensors so $\bTheta$ is square and invertible. EIM and subsequent discrete variants such as DEIM~\cite{Chaturantabut2010siamjsc} and Q-DEIM~\cite{drmac2016siam} are greedy procedures that control the condition number of $(\bC\bPsi_r)^{-1}$ by minimizing the spectral norm (i.e. largest singular value)
\begin{equation}
\label{eqn:DEIM_opt}
\bC_\star = \argmin_{\bC} \|(\bC\bPsi_r)^{-1}\|_2,
\end{equation}
where $\bC$ is subject to the same structural constraints as above. 
Related strategies include selecting measurements at maxima of successive POD modes or iteratively seeking measurements that decrease the condition number of $\bC\bPsi_r$~\cite{Yildirim2009oceanmod,Willcox2006compfl}.

This strategy is particularly advantageous when only a few modes are required to characterize the data. When singular values decay slowly, we may require $p>r$ measurements for well-conditioned reconstruction. In~``\nameref{Sec:QR}", we generalize one particular empirical interpolation method, Q-DEIM, which uses QR pivoting to determine sensor locations.  In particular, we extend this method to the case of $p> r$ sensors, making it viable for more general sensor selection and significantly speeding up its computation.

\subsubsection{General formulation}
Sensor selection is more generally framed as follows: given a set of $n$ possible measurement indices $V$, select a subset of sensor indices $S$ to optimize a carefully chosen function evaluating their quality. The problem is mathematically realized as 
\begin{equation}
\max_{S\subseteq V} f(S) \mbox{ subject to constraints on }S.
\end{equation}
This problem belongs to an area of research called submodular function optimization, surveyed in~\cite{Krause2012submodular}. The so-called submodular function, ${f(S):2^V\mapsto\reals}$, evaluates these subsets to some scalar value, thus rendering a brute force search over its entire domain (the power set of all possible subsets $2^V$) computationally intractable even for moderate values of $n$. 

Possible choices of $f$ include mutual information, entropy, estimation error covariance, and spatial coverage of sensors. Often, even a single evaluation of $f$ is expensive when matrix determinants or factorizations are involved. The present work aims to reconstruct states from measurements with minimal variance of reconstruction error, which requires maximizing the determinant of error covariance matrices. The minimization of various error metrics is an active topic in optimal experiment design, for instance, the determinant or volume (D-optimal), the trace or mean-squared error (A-optimal), and maximum eigenvalue or worst-case variance (e-optimal) of error covariance matrices. 

Joshi and Boyd~\cite{Joshi2009ieee} frame sensor selection as selecting $p$ rows of $\bPsi_r$ that minimize the volume of the estimation error covariance matrix, $(\bC\bPsi_r)^T\bC\bPsi_r=\bTheta^T\bTheta$. Since matrix volume is the absolute value of its determinant, the naive approach requires evaluating $n\choose p$ determinants over a combinatorially large set of row-selected submatrices. 
Instead, the authors approach this search with a computationally tractable convex optimization,
\begin{align}
\label{eqn:sensor_opt_l1}
\bC_1 = &\argmax_\bC \logdet (\bC\bPsi_r)^T\bC\bPsi_r \nonumber \\
&\subjto \|\bc^+\|_0=p,~  c_i^+ \in [0,1].
\end{align}
Each intermediate iteration of this method operates on an $n\times n$ matrix, so that storage requirements scale as $\mathcal{O}(n^2)$. Thus this approach is costly for high-dimensional states. The authors frame $\bPsi_r$ as a general matrix whose rows are possible measurements from which they select a subset of measurements. However, this matrix is not explicitly restricted to a POD basis, as proposed here.
} 

\subsection{Sparse sensor placement with QR pivoting}\label{Sec:QR}

An original contribution of this work is extending Q-DEIM to the case of oversampled sensor placement, where the number of sensors exceeds the number of modes used in reconstruction ($p>r$). The key computational idea enabling oversampling is the QR factorization with column pivoting applied to the POD basis. QR pivoting itself dates back to the 1960s by Businger and Golub to solve least squares problems~\cite{Businger1965nm}, and it has found utility in various measurement selection applications~\cite{Sommariva2009qr,Heck1998qr,Seshadri2016qr}. Similar to empirical interpolation methods such as DEIM, pivots from the QR factorization optimally condition the measurement or row selected POD basis, as described below.  

	The reduced matrix QR factorization with column pivoting decomposes a matrix $\bA\in\reals^{m\times n}$ into a unitary matrix $\bQ$, an upper-triangular matrix $\bR$ and a column permutation matrix $\bC$ such that $\bA\bC^T = \bQ\bR$. The pivoting procedure provides an approximate greedy solution method for the optimization in~\eqref{eqn:d_opt}, which is also known as submatrix volume maximization because matrix volume is the absolute value of the determinant. QR column pivoting increments the volume of the submatrix constructed from the pivoted columns by selecting a new pivot column with maximal 2-norm, then subtracting from every other column its orthogonal projection onto the pivot column (see Algorithm \ref{alg:qrpivot}). Pivoting expands the submatrix volume by enforcing a diagonal dominance structure~\cite{drmac2016siam}
	\begin{equation}
		\sigma_i^2 = |r_{ii}|^2 \ge \sum_{j=i}^k |r_{jk}|^2; \quad 1\le i \le k \le m.
	\end{equation}
	This works because matrix volume is also the product of diagonal entries $r_{ii}$
	\begin{equation}
		\left| \det \bA \right| = \prod_i \sigma_i = \prod_i |r_{ii}|.
	\end{equation}
	Furthermore, the oversampled case $p>r$ may be solved using the pivoted QR factorization of $\bPsi_r\bPsi_r^T$, where the column pivots are selected from $n$ candidate state space locations based on the observation that
	\begin{equation}
		\det \bTheta^T\bTheta  = \prod_{i=1}^r \sigma_i(\bTheta\bTheta^T),
	\end{equation}
	where we drop the absolute value since the determinant of $\bTheta\bTheta^T$ is nonnegative.

	\begin{algorithm}
		\caption{Greedy sensor selection using a given tailored basis $\bPsi_r$ and number of sensors $p$}\label{alg:greedy}
		\begin{algorithmic}[1]
			\If{$p==r$}
			\State $\Ind \gets \mbox{qrPivot}(\bPsi_r,r)$
			\ElsIf{$p>r$}
			\State $\Ind \gets \mbox{qrPivot}(\bPsi_r\bPsi_r^T,p)$
			\EndIf				
			\State $\bC \gets \begin{bmatrix} 
			\be_{\ind_1} &
			\be_{\ind_2} &
			\dots &
			\be_{\ind_p}
			\end{bmatrix}^T$
		\end{algorithmic}
	\end{algorithm}

\begin{figure}
\begin{Sidebar}{Sidebar: QR pivoting code}{Sidebar: QR pivoting code}
\label{sb:qr_code}
The QR pivoting procedure of Algorithm~\ref{alg:greedy} can be concisely implemented using the given Matlab code. The pivoted QR factorization (Algorithm~\ref{alg:qrpivot}) is called within this code using the native Matlab \verb|qr()| subroutine. The variable \verb|Psi_r| stores the given data-driven tailored basis $\bPsi_r$, and the desired list of sensor indices (QR pivots) is returned within the \verb|pivot| variable. 
\begin{lstlisting}
if (p==r)
	% QR sensor selection, p=r (Q-DEIM)
	[Q,R,pivot] = qr(Psi_r','vector');
elseif (p>r)
	% Oversampled QR sensors, p>r
	[Q,R,pivot] = qr(Psi_r*Psi_r','vector');
end
pivot = pivot(1:p)
\end{lstlisting}
\end{Sidebar}
\end{figure}

	\begin{algorithm*}
		\caption{QR factorization with column pivoting of $\bB\in\reals^{n\times m}$}\label{alg:qrpivot}
		\begin{algorithmic}[1]		
			\Procedure{qrPivot}{ $\bB,~p$ }
			
			\State $\Ind \gets [~~]$
			\For{$k=1,\dots,p$}
			\State $\ind_k = \argmax_{j\notin \Ind} \|\mathbf{b}_j\|_2 $
			\State Find Householder $\tilde{\bQ}$ such that $\tilde{\bQ} \cdot \begin{bmatrix} b_{kk} \\ ~\\ \vdots \\ b_{nk} \end{bmatrix} = \begin{bmatrix}
			\star \\ 0 \\ \vdots \\ 0
			\end{bmatrix}.$  \hfill \mbox{\color{purple}$\star$'s are the diagonal entries of $\bR$}
			\State $ \bB \gets \mbox{diag}(I_{k-1},\tilde{\bQ}) \cdot \bB$ \hfill \mbox{
				\color{purple}
				remove from all columns the orthogonal projection onto $\mathbf{b}_{\ind_k}$}
			\State $\Ind \gets [\Ind,~\ind_k]$
			\EndFor
			\Return $\Ind$
			\EndProcedure
		\end{algorithmic}
	\end{algorithm*}
	
Thus the QR factorization with column pivoting yields $r$ point sensors (pivots) that best sample the $r$ basis modes $\bPsi_r$
\begin{equation}
\bPsi_r^T\bC^T = \bQ\bR.
\end{equation}
Based on the same principle of pivoted QR, which controls the condition number by minimizing the matrix volume, the oversampled case is handled by the pivoted QR factorization of $\bPsi_r\bPsi_r^T$,
\begin{equation}
(\bPsi_r\bPsi_r^T)\bC^T = \bQ\bR.
\end{equation}
\revision[]{Algorithm~\ref{alg:qrpivot}, the QR factorization, is natively implemented in most scientific computing software packages. A Matlab code implementation of Algorithm~\ref{alg:greedy} is provided in~``\nameref{sb:qr_code}".}{}
The oversampled case requires an expensive QR factorization of an $n\times n$ matrix, whose storage requirements scale quadratically with state dimension. 
However, this operation may be advantageous for several reasons. 
The row selection given by the first $p$ QR pivots increase the leading $r$ singular values of $\bTheta\bTheta^T$, hence increasing $\det\bTheta^T\bTheta$.  
This is the same maximization objective used in D-optimal experiment design~\cite{doptimal}, which is typically solved with Newton iterations using a convex relaxation of the subset selection objective. 
These methods require one matrix factorization of an $n\times n$ matrix per iteration, leading to a runtime cost per iteration of at least $O(n^3)$. The entire procedure must be recomputed for each new choice of $p$. 
Our proposed method only requires a single $O(n^3)$ QR factorization and results in a hierarchical list of all $n$ total pivots, with the first $p$ pivots optimized for reconstruction in $\bPsi_r$ for any $p>r$.  
Thus, additional sensors may be leveraged if available.  

\begin{figure*}
\begin{Sidebar}{Sidebar: QR pivoting fails to exploit sparsity}{Sidebar: QR pivoting fails to exploit sparsity}
		\label{Sidebar_6}
		The QR pivoting algorithm can be used to efficiently find a sparse solution of an underdetermined linear system $\by = \bTheta\bs$.  In fact, in Matlab, this is as simple as using the `$\backslash$' command: 
		\begin{lstlisting}
>>s = Theta\y;
		\end{lstlisting}			
However, QR pivoting does not always provide the \emph{sparsest} solution.  In addition, the pivot locations, and hence the nonzero entries in $\bs$, are determined entirely by $\bTheta$.  
Therefore, the structure of the nonzero coefficients in $\bs$ have nothing to do with the specific signal or measurements.  
Given $p>K$ measurements of an unknown signal in a universal basis, the QR pivot algorithm will fail in two ways:  1) it will return a sparse vector with $p$ nonzero elements (instead of the desired $K$ nonzero elements), and 2) the nonzero elements will not depend on the frequency content of the actual signal.  This observation is illustrated on two test images in Figure~\ref{fig:backslash}.  Notice that the sparse coefficients for both the mountain and cappuccino are the same, resulting in poor image reconstruction.  	
\begin{figure}[H]
		\begin{center}
		\begin{subfigure}[b]{\textwidth}
			\begin{tabular}{  c c }
				 {\bf Mountain} & {\bf Cappuccino } \\
				\begin{sideways} ~~~\quad\bf True \end{sideways}		
				\begin{overpic}[width=0.12\textwidth]{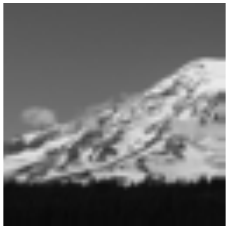}  
				\end{overpic}
				\begin{overpic}[width=0.3\textwidth]{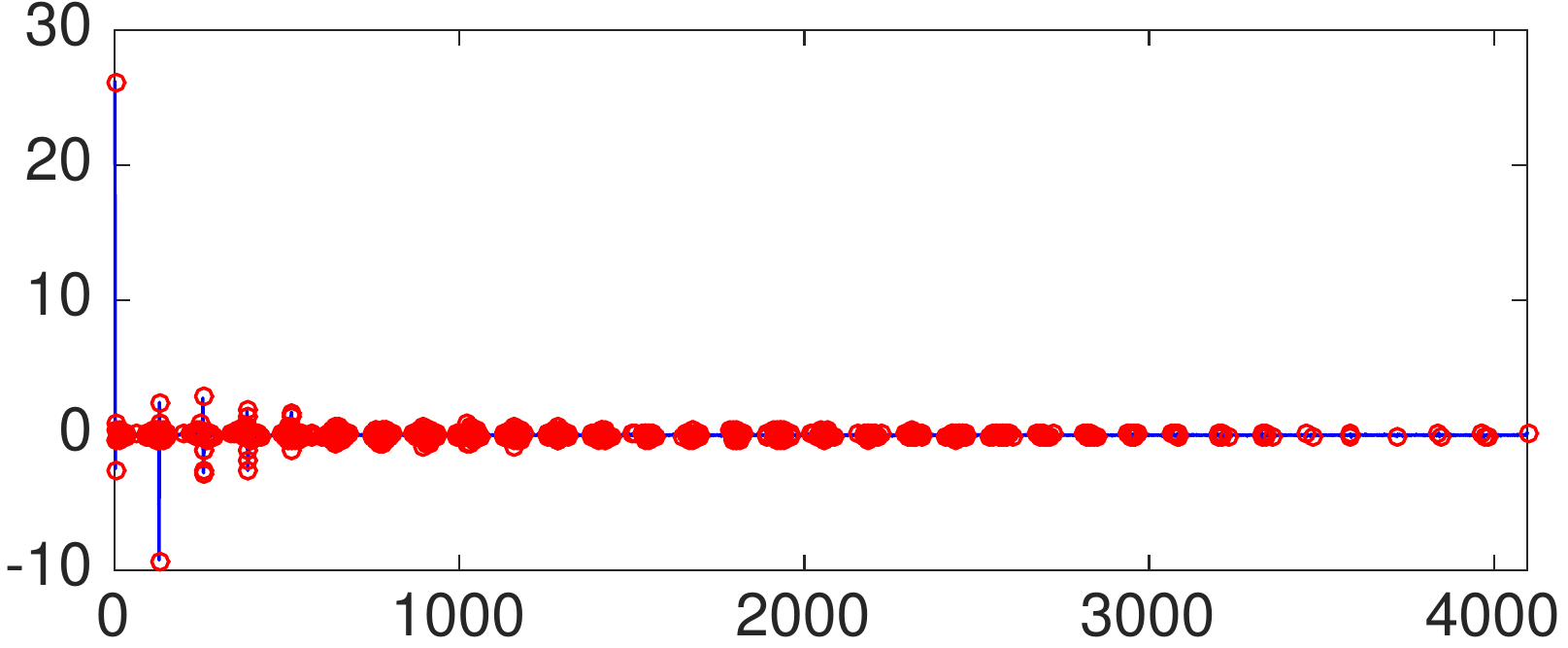} 
				\put(85,30){$\bs$}
				\end{overpic} & \qquad
				\includegraphics[width=0.12\textwidth]{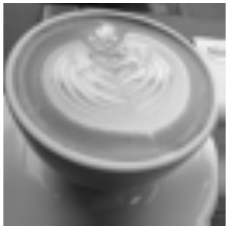}  
				\begin{overpic}[width=0.3\textwidth]{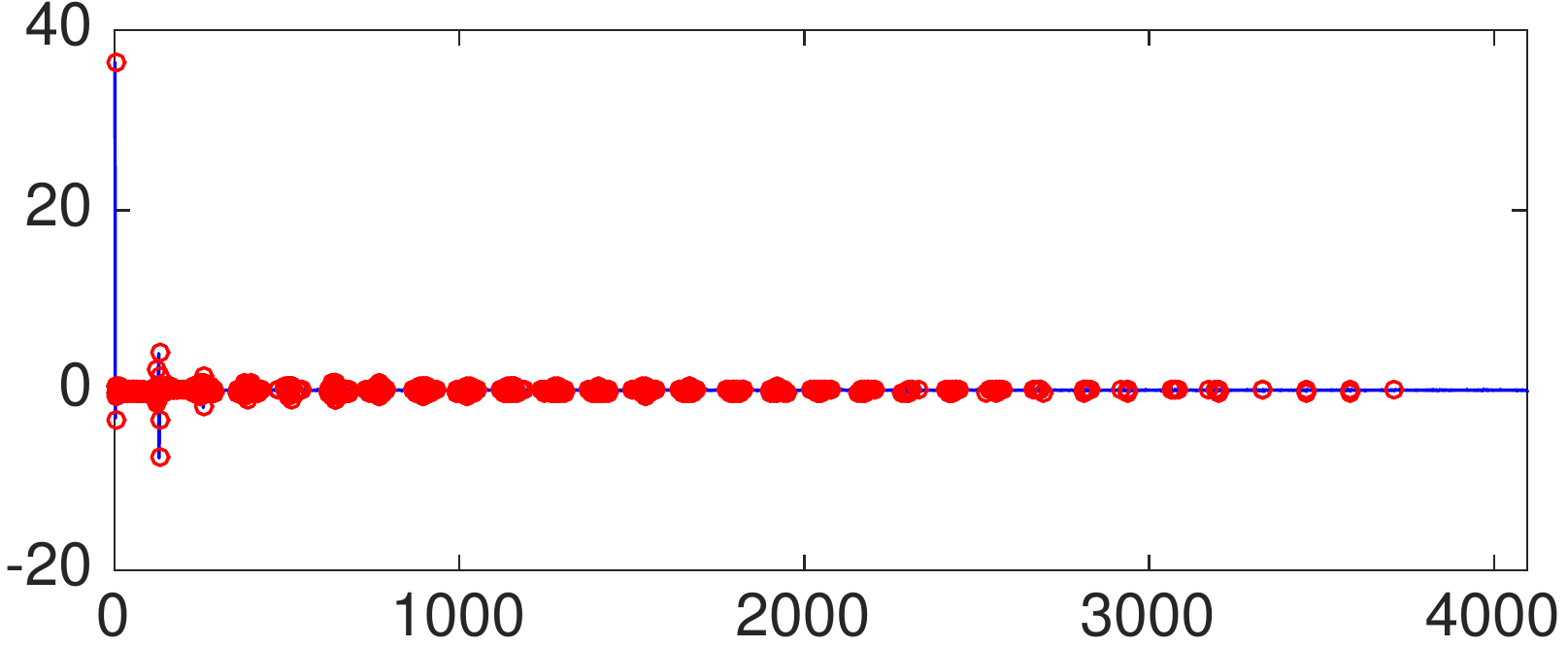}
					\put(85,30){$\bs$}
				\end{overpic} \\			
			
				\begin{sideways} \bf Compressed \end{sideways}	
				\begin{overpic}[width=0.12\textwidth]{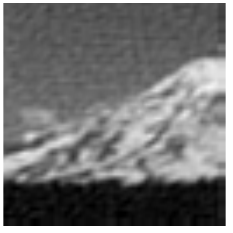}  
				\end{overpic}
				\begin{overpic}[width=0.3\textwidth]{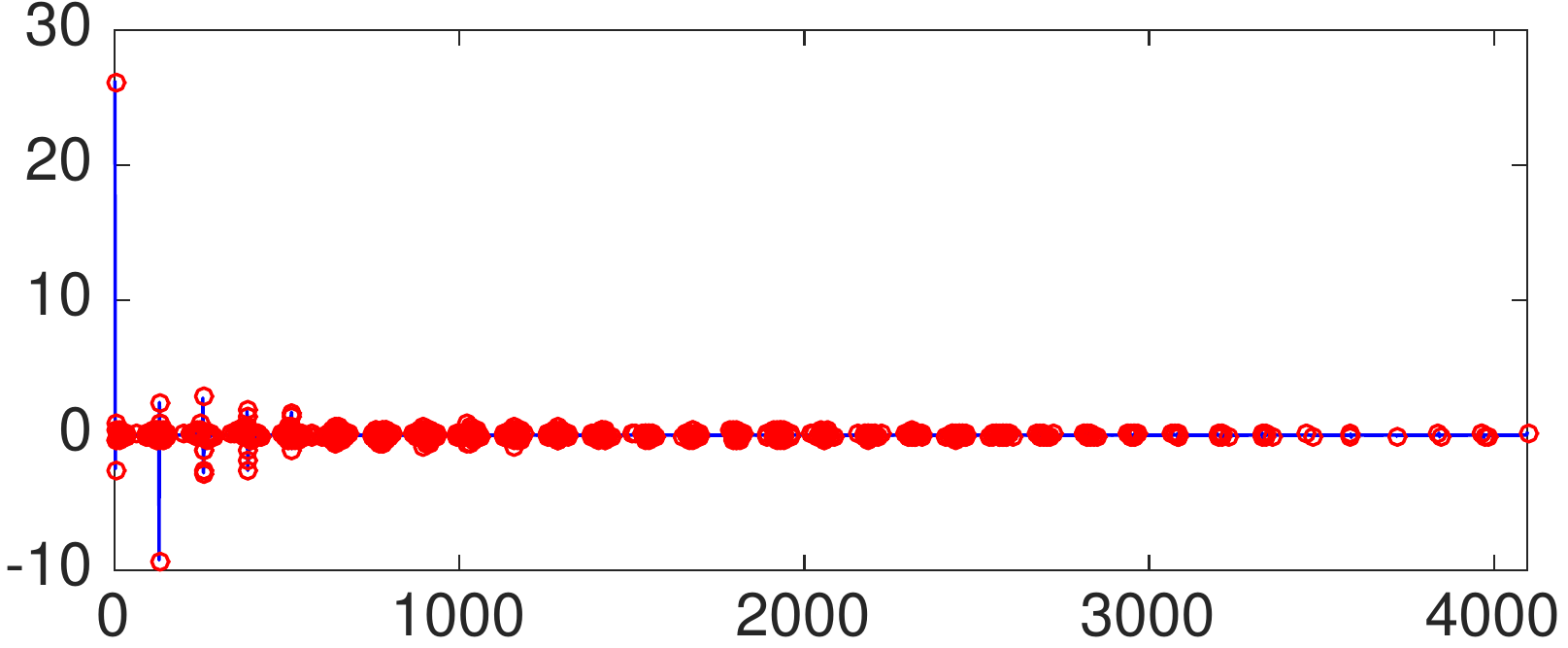}
					\put(75,30){$\bs_{\text{thresh}}$} 
				\end{overpic} & \qquad
				\includegraphics[width=0.12\textwidth]{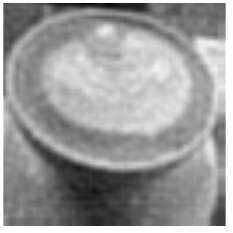}  
				\begin{overpic}[width=0.3\textwidth]{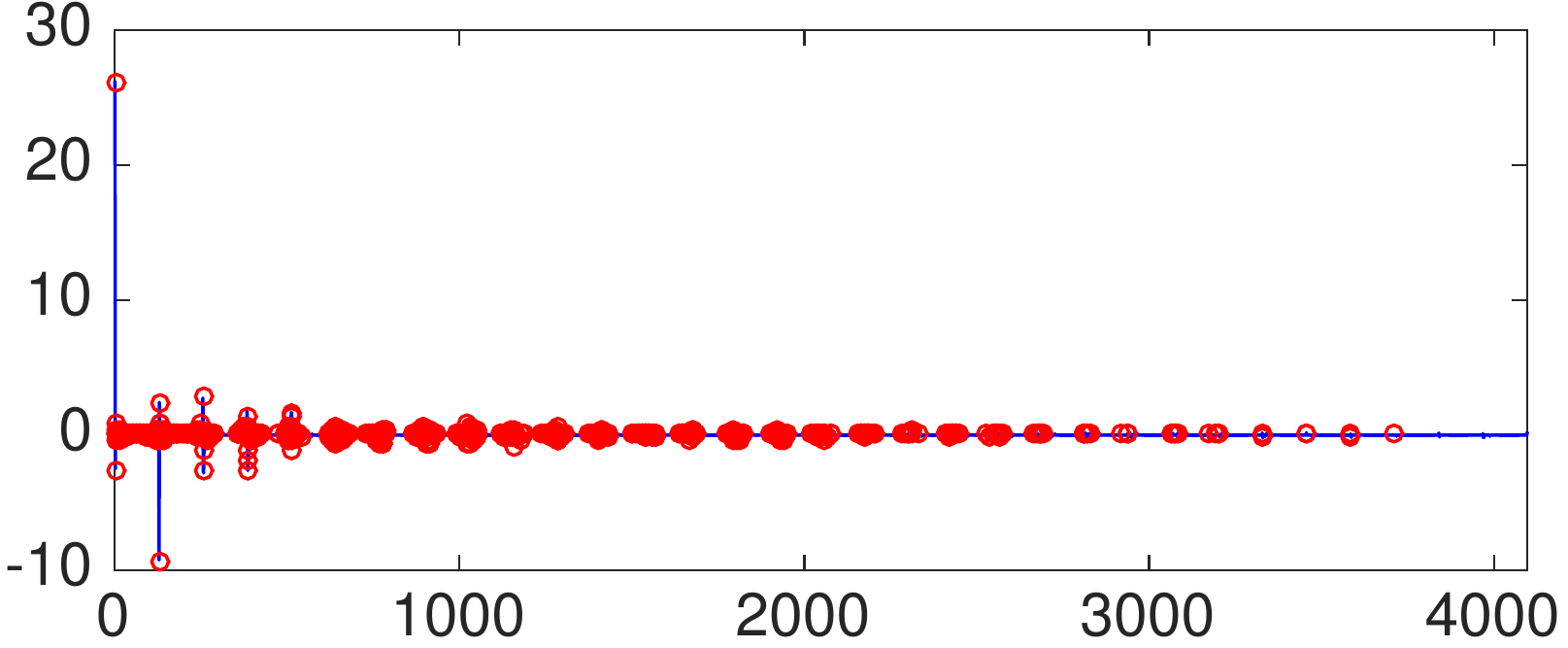}
					\put(75,30){$\bs_{\text{thresh}}$} 
				\end{overpic} \\ 	
				\end{tabular}
		\caption{{\bf Fourier coefficients}. The compressed images are recovered from a very small percentage of nonzero Fourier coefficients (red) that represent the spatial frequency signature of the image. }
		\end{subfigure}
		
		\begin{subfigure}[b]{\textwidth}
			\begin{tabular}{  c c }
 \\ \tabularnewline  \multicolumn{2}{ c }{\bf QR pivoting solution} \\

				\begin{sideways}  $~~~~~p=K$  \end{sideways}
				\begin{overpic}[width=0.12\textwidth]{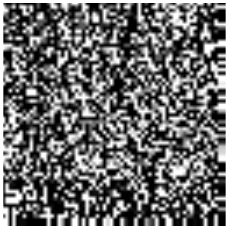}  
				\end{overpic}
				\begin{overpic}[width=0.3\textwidth]{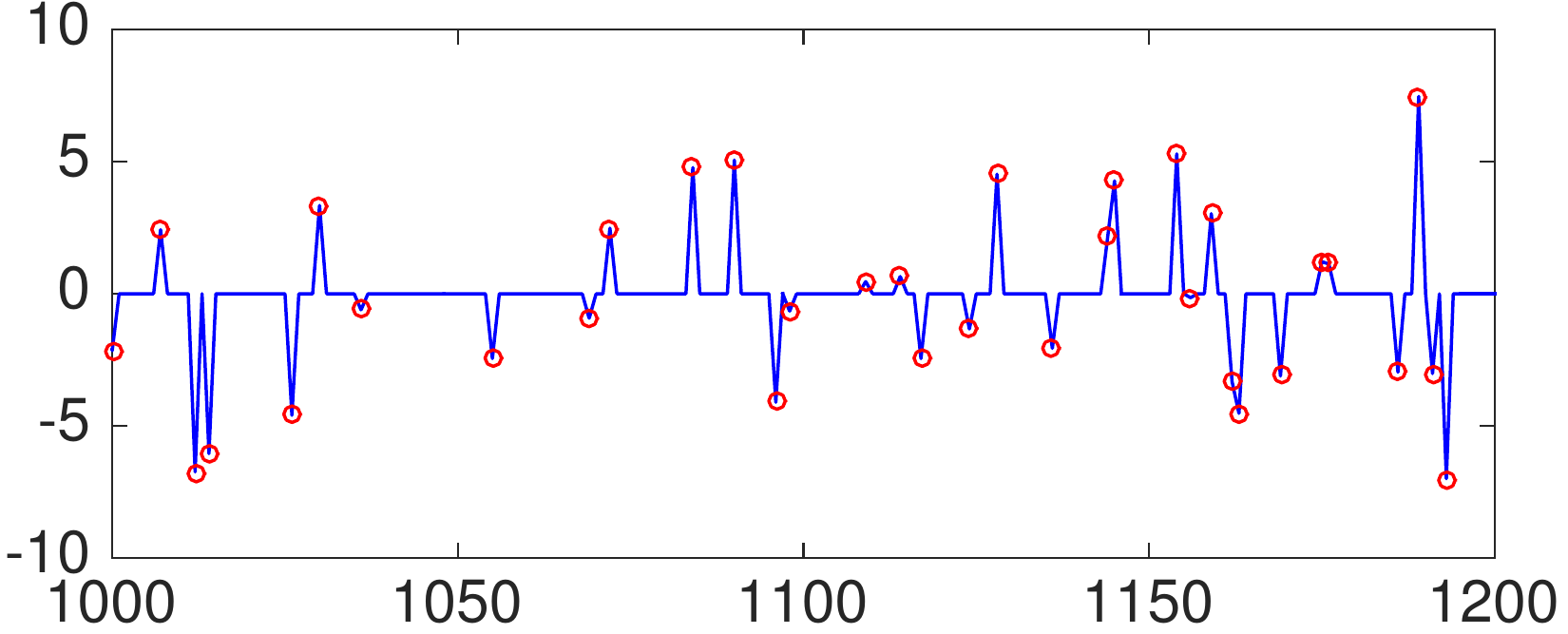} 
				\end{overpic} & \qquad
				\includegraphics[width=0.12\textwidth]{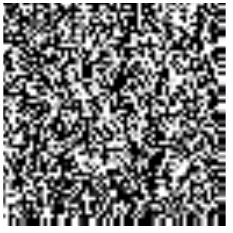}  
				\begin{overpic}[width=0.3\textwidth]{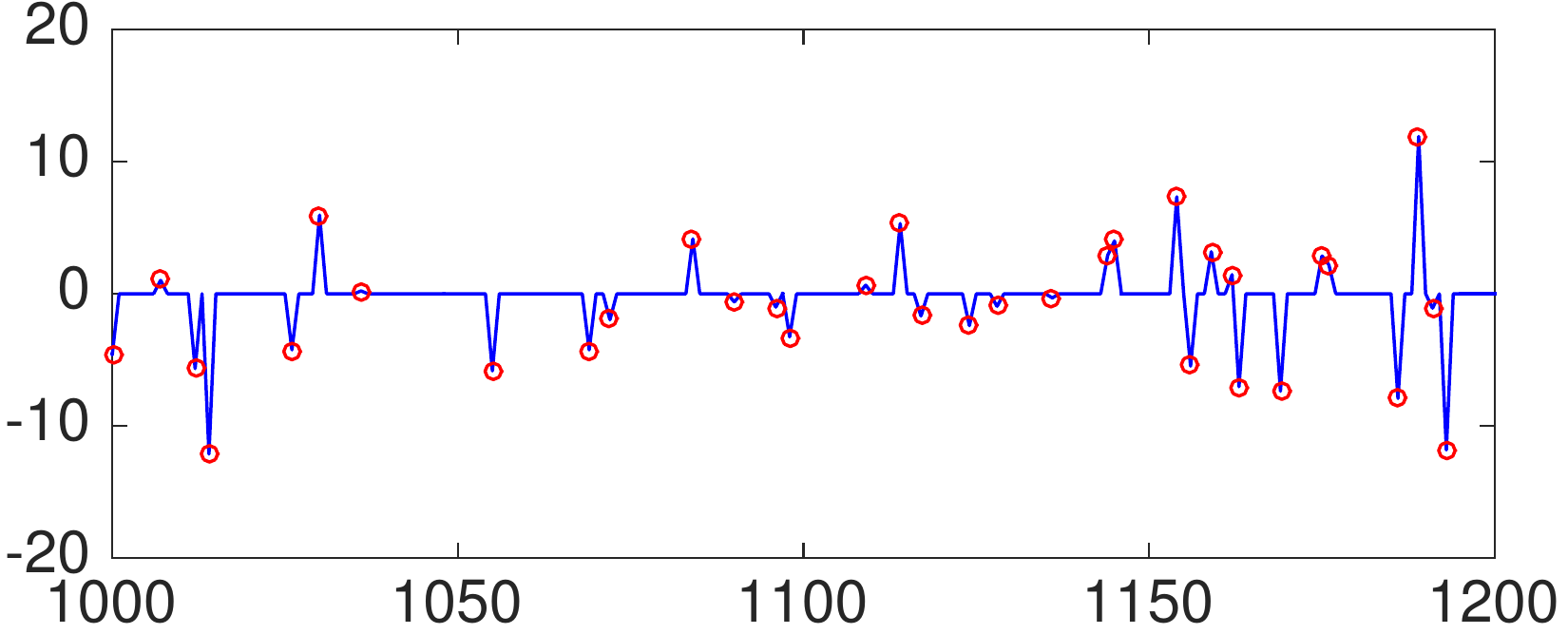}
				\end{overpic}				
				 \\

				\begin{sideways} $~~~p=2000$ \end{sideways}
				\begin{overpic}[width=0.12\textwidth]{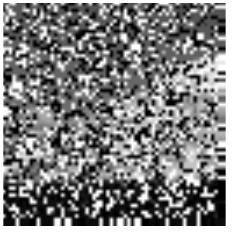}  
				\end{overpic}
				\begin{overpic}[width=0.3\textwidth]{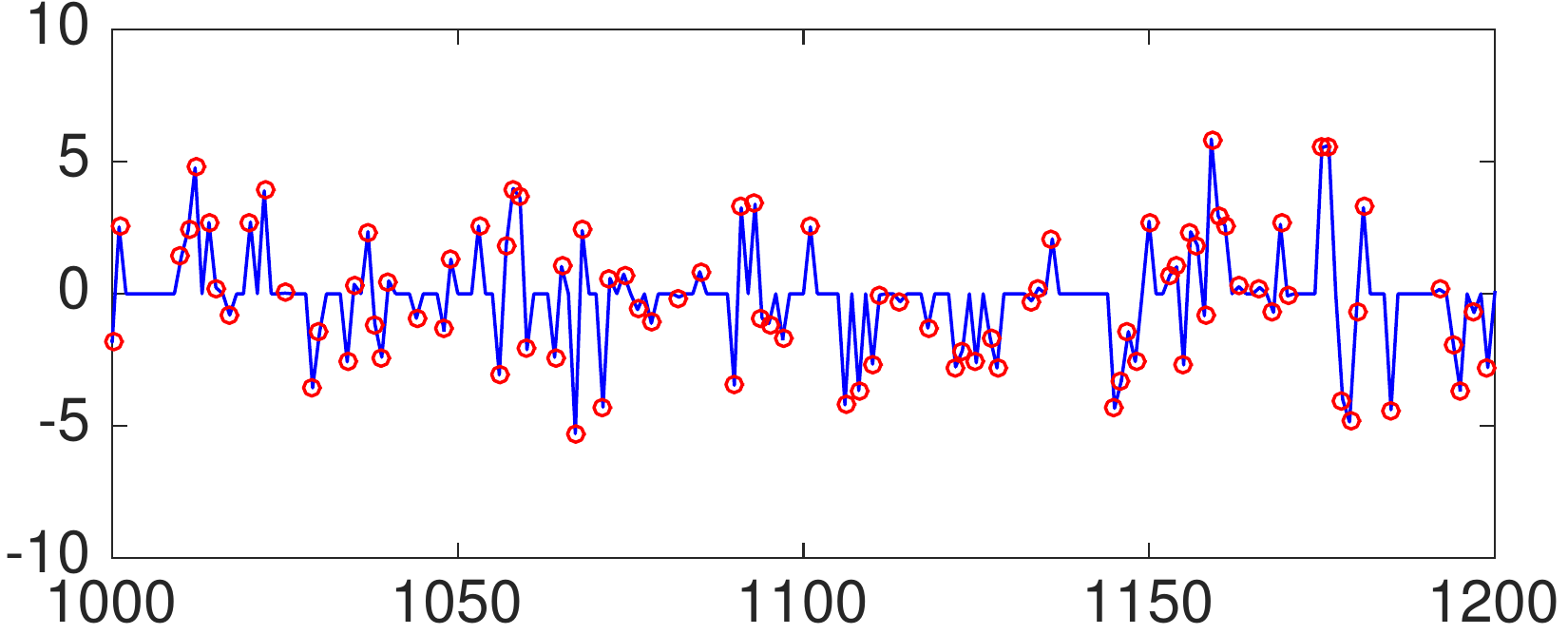} 
				\end{overpic} & \qquad
				\includegraphics[width=0.12\textwidth]{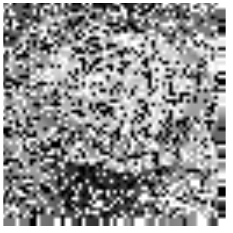}  
				\begin{overpic}[width=0.3\textwidth]{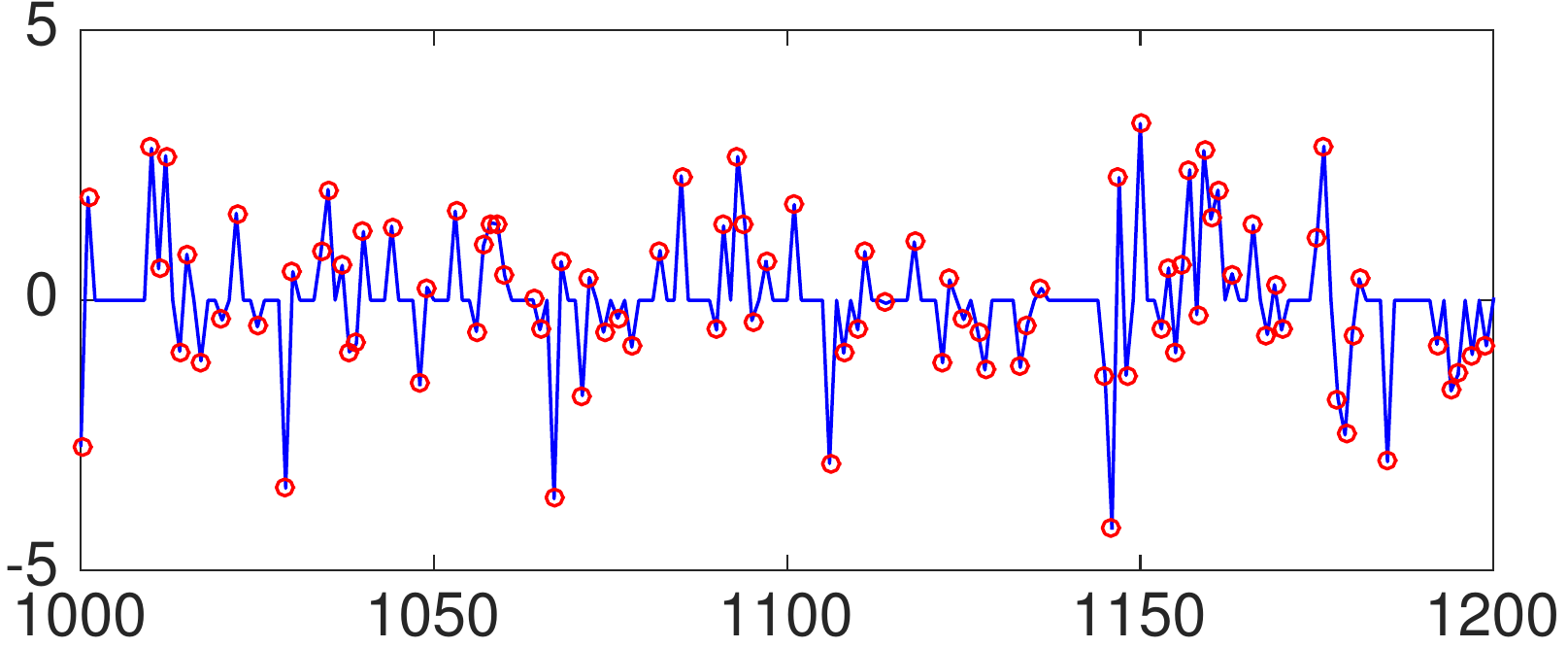}
				\end{overpic}				
				 \\ 		
			\end{tabular}
			\caption{{\bf Reconstruction}. Even at nearly 50\% pixel sampling, least squares regression using the universal Fourier basis is unable to recover the characteristic frequencies. This is because backslash sets the number of nonzero entries in the solution to the number of samples. Additional constraints are needed to recover the true coefficients, which is possible with compressed sensing.} 
			\end{subfigure}	
			\caption{\label{fig:backslash}}		
		\end{center}	
\end{figure}
\end{Sidebar}
\end{figure*}

The QR factorization is implemented and optimized in most standard scientific computing packages and libraries, including Matlab, LAPACK, NumPy, among many others. In addition to software-enabled acceleration, QR runtime can be significantly reduced by terminating the procedure after the first $p$ pivots are obtained. The operation can be accelerated further using randomized techniques, for instance, by the random selection of the next pivot~\cite{drmac2016siam} or by using random projections to select blocks of pivots~\cite{Martinsson2017siamjsc,martinsson2015blocked,duersch2015true}.  
``\nameref{Sidebar_6}" shows why QR pivoting does not find the sparsest vector in an \emph{universal} basis for compressed sensing.

The sparse sensor placement problem we have posed here is related to machine learning
concepts of variable and feature selection~\cite{Tibshirani1996lasso,Guyon2003jmlr}. 
Such sensor (feature) selection concepts generalize to data-driven classification. For image classification using linear discriminant analysis (LDA), sparse sensors may be selected that map via POD modes into the discriminating subspace~\cite{Brunton2016siap}.  Moreover, sparse classification within libraries of POD modes ~\cite{Wright2009ieeetpami,Bright2013pof,Brunton2014siads} can be improved by augmenting DEIM samples with a genetic algorithm~\cite{Sargsyan2015pre} or adapting QR pivots for classification~\cite{Manohar2016jfs}. Sparse sensing has additionally been explored in signal processing for sampling and estimating signals over graphs~\cite{Ribeiro2010sigcomm,DiLorenzo2016ieee,Chen2016ieee,Chepuri2016sam}.

\subsection{Relation to optimal experimental design}
The matrix volume objective described above is closely related to D-optimal experiment design~\cite{doptimal}; in fact, the two problems are identical when regarding the tailored basis $\bPsi_r$ as a set of $n$ candidate experiments of a low-dimensional subspace. Classical experimental design selects the best $p$ out of $n$ candidate experiments for estimating $r$ unknown parameters $\ba \in \reals^r$. Each experiment, denoted $\btheta_i$, produces one output $y_i$ that may include zero-mean i.i.d. Gaussian noise $\boldsymbol{\xi}\sim \mathcal{N}(0,\eta^2)$. Again, we wish to estimate the parameters from $p$ experiment outputs $\by\in\reals^p$ in the following linear system,
\begin{equation}
\label{eqn:doe_est}
\by = \begin{bmatrix}
\rule{0.4cm}{0.4pt} \btheta_{1}\rule{0.4cm}{0.4pt} \\
\rule{0.4cm}{0.4pt} \btheta_{2}\rule{0.4cm}{0.4pt}  \\ 
\vdots \\
\rule{0.4cm}{0.4pt} \btheta_{p}\rule{0.4cm}{0.4pt} 
\end{bmatrix}
\cdot \ba + \boldsymbol{\xi} = 
\sum_{k=1}^r\begin{bmatrix}
\Psi_{\ind_1,k}  \\
 \Psi_{\ind_2,k} \\ 
\vdots \\
 \Psi_{\ind_p,k} 
\end{bmatrix}a_k + \boldsymbol{\xi}
 =
\bC\bPsi_r\ba + \boldsymbol{\xi},
\end{equation}
which is equivalent to the state reconstruction formulation of gappy POD~\cite{Everson1995gappy}. In Matlab notation we sometimes refer to $\bC\bPsi_r$ as $\bPsi_r(\Ind,:)$.
Each possible experiment $\btheta_i$ may be regarded as a row of $\bTheta$ or of the tailored basis $\bPsi_r$ such that $\btheta_i = \bPsi_r(\ind_i,:)$.
Equivalently each $\btheta_i$ is a weighted ``measurement'' of the lower dimensional POD parameter space (not to be confused with the point measurement operation $\bC$). Note that when all experiments are selected the output is simply the state vector $\bx$ since $\bx = \bPsi_r\ba+\boldsymbol{\xi}$.

Given experiment selections indexed by $\Ind$, the estimation error covariance is given by 
\begin{equation}
	\mbox{Var}(\ba-\hat\ba) = \eta^2(\bTheta^T\bTheta)^{-1} =  \eta^2((\bC\bPsi_r)^T\bC\bPsi_r)^{-1}.
\end{equation}
D-optimal subset selection minimizes the error covariance by maximizing the matrix volume of ${\bM_{\Ind}=\bTheta^T\bTheta}$: 
\begin{align}
	\Ind_{\star} &=
	\argmax_{\Ind,|\Ind|=p} \log  \det \sum_{i=1}^p \btheta_{i}^T\btheta_{i}  \nonumber \\
	&= \argmax_{\Ind,|\Ind|=p}  \det (\bC\bPsi_r)^T\bC\bPsi_r,
\end{align}
which is equivalent to~\eqref{eqn:d_opt}. 
Similarly, A-optimal and E-optimal design criteria optimize the trace and spectral radius of $\bTheta^T\bTheta$, and are equivalent to~\eqref{eqn:a_opt} and~\eqref{eqn:e_opt}, respectively.  
The exact solutions of these optimization problems are intractable, and they are generally solved using heuristics. This is most commonly accomplished by solving the convex relaxation with a linear constraint on {\em sensor weights} $\boldsymbol\beta$,
\begin{align}
	\boldsymbol{\beta}_{\star} = \argmax_{\boldsymbol{\beta}\in\reals^n} &\log  \det \sum_{i=1}^n \beta_i\btheta_{i}^T\btheta_{i},  \nonumber \\
	&\subjto \sum_{\substack{i=1\\ 0\le\beta_i\le 1}}^n\beta_i = p.
\end{align}
\begin{table*}
\centering
\captionof{table}{Summary of sensor selection methods} \label{tab:LpTable} 	
		\noindent \begin{tabular}{|l| c| c| c| } 
			\hline
			{\bf Method} & {\bf Objective} & {\bf Runtime} \\
			\hline 			
			D-optimal subset selection  & $\begin{aligned}
			\Ind_{\star} &=
			\argmax_{\Ind,|\Ind|=p} \log \det \sum_{i=1}^p \btheta_{i}^T\btheta_{i} \nonumber \\
			&= \argmax_{\Ind,|\Ind|=p} \det (\bC\bPsi_r)^T\bC\bPsi_r
			\end{aligned}$ & $n\choose p$ determinant evaluations \tabularnewline \hline  
			Convex optimization  &$\begin{aligned}
			\boldsymbol{\beta}_{\star} = \argmax_{\boldsymbol{\beta}} &\log \det \sum_{i=1}^n \beta_i\btheta_{i}^T\btheta_{i} \nonumber \\
			&\subjto \sum_{\substack{i=1 \\ 0\le\beta_i\le 1}}^n \beta_i = p
			\end{aligned}$ & $\mathcal{O}(n^3)$ per iteration  \tabularnewline \hline 
			QR pivoting (greedy)  & $\begin{aligned}[t]\mbox{Case } p&=r: &\bPsi_r^T\bC^T = \bQ\bR \\ \mbox{Case } p&>r : &(\bPsi_r\bPsi_r^T)\bC^T = \bQ\bR\end{aligned}$ & $\begin{aligned}[t]\mathcal{O}(nr^2)  \\ \mathcal{O}(n^3)\end{aligned}$ \\
			\hline
		\end{tabular}
\end{table*}
The optimized sensors are obtained by selecting the largest sensor weights from $\boldsymbol\beta$. 
The iterative methods employed to solve this problem, \ie~convex optimization and semidefinite programs~\cite{Boyd2004convexbook,Joshi2009ieee}, require matrix factorizations of $n\times n$ matrices in each iteration. Therefore they are computationally more expensive than the QR pivoting methods, which cost one matrix factorization in total. Greedy sampling methods such as EIM and QR are practical for sensor placement within a large number of candidate locations in fine spatial grids; hence, they are the methods of choice in reduced-order modeling~\cite{Benner2015siamreview}. The various optimization methods for data-driven sensor selection are summarized in Table~\ref{tab:LpTable}.

\begin{figure*}[t]
\begin{Sidebar}{Sidebar: Other tailored bases -- polynomial interpolation}{Sidebar: Other tailored bases -- polynomial interpolation}
		\label{Sidebar_rev}

Suppose $\mathbf{x}\in\reals^n$ is a suitably fine discretization of the interval $[0,1]$. We may construct a degree $r$ polynomial interpolant of a function $\mathbf{f}(\mathbf{x})$ on this interval by forming the $n\times r$ Vandermonde basis 
\begin{equation}
\bPsi_r = [\mathbf{1} ~|~ \mathbf{x} ~|~ \mathbf{x}^2 ~|~ \dots ~|~ \mathbf{x}^{r-1}].
\end{equation}
It is well-known that equispaced interpolation in the Vandermode basis is ill-conditioned; the most commonly used alternatives are Chebyshev or Legendre bases with non-equispaced points that satisfy a non-uniform density on the interval.
Similar to the discussion above, we seek the (near) best interpolation samples (sensors) for interpolating arbitrary functions within this basis. This is an equivalent formulation to the above, except the basis under consideration is a tailored basis and not a data-driven one, however, the same sampling methodologies apply. The resulting samples, in the non data-driven case, are general enough for sampling any well-behaved univariate function. Below we apply QR sampling on the polynomial basis functions to demonstrate the power of optimized sampling.
\begin{equation}
	\mathbf{f}(\mathbf{x}) = \left|\mathbf{x}^2-\frac{1}{2}\right|
\end{equation} 
This application of QR pivoting to find near-optimal Fekete points for polynomial interpolation was first introduced in~\cite{Sommariva2009qr}. 
Matlab implementation code for polynomial interpolation is listed below, and the comparison between equispaced and QR interpolation samples is shown in Figure~\ref{fig:fekete_pts}.%

\begin{minipage}{.45\textwidth}
\begin{lstlisting}[linerange={4-18,21-21,24-27}]
r = 11; n = 1000;
x = linspace(0,1,n+1)';

% Construct Vandermonde matrix on [0,1]
Vde = zeros(n+1,r);
for i = 1:r
    Vde(:,i) = x.^(i-1);
end

% approximate Fekete points
[Q,R,pivot] = qr(Vde','vector');
fekete = pivot(1:r);

% equispaced points
equi = 1:n/(r-1):(n+1);
pts = fekete; color = 'b';
% build interpolant
f = abs(x.^2-0.5);
coefs = Vde(pts,:)\f(pts);
pstar = Vde*coefs;
\end{lstlisting}
\end{minipage}\quad\quad
\begin{minipage}{.45\textwidth}
\begin{figure}[H]
\centering
	\begin{subfigure}[b]{\textwidth}
		\begin{overpic}[width=\textwidth]{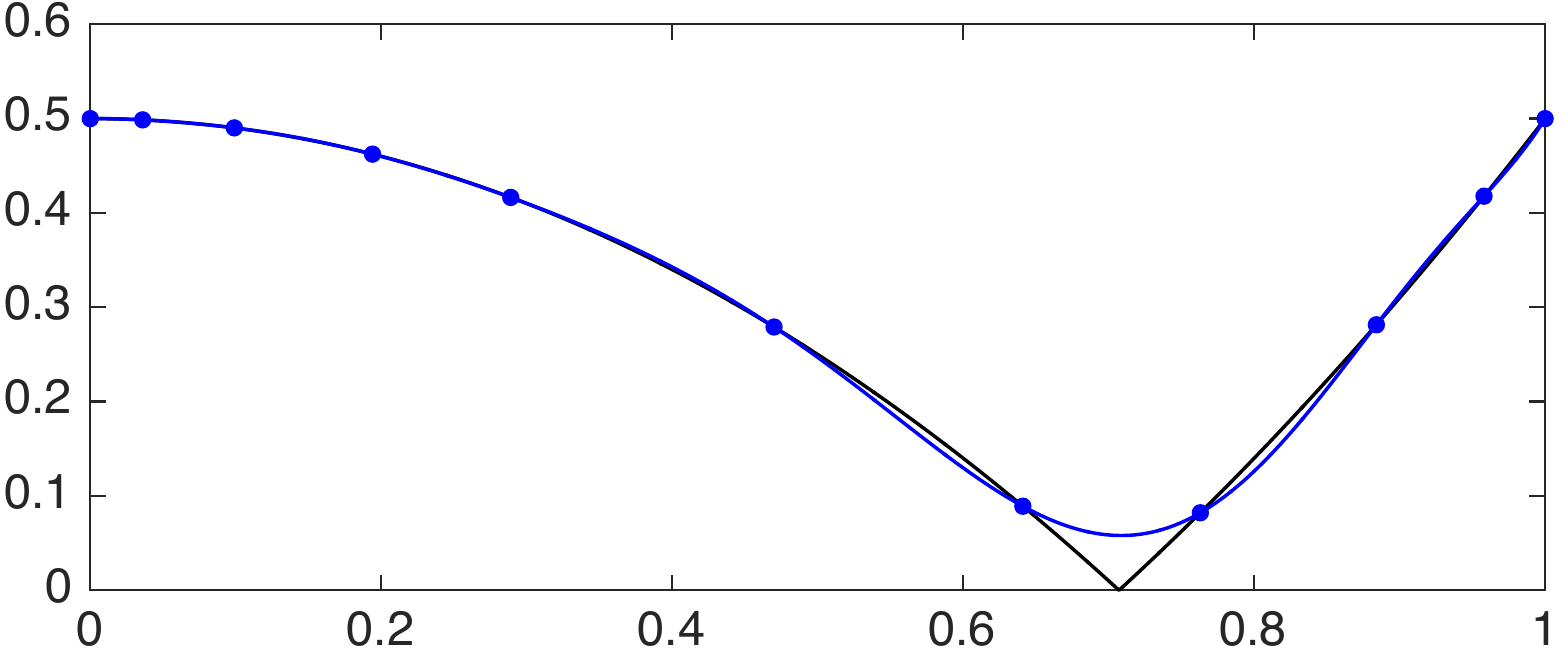}
			\put(7,36){Approximate Fekete points (QR pivots)}
		\end{overpic}
		\begin{overpic}[width=\textwidth]{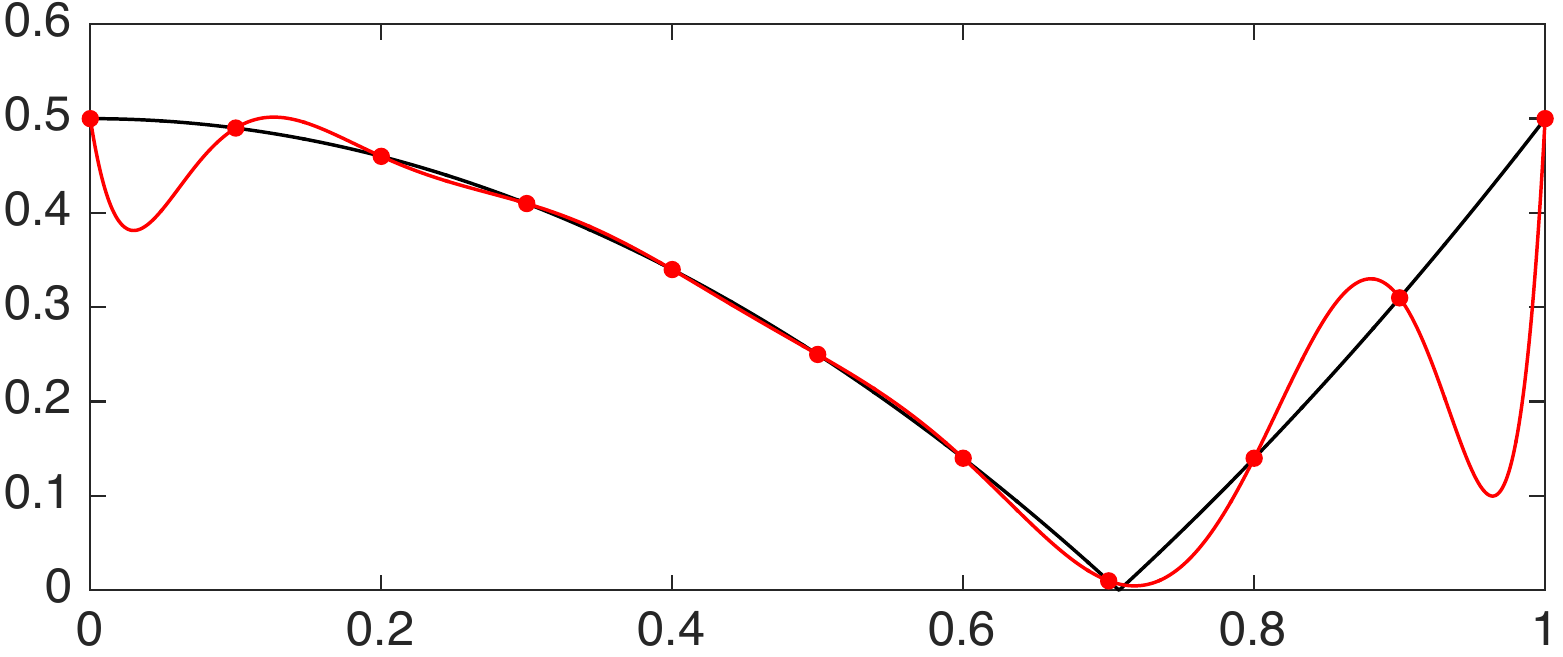}
			\put(7,36){Equispaced points}
			\put(50,-2){$\bx$}
		\end{overpic}		
	\end{subfigure}	
	\caption{Polynomial interpolation with selected QR pivots overcomes the Gibbs phenomenon observed with ill-conditioned equispaced points.	
\label{fig:fekete_pts}}
\end{figure}
\end{minipage}

\end{Sidebar}
\end{figure*}

\section{Comparison of methods}\label{Sec:Methods}
Sensor selection and signal reconstruction algorithms are implemented and compared on data from fluid dynamics, facial images, and ocean surface temperatures.  
The examples span a wide range of complexity, exhibit both rapid and slow singular value decay, and come from both static and dynamic systems. 

In each example, optimized sensors obtained in a tailored basis with QR pivots (See~``\nameref{Sec:QR}") outperform random measurements in a universal basis using compressed sensing (See~``\nameref{Sec:CS}") for signal reconstruction.  
Moreover, for the same reconstruction performance, many fewer QR sensors are required, decreasing the cost associated with purchasing, placing, and maintaining sensors, as well as reducing the latency in computations.  
Thus, for a well-scoped reconstruction task with sufficient training data, we advocate principled sensor selection rather than compressed sensing. 
For example, the QR-based sampling method is demonstrated with yet another tailored basis commonly encountered in scientific computing -- the Vandermonde matrix of polynomials. ``\nameref{Sidebar_rev}" compares polynomial interpolation with QR pivots for the ill-conditioned set of equispaced points on an interval.
When the structure of the underlying signal is unknown, then compressed sensing provides more flexibility with an associated increase in the number of sensors.  

\subsection{Flow past a cylinder}

Fluid flow past a stationary cylinder is a canonical example in fluid dynamics that is high-dimensional yet reveals strongly periodic, low-rank phenomena. 
It is included here as an ideal system for reduction via POD and hence, minimal sensor placement. 
The data is generated by numerical simulation of the linearized Navier-Stokes equations using the immersed boundary projection method (IBPM) based on a fast multi-domain method~\cite{taira:07ibfs,taira:fastIBPM}. 

The computational domain consists of four nested grids so that the finest grid covers a domain of $9\times 4$ cylinder diameters and the largest grid covers a domain of $72\times 32$.  Each grid has resolution $450\times 200$, and the simulation consists of 151 timesteps with ${\delta t=0.02}$.  The Reynolds number is 100, and the flow is characterized by laminar periodic vortex shedding~\cite{Noack2003jfm}.

\begin{figure*}[t]
\begin{Sidebar}{Sidebar: Selecting number of sensors $p$ and rank $r$ for flow past a cylinder}{Sidebar: Selecting number of sensors $p$ and rank $r$ for flow past a cylinder}\
		\label{Sidebar_7}
		An overarching goal of optimized sensor placement is choosing the fewest $p$ sensors for reconstruction.  
		This sparse sensor optimization is facilitated by low-rank structure in the data (Fig.~\ref{fig:ibpm_spect}) and inherently involves a trade-off between the number of sensors and reconstruction accuracy.  
		As seen in Fig.~\ref{fig:ibpm_conv}, effective sensor optimization moves the elbow of the reconstruction error curve down and to the left, indicating accurate reconstruction with few sensors.  
		Reducing the number of sensors may be critically enabling when sensors are expensive or when low computational latency is desired.  
		\begin{figure}[H]
		\centering
			\begin{subfigure}[b]{0.45\textwidth}
				\centering
				\begin{overpic}[width=\textwidth]{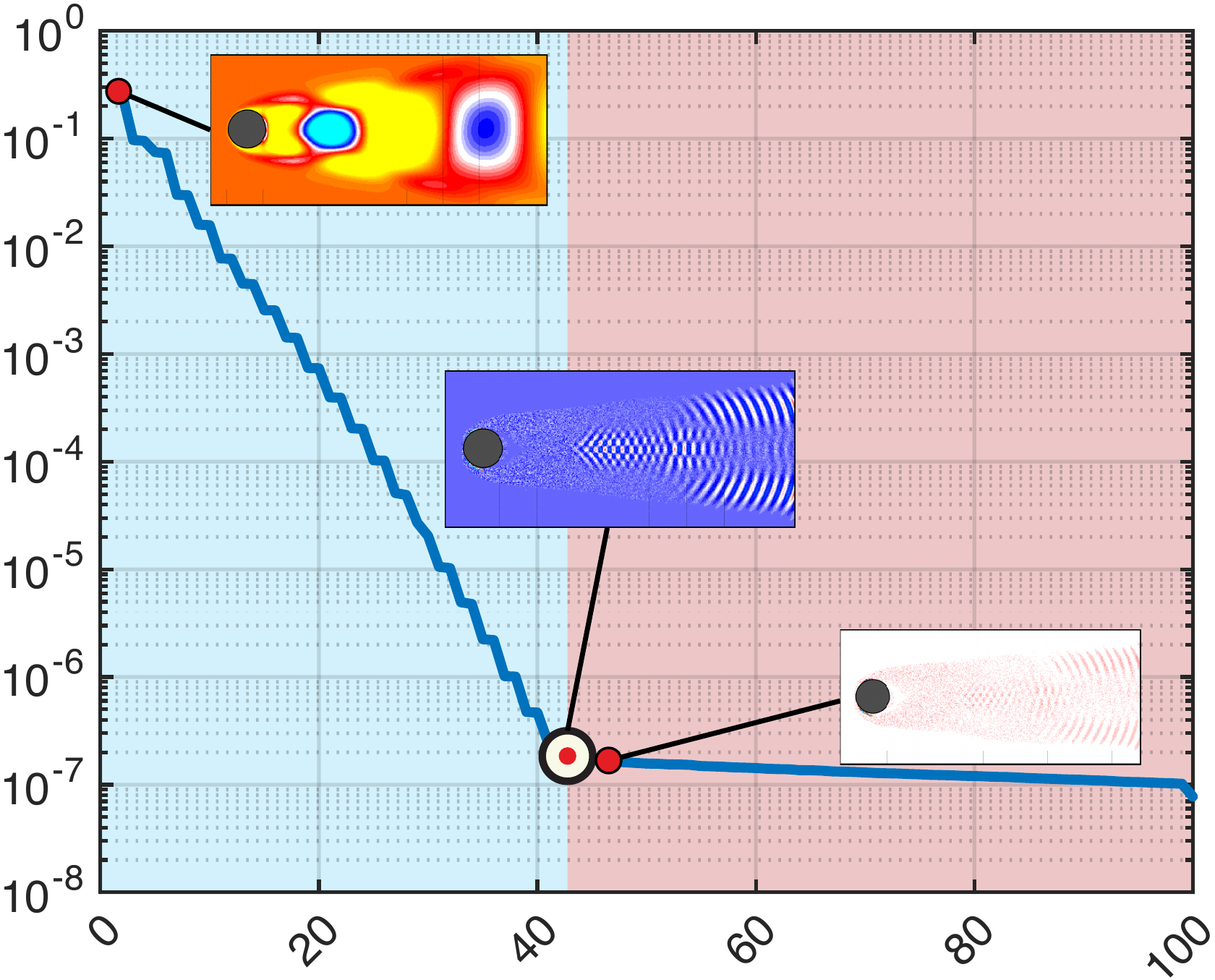}
					\put(-2.5,40){\makebox(0,0){\rotatebox{90}{\small Normalized singular values}}}
					\put(45,-2){\small $r$}
					\put(22,58){Mode 1}
					\put(43,52){Mode 42}
					\put(72,30){Mode 43}
				\end{overpic}
				\caption{Singular value spectrum. \label{fig:ibpm_spect}}
			\end{subfigure}
			\quad
			\begin{subfigure}[b]{0.45\textwidth}		
				\begin{overpic}[width=\textwidth]{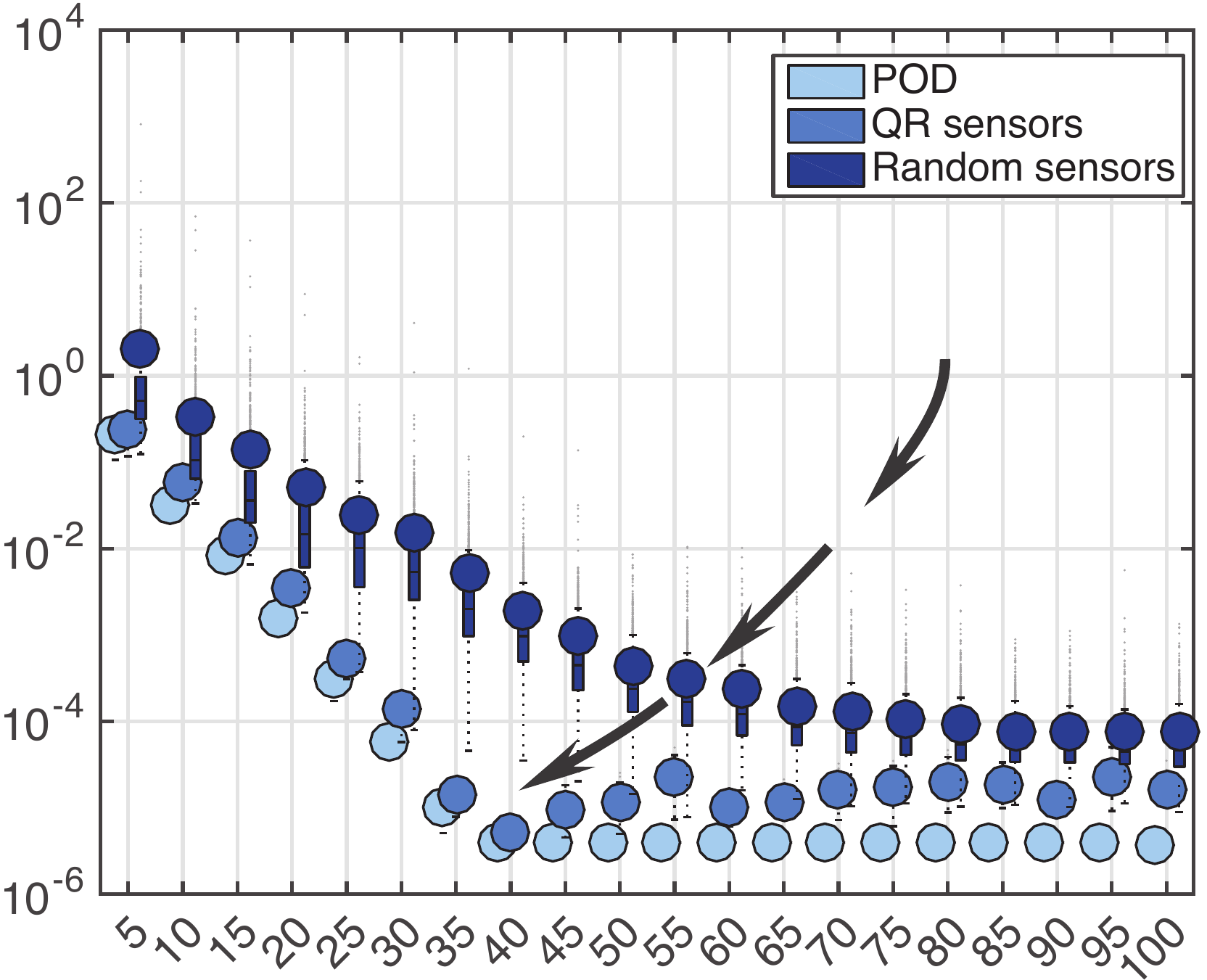}						
					\put(-2.5,40){\makebox(0,0){\rotatebox{90}{\small Reconstruction error}}}
					\put(45,-2){\small $r$}
					\put(53,59){Goal: Fewer sensors}
					\put(63,54){and lower error}
				\end{overpic}
				\caption{Reconstruction error with QR pivot sensors.\label{fig:ibpm_conv}}	
			\end{subfigure}	
			\caption{Illustration of singular value decay and sparse signal reconstruction in the fluid flow past a cylinder example.  (a) The exponential decay of singular values indicates low-rank dynamics with a modal truncation at $r=42$ modes, the last modes containing any meaningful spatial structure. (b) Recovery with a minimal number of QR sensors with increasing POD basis rank is nearly as good as a snapshot's approximation in the same POD basis, and overfitting occurs beyond $r>42$.  Optimized sensing remains orders of magnitude more accurate than reconstruction with random sensors.} 
\end{figure}
\end{Sidebar}
\end{figure*}

\begin{figure}[H]
	\centering
	\begin{overpic}[width=0.325\textwidth]{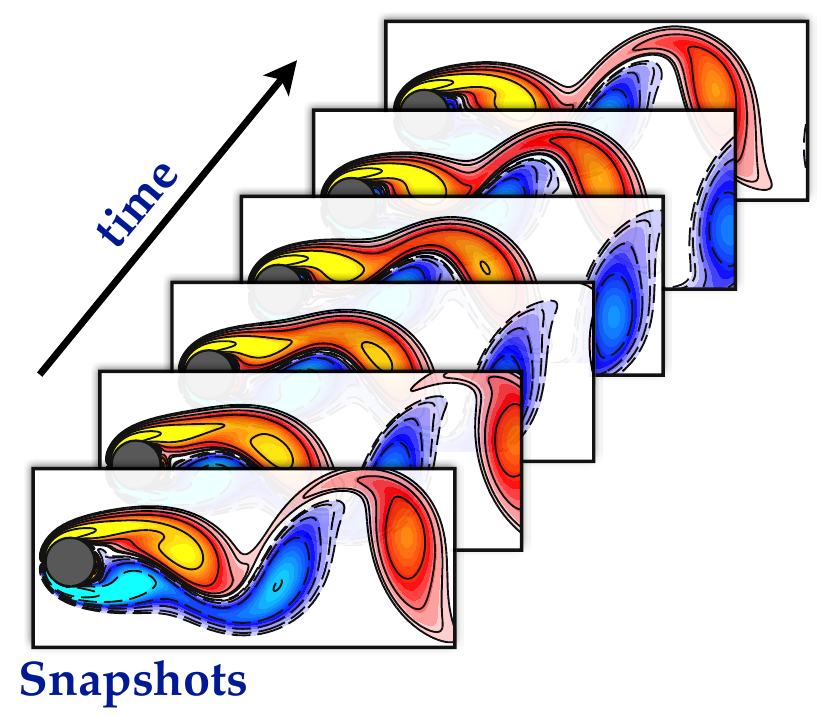}
		\put(57,9){$\bx_1$}
		\put(65.,22){$\bx_2$}
		\put(99.5,64){$\bx_m$}
	\end{overpic}
	\caption{Fluid flow past a cylinder at Reynolds number $100$, characterized by vortex shedding.}\label{Fig:VortexShedding}
\end{figure}

\begin{figure*}[t!]
		\centering
		\begin{overpic}[width=0.8\textwidth]{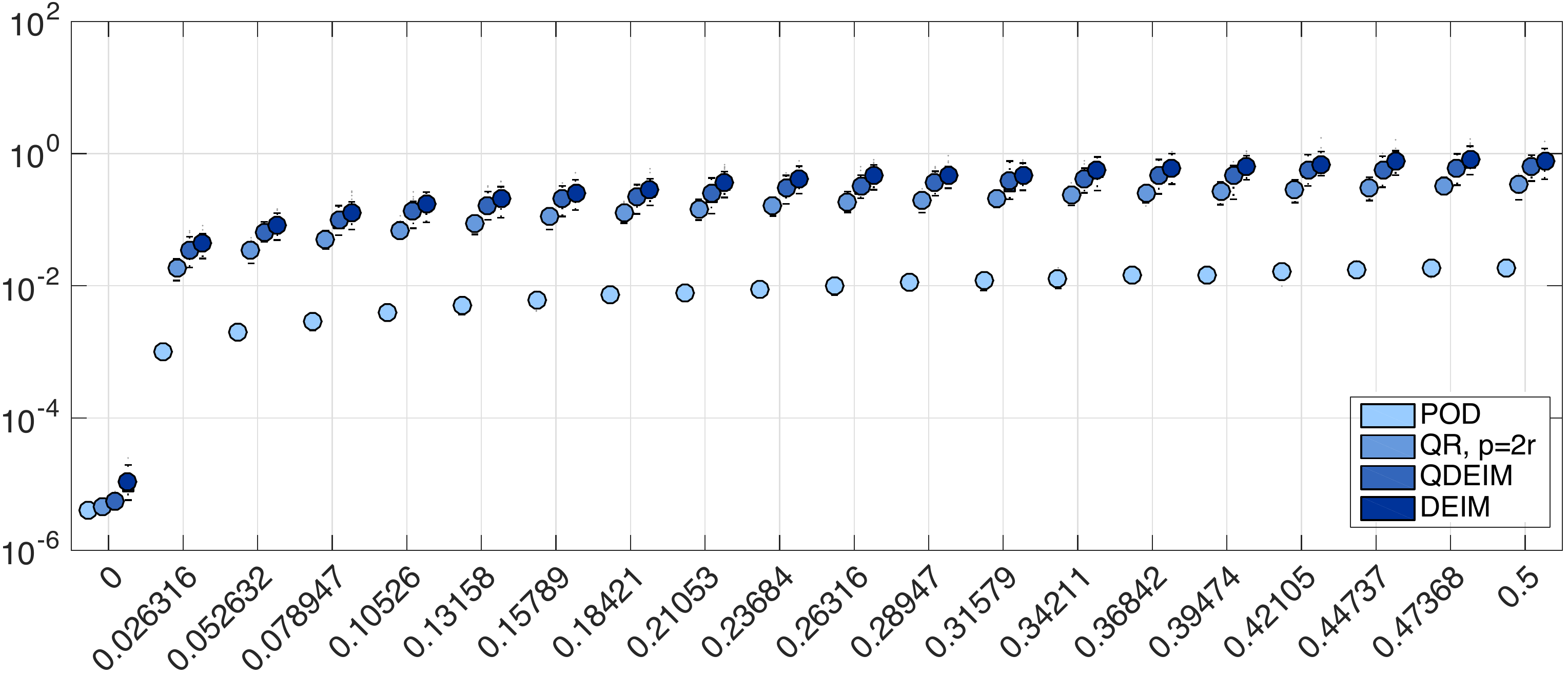}
			\put(-2.5,20){\makebox(0,0){\rotatebox{90}{\small Reconstruction error}}}
			\put(45,-2){Sensor noise variance $\eta$}
		\end{overpic}
		\vspace{.1in}
		\caption{Comparison with DEIM. Oversampled QR permits additional sensors $p=2r$ and a 4x reduction in reconstruction error compared to discrete empirical interpolation method (DEIM). In the $p=r=40$ case DEIM and QR pivoting (Q-DEIM) perform comparably, with a slight improvement observed with QR. All sampling methods are within some constant factor of the optimal POD approximation with the full states. Hence POD-based sampling methods demonstrate robust, bounded growth of reconstruction error as sensor noise increases. \label{fig:ibpmnoise}}
		\vspace{.1in}
\end{figure*}

Vorticity field snapshots are shown in Fig.~\ref{Fig:VortexShedding}.  
In the cylinder flow and sea surface temperature examples, each snapshot $\bx_i=\bx(t_i)$ is a spatial measurement of the system at a given time $t_i$. Thus POD coefficients $a_k(t_i)$ are time dependent, and $\bpsi_k(x)$ are spatial eigenmodes. The first 100 cylinder flow snapshots are used to train POD modes and QR sensors, and reconstruction error bars are plotted over 51 remaining validation snapshots in Figures~\ref{fig:ibpm_conv} and \ref{fig:ibpmnoise}.
The POD modes of this data reflect oscillatory dynamics characterized by periodic vortex-shedding.  
The data is low-rank, and the singular values decay rapidly, as shown in Fig.~\ref{fig:ibpm_spect}.  
The singular values occur in pairs corresponding to harmonics of the dominant vortex shedding frequency.  
Most of the spectral energy in the dataset is captured by the first 42 POD modes. Thus the intrinsic rank of the dataset in POD feature space is $r=42$, and the minimal number of QR pivots is $p=42$. This yields a dramatic reduction from the initial state dimension of $n\approx 90000$ spatial gridpoints. Here, QR pivoting of $\bPsi_r^T$ (with $O(nr^2)$ operations) bypasses expensive $O(n^3)$ factorizations of large $n\times n$ matrices with alternate sampling or convex optimization methods. 
	
Reconstruction from QR sensors (Fig.~\ref{fig:ibpm_conv}) successfully captures modal content with only ${p=r}$ sensors when fitting to the first 42 POD modes. The first 42 POD modes characterize nearly 100\% of the system's energy, the normalized sum of the singular values. 
Using modes beyond $r>42$ results in overfitting, and QR pivoting selects sensors based on uninformative modes.  
Thus, accuracy stops improving beyond ${r=42}$ target modes, whereupon sensor data amplifies noise. However, these tailored sensors perform significantly better than random sensors due to the favorable conditioning properties of QR interpolation points. 
``\nameref{Sidebar_7}" illustrates the ability of optimized sensing to significantly reduce the number of sensors required for a given performance. 
	
\subsection{Noise comparison study}
\begin{figure*}[t!]
	\centering
	\begin{overpic}[width=0.8\textwidth]{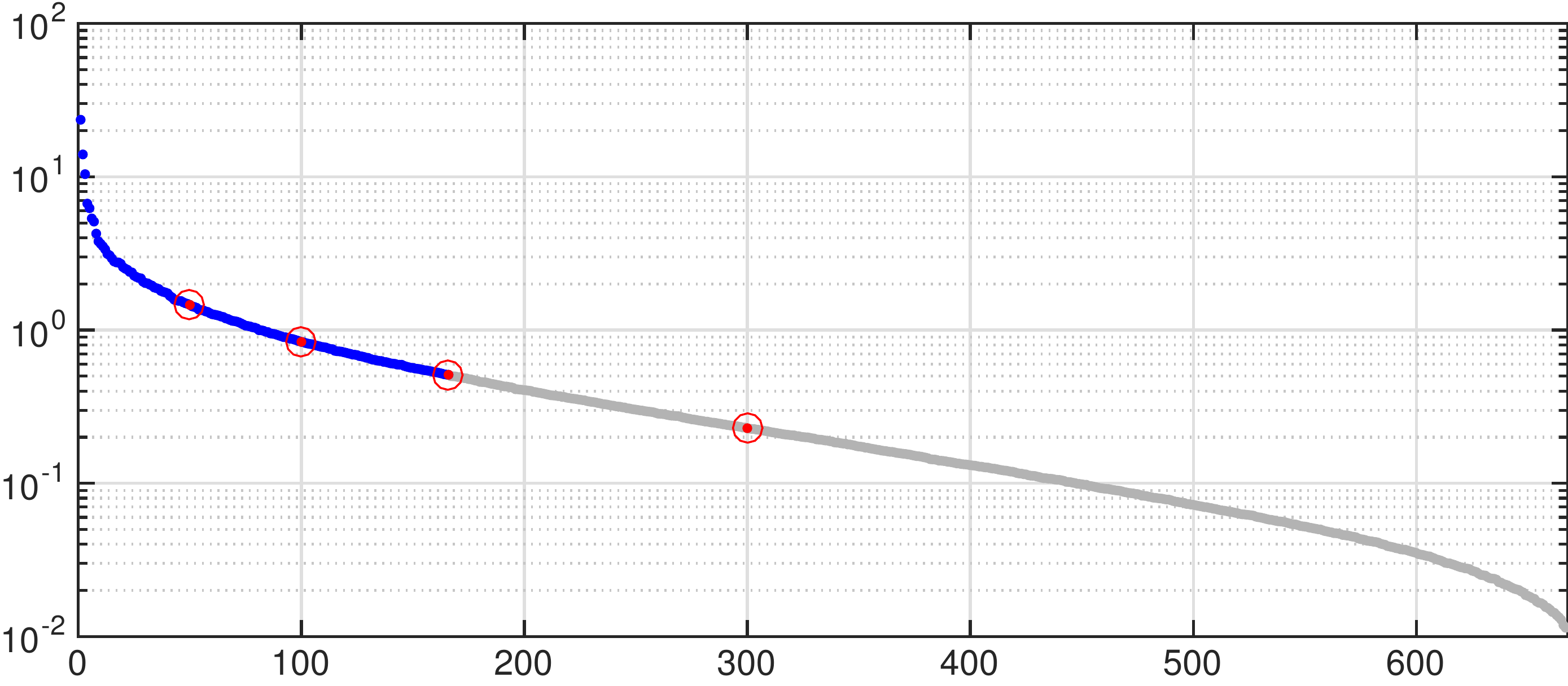}
		\put(7,10){\includegraphics[width=0.07\textwidth]{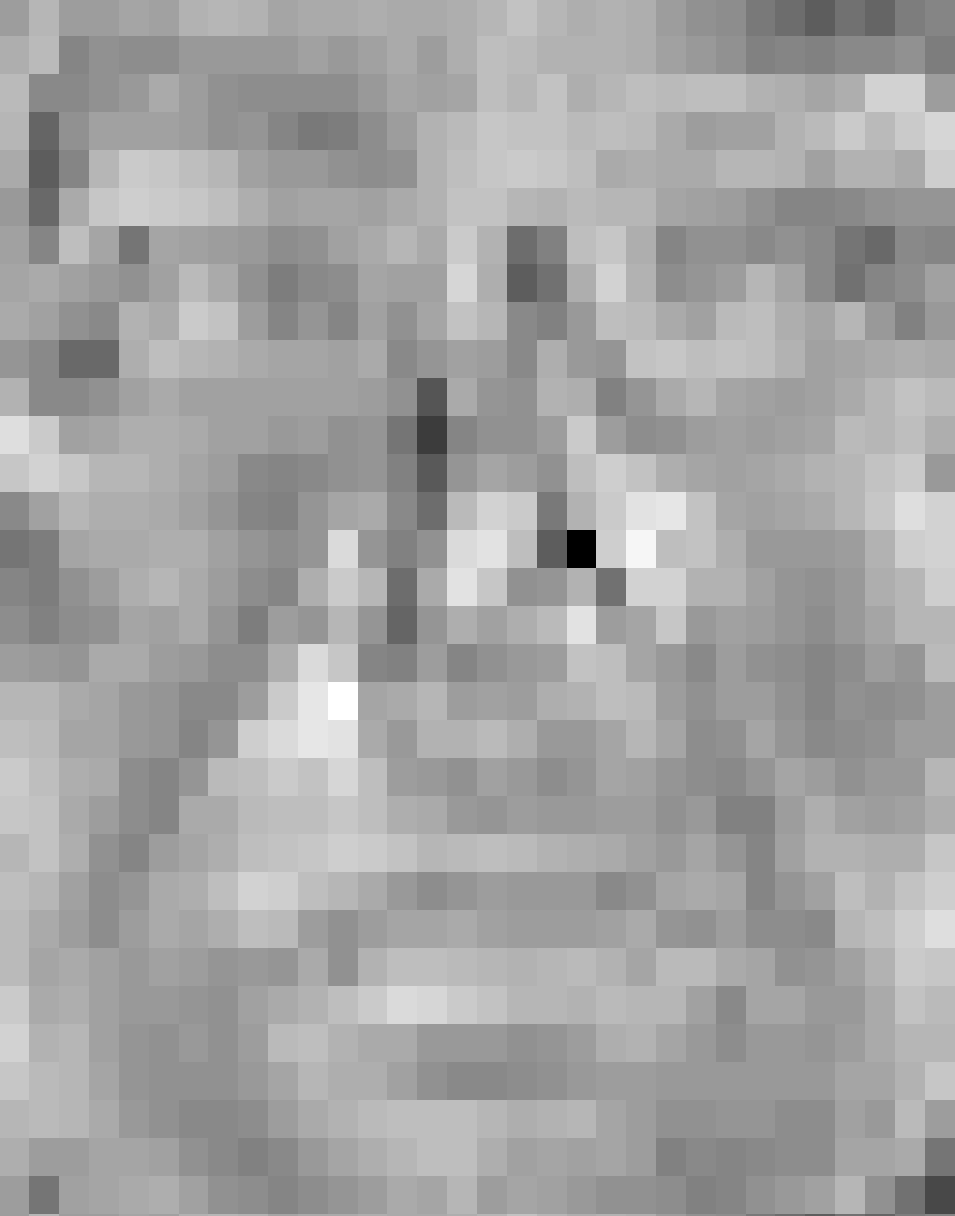}}
		\put(8,28){\small Mode 50}
		\put(16,7){\includegraphics[width=0.07\textwidth]{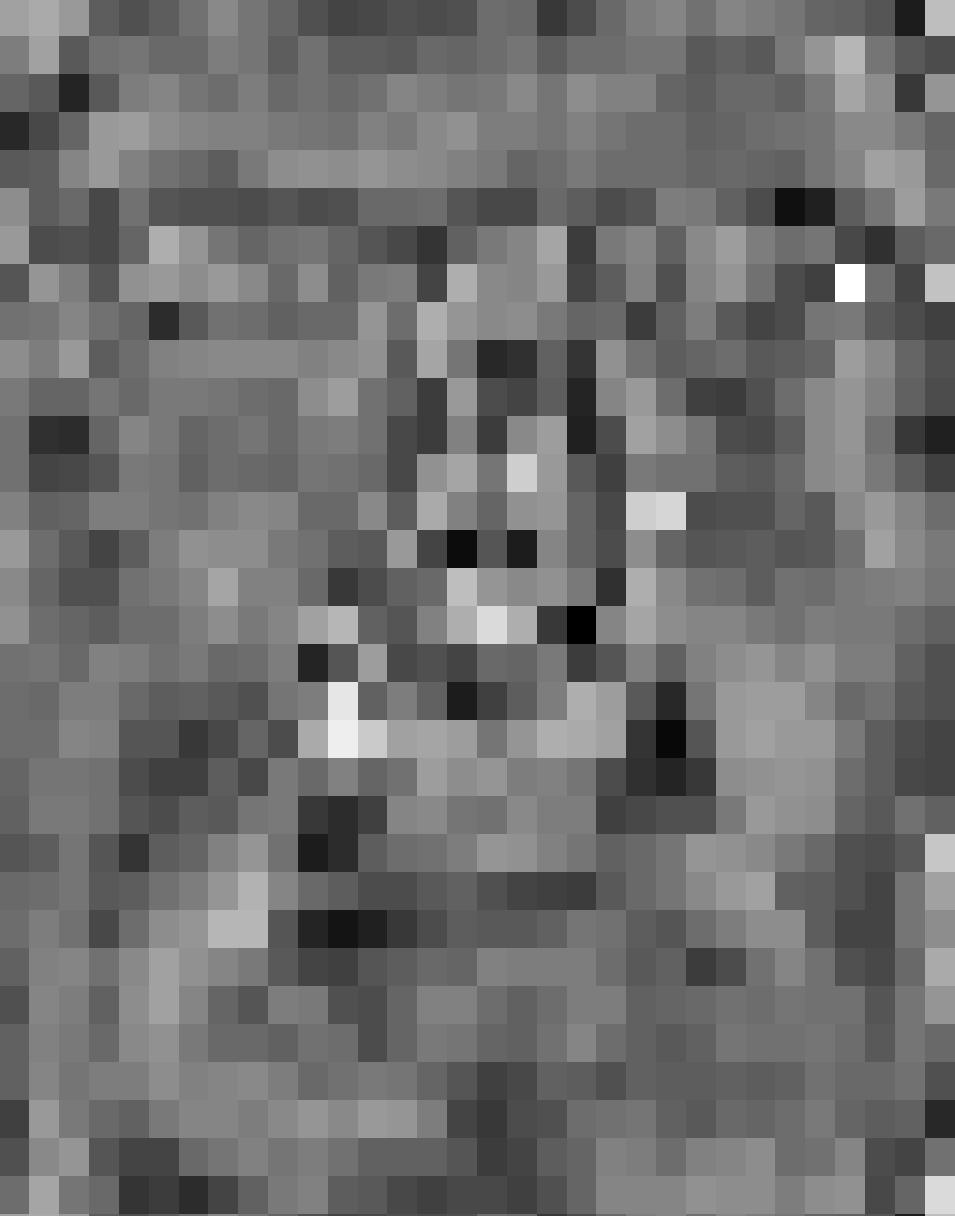}}
		\put(20,24){\small Mode 100}
		\put(25,5){\includegraphics[width=0.07\textwidth]{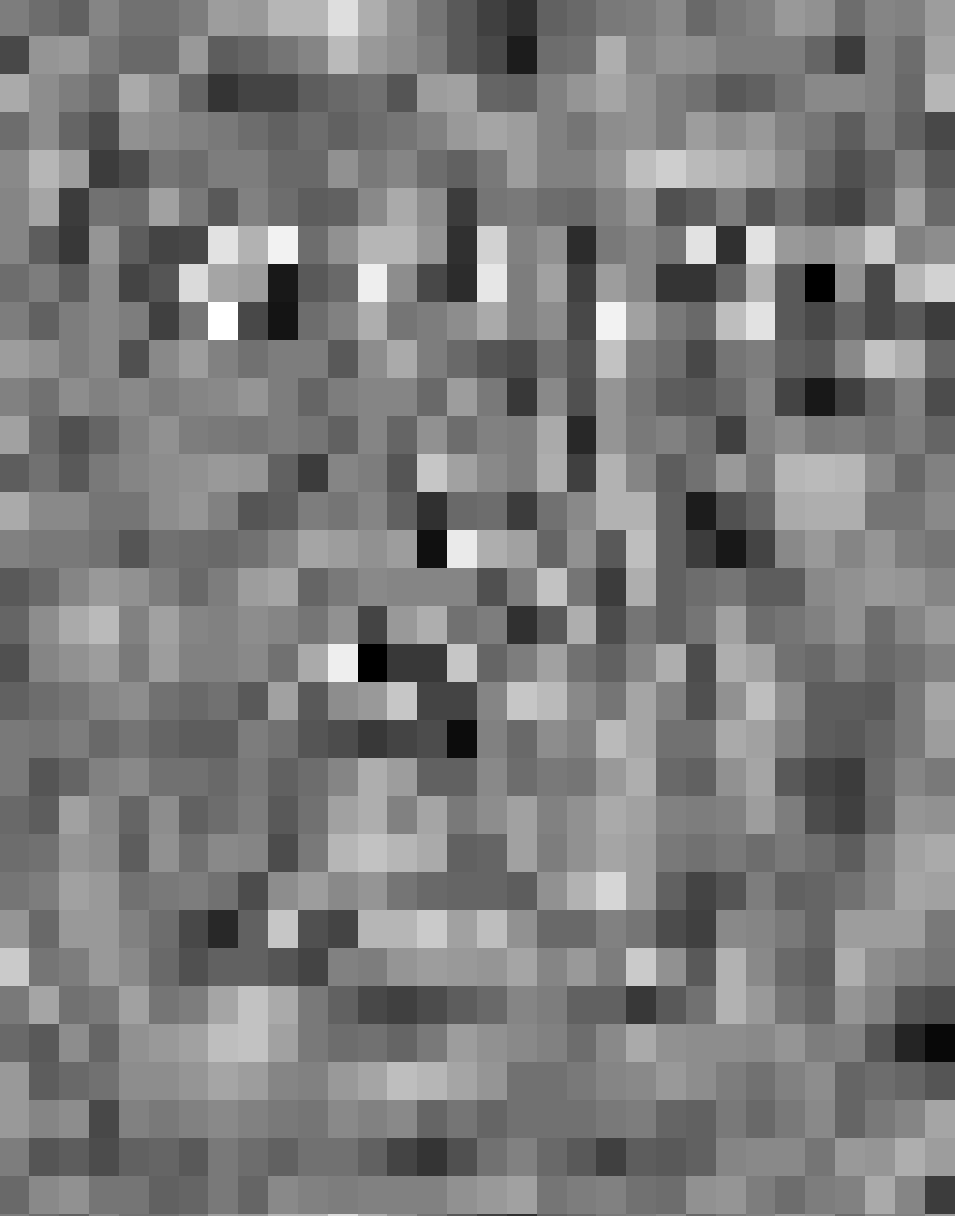}}
		\put(30,20){\small Mode 166}
		\put(43,4){\includegraphics[width=0.07\textwidth]{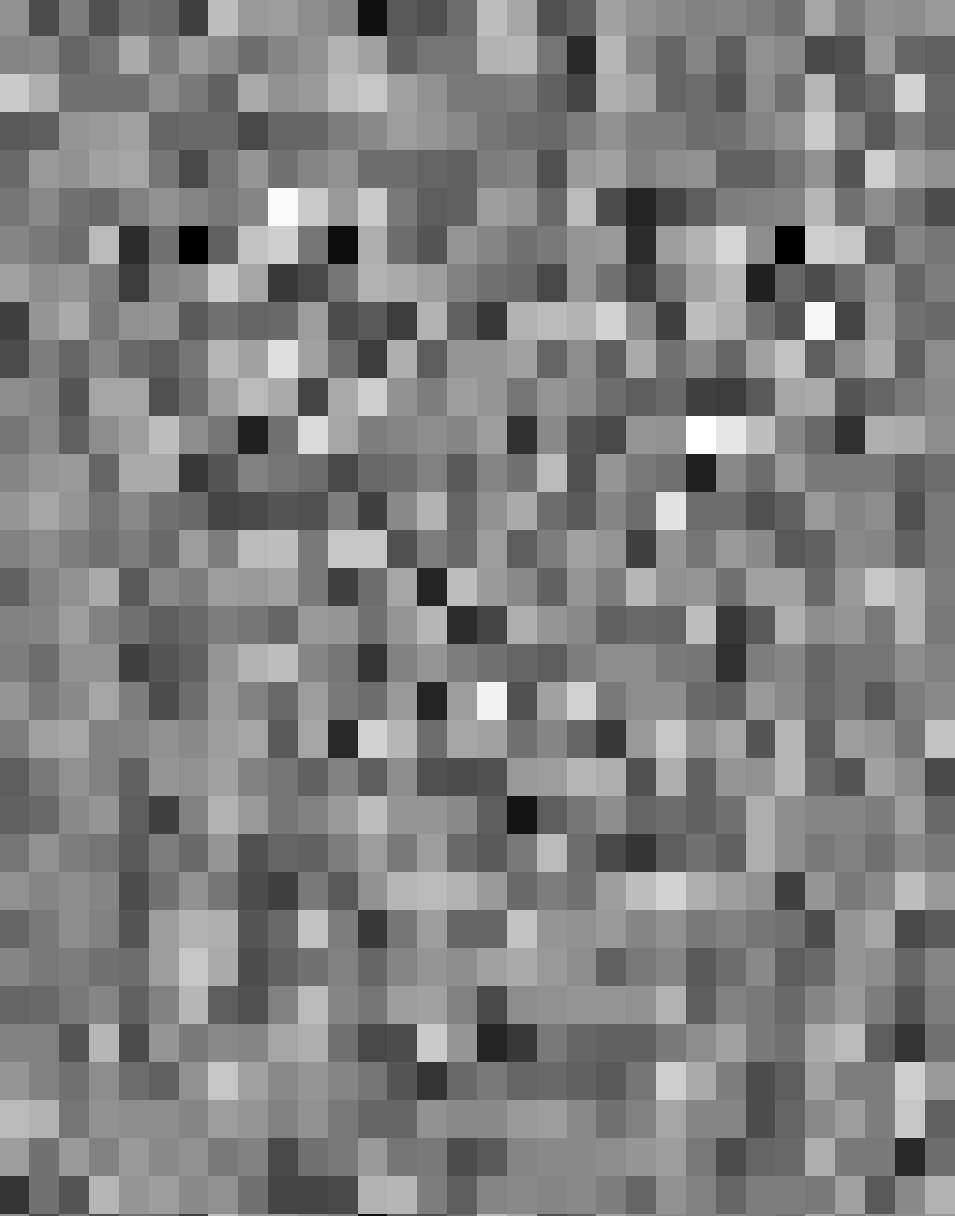}}
		\put(43,18){\small Mode 300}
		\put(50,-1.5){$r$}
		\put(-2,22){\large$\sigma_r$}
	\end{overpic}
	\vspace{.1in}
	\caption{Normalized singular values and selected eigenfaces. Facial information progressively decreases across selected eigenfaces, and no facial features can be readily discerned beyond eigenface $r=166$, the optimal modal truncation value according to~\cite{Gavish2014ieeetit}. \label{fig:yalebspectrum}}
	\vspace{.1in}
\end{figure*}
	Measurements of real-world data are often corrupted by sensor noise.
	The POD-based sensor selection criteria, as well as A,D and E-optimal experimental design criteria, are optimal for estimation with measurements corrupted by zero-mean Gaussian white noise. 
	We empirically demonstrate this on the cylinder flow data with increasing additive white noise. 
	Here we assume sensor noise only in the test measurements and not in the training data or features, see Eqn.~\eqref{eqn:doe_est}.
	The POD modes and the different sensor sets are trained on the first 100 snapshots, and these different sensor sets are used to reconstruct the remaining 50 validation snapshots, which were not used for training features. 
	
	The reconstruction accuracy of the various sampling methods are compared for increasing sensor noise in Fig.~\ref{fig:ibpmnoise},  alongside the full-state POD approximation for illustration. 
	Here we truncate the POD expansion to $r=40$ eigenmodes and compare the $p=r$ reconstruction computed with the discrete empirical interpolation method (DEIM)~\cite{Chaturantabut2010siamjsc} against the QR pivoting reconstruction (Q-DEIM, $p=r$). 
	The DEIM greedy strategy places sensors at extrema of the residual computed from approximating the $k$-th mode with the previous $k-1$ sensors and eigenmodes. 
	It can be seen that QR reconstruction is slightly more accurate than that of DEIM, which is the leading sampling method currently in use for reduced-order modeling~\cite{Benner2015siamreview}. 
	
	QR pivoting is competitive in both speed and accuracy. 
	The speed of QR pivoting is enabled by several implementation accelerations; for example, the column norms in line 4 of Algorithm 2 are only computed once and are then reused. 
	Unlike QR pivoting, DEIM and related methods add successive sensors per iteration by similarly optimizing some metric over all locations. 
	However, this metric (e.g., the approximation residual or the largest singular value) is recomputed at every iteration. 
	The QR factorization is significantly faster than convex optimization methods used in optimal design of experiments, which typically require one matrix factorization per iteration.

	In fact, convex optimization methods that relax the subset selection to weighted sensor placement provide no bounds for deviation from the global optimum, partly because rounding procedures are employed on the weights to decide the final selection. However, reconstruction error bounds for the globally optimal selection are known for DEIM~\cite{Chaturantabut2010siamjsc}, Q-DEIM~\cite{drmac2016siam} and related POD sampling methods~\cite{Willcox2006compfl,Yildirim2009oceanmod}. 
	Furthermore, QR pivoting can achieve significant accuracy gains over DEIM by oversampling -- when $p=2r$ QR reconstruction error is 4x smaller than that of DEIM. It should be noted that while DEIM has not yet been extended to the $p>r$ case, oversampling methods exist for other POD-sampling methods~\cite{Willcox2006compfl,Yildirim2009oceanmod}. 
	However, the iterative procedures involved in the latter are typically more expensive. 
	Recent accelerated variants of greedy principled sampling~\cite{Zimmermann2016siamjsc} may permit oversampling for large $n$, when oversampled QR storage requirements would be excessive. 
	In the cylinder flow case, we bypass this storage requirement by uniformly downsampling the fine grid by a factor of 5 in each spatial direction, thus reducing the number of candidate sensor locations to $n=3600$ instead of $n=89351$.
	
\subsection{Extended Yale B eigenfaces}\label{Sec:Results:Eigenfaces}

\begin{figure*}[t]
	\centering
	\begin{tabular}{|c c c|}
		\hline {\bf $\ell_2$ QR sensors} & {\bf $\ell_2$ Random sensors} & {\bf Compressed sensing (random)}\\
		{\bf \quad 50 \hspace{0.5in} 166} &{\bf 50 \hspace{0.5in} 166} & {\bf 50 \hspace{0.5in} 166 \hspace{0.5in} 300 \hspace{0.5in} 600} \\
		\begin{sideways} \quad~ \bf Sensors \end{sideways}~
		\includegraphics[width=.105\textwidth]{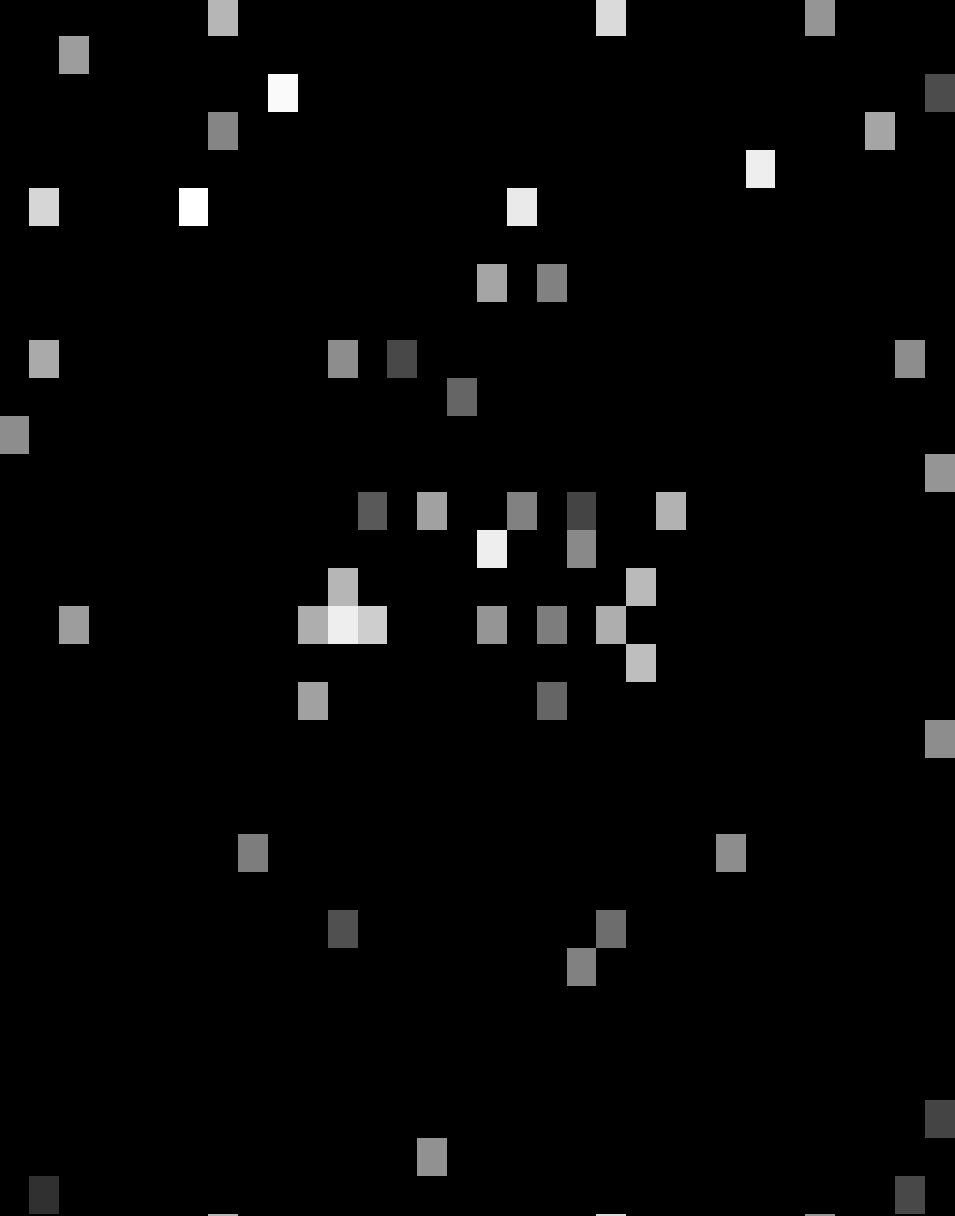} 
		\includegraphics[width=.105\textwidth]{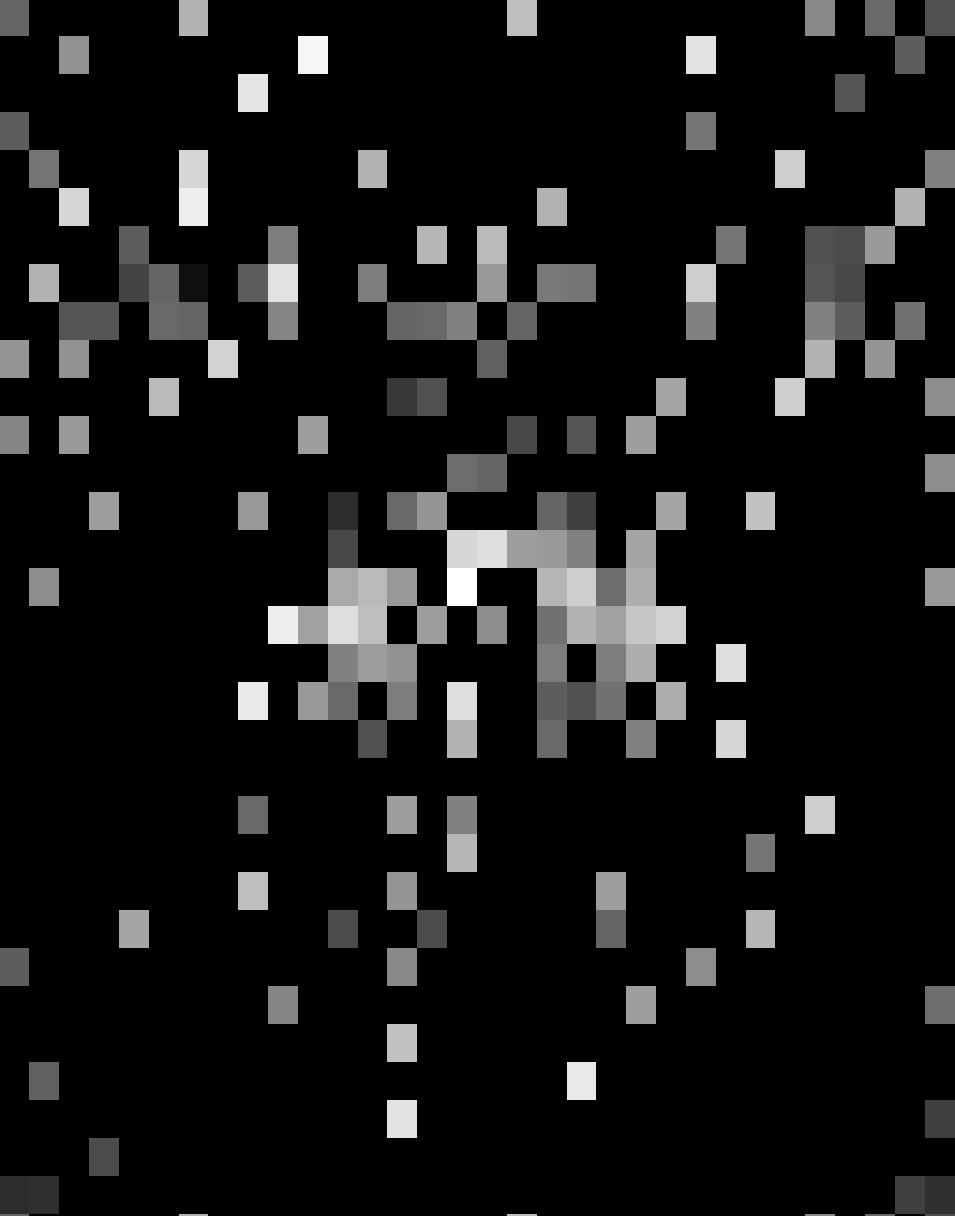} &~
		\includegraphics[width=.105\textwidth]{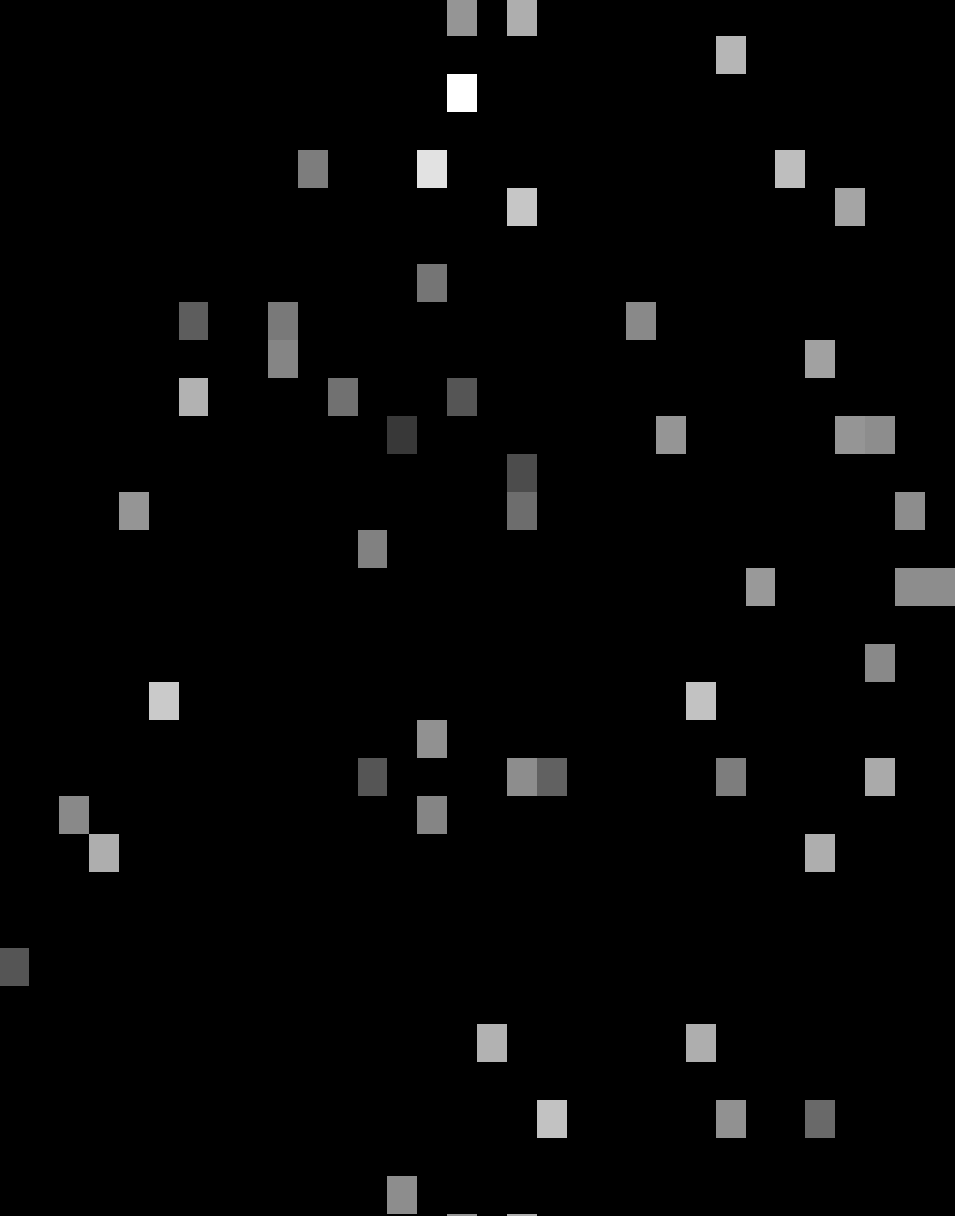}
		\includegraphics[width=.105\textwidth]{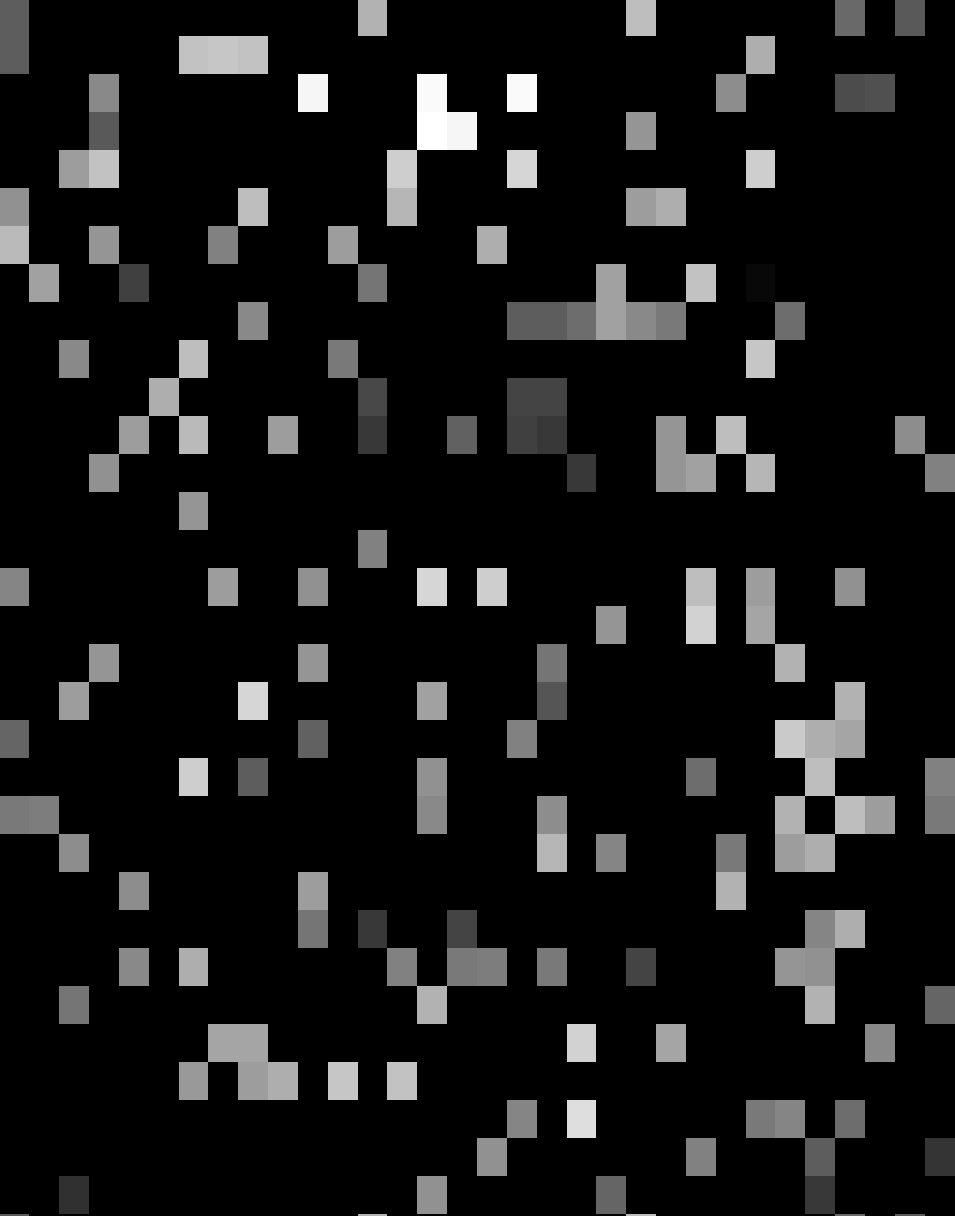} &~
		\includegraphics[width=.105\textwidth]{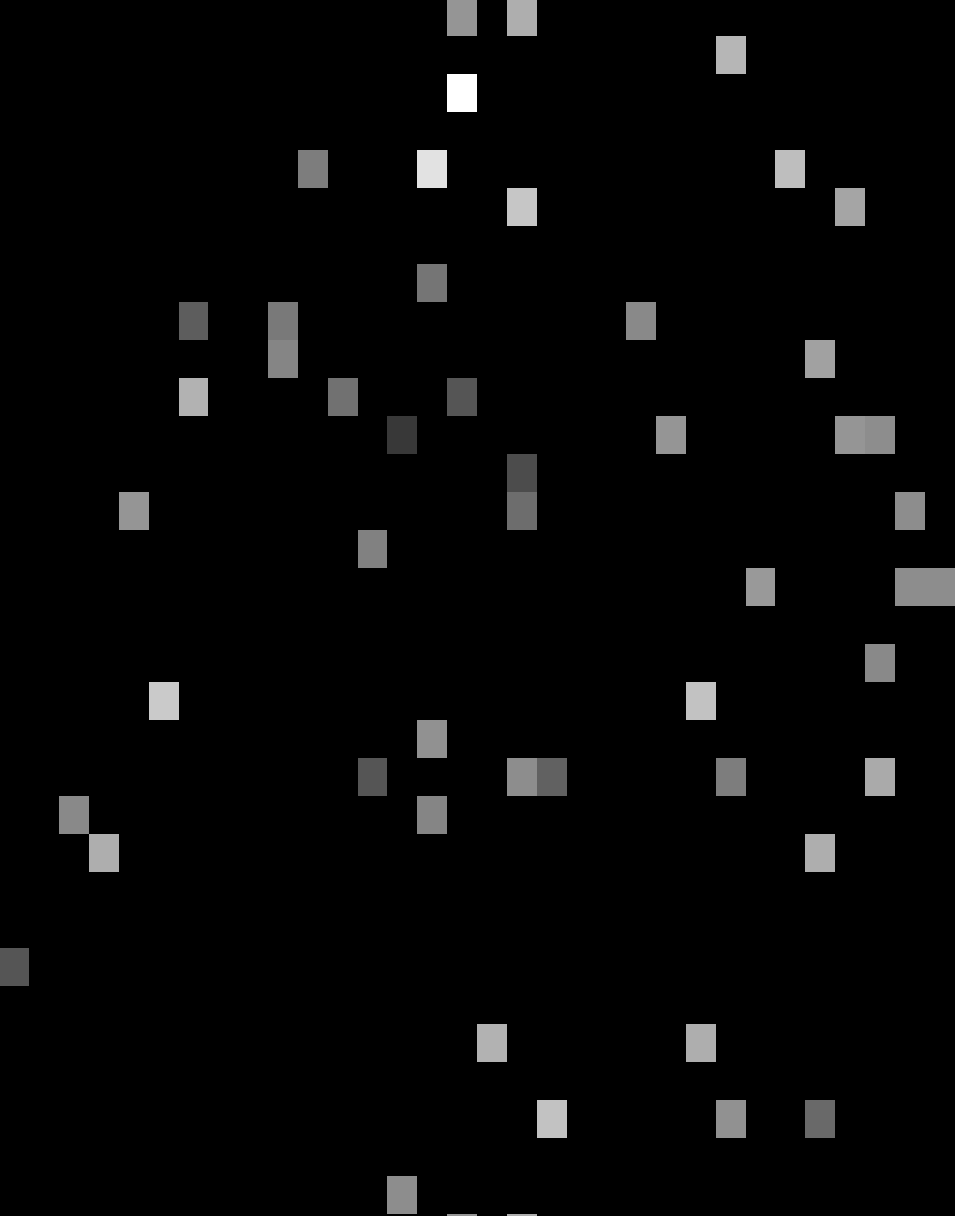} 
		\includegraphics[width=.105\textwidth]{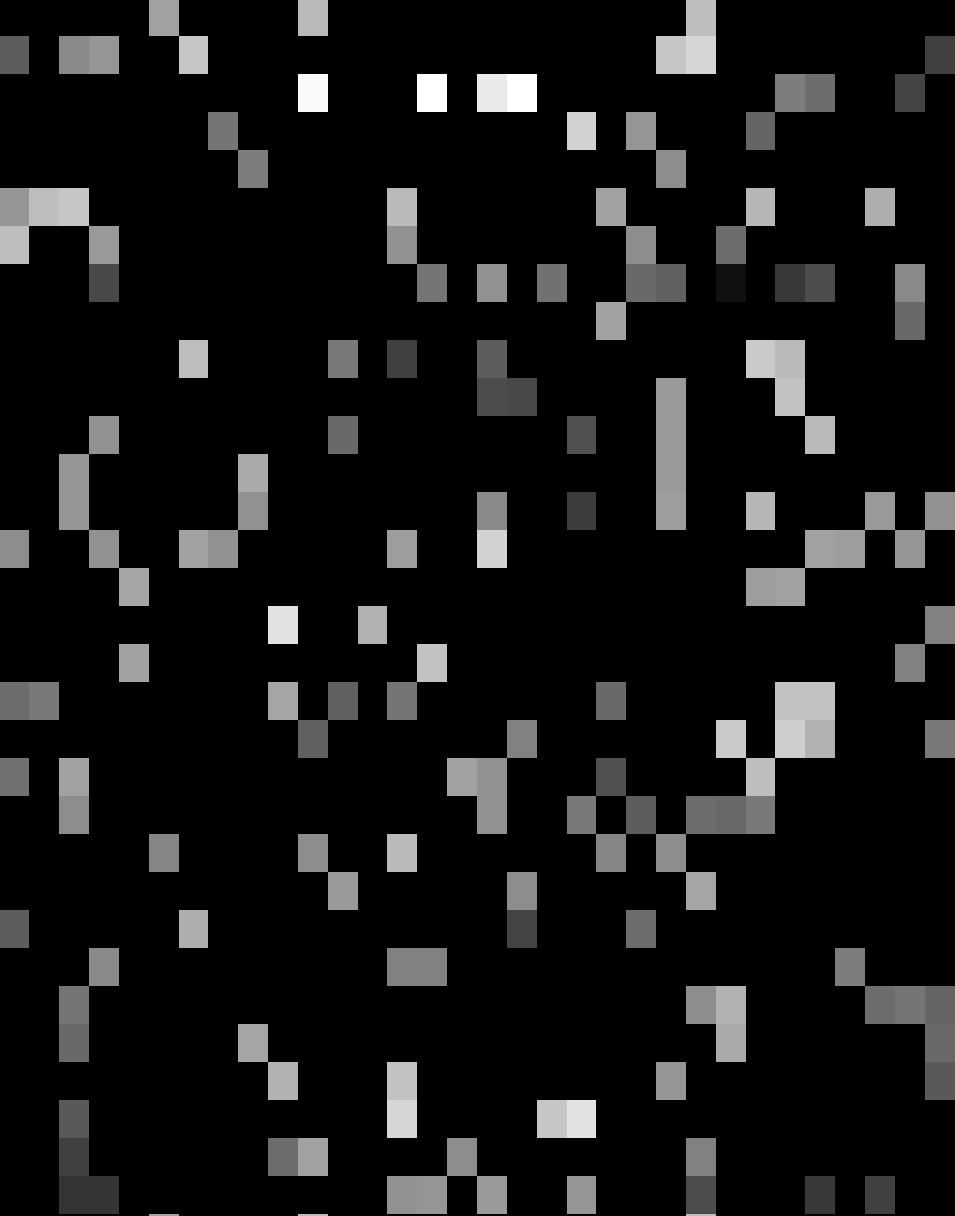}
		\includegraphics[width=.105\textwidth]{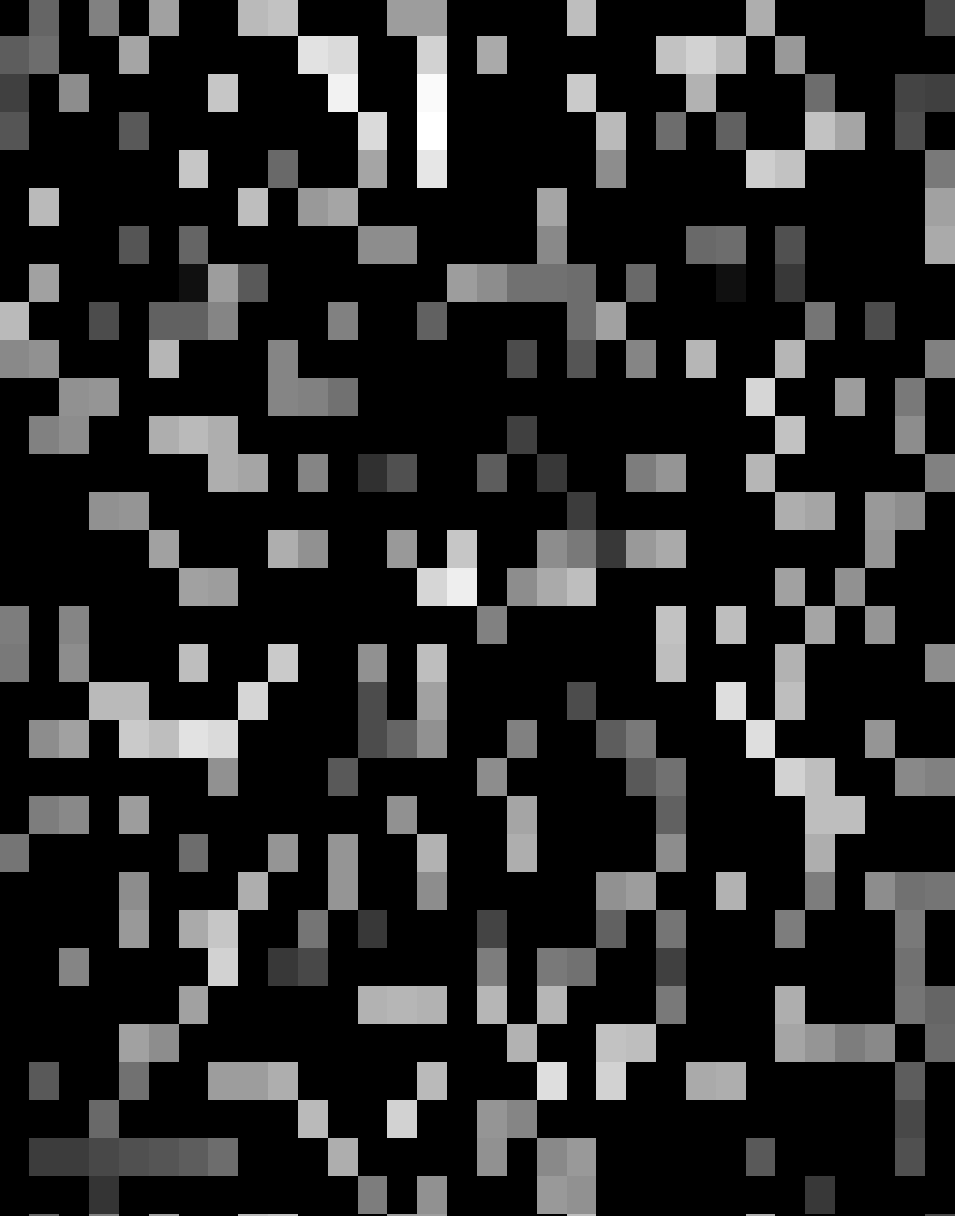}
		\includegraphics[width=.105\textwidth]{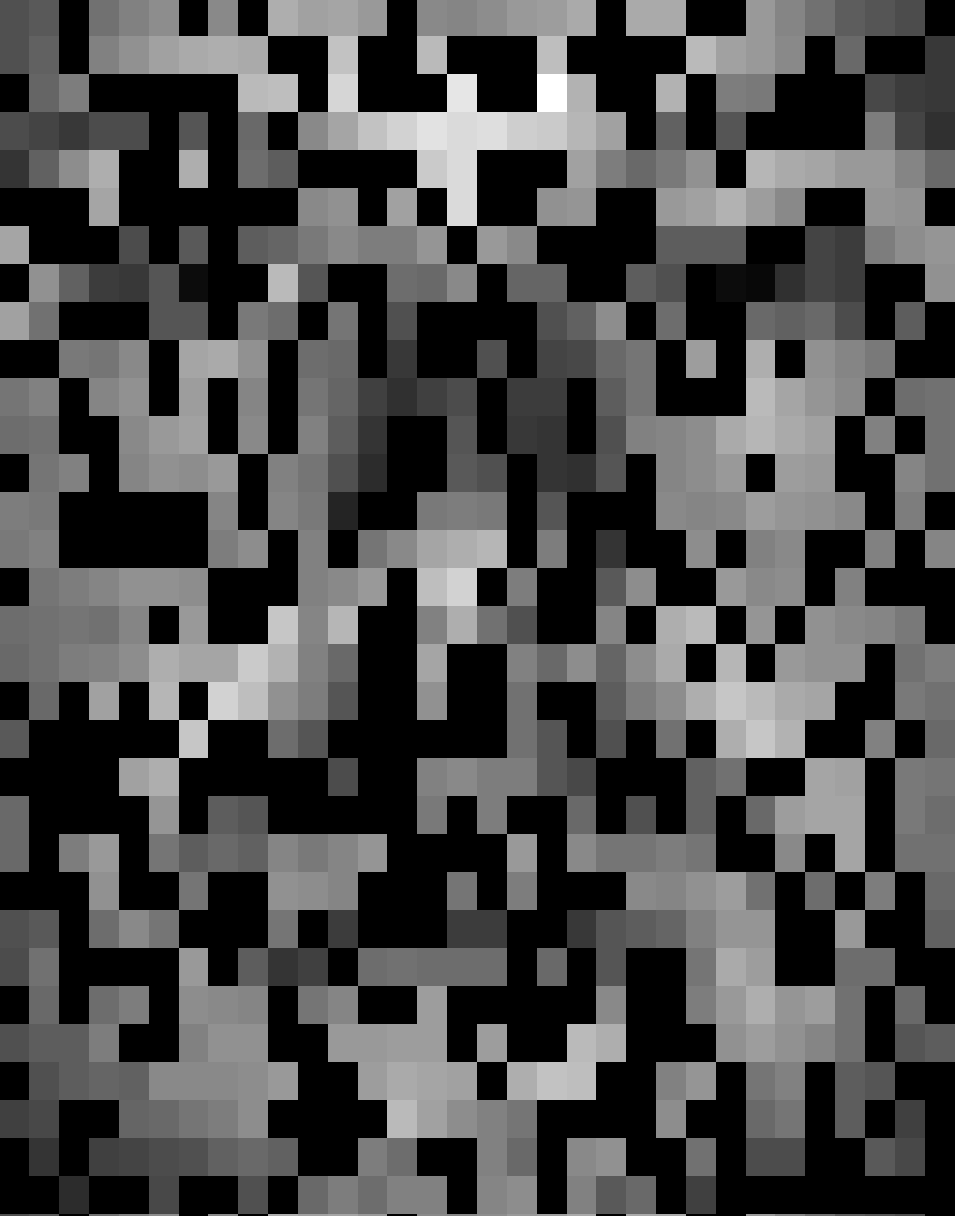} \\ 
		\begin{sideways} \bf Reconstruction \end{sideways}~
		\includegraphics[width=.105\textwidth]{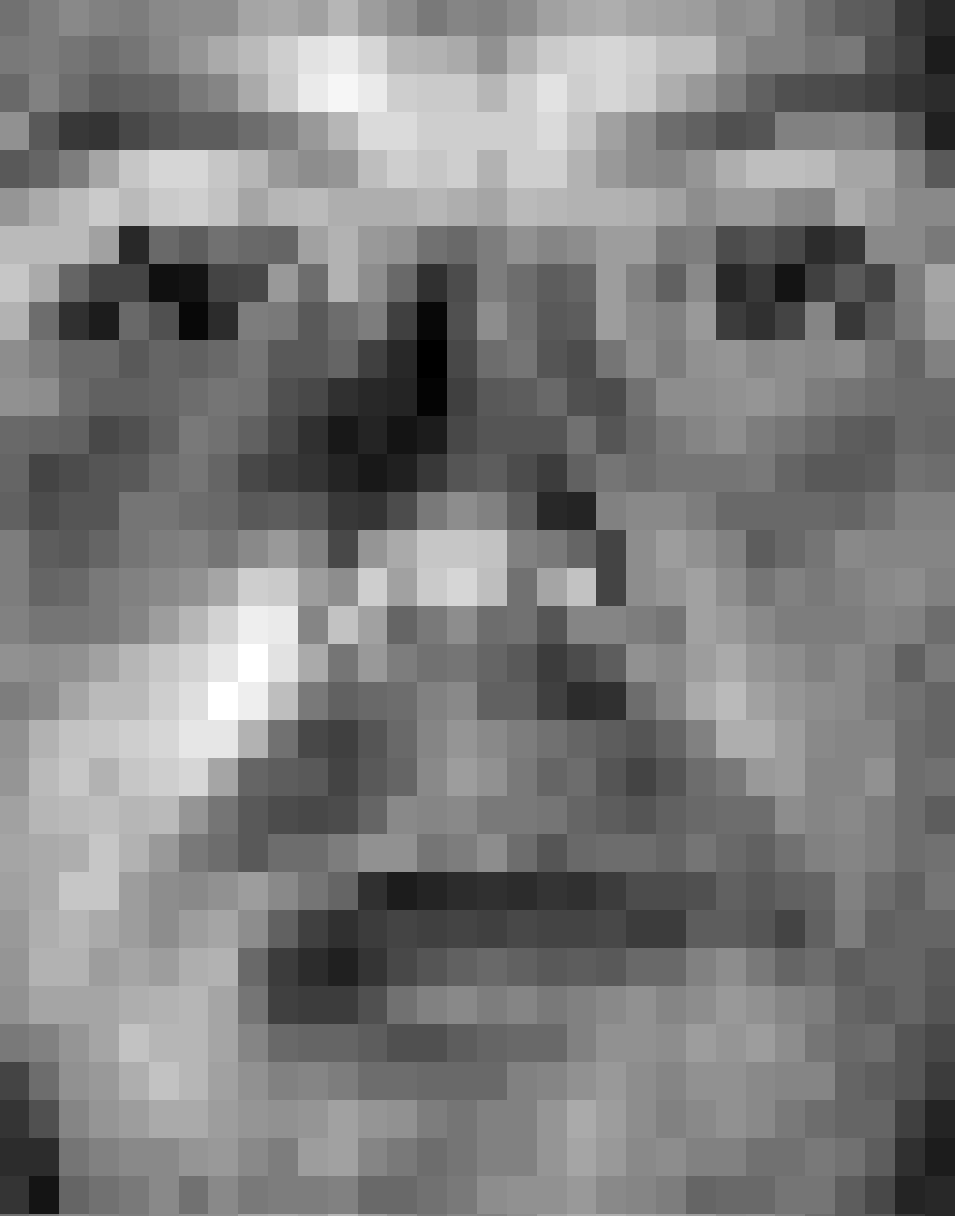}
		\includegraphics[width=.105\textwidth]{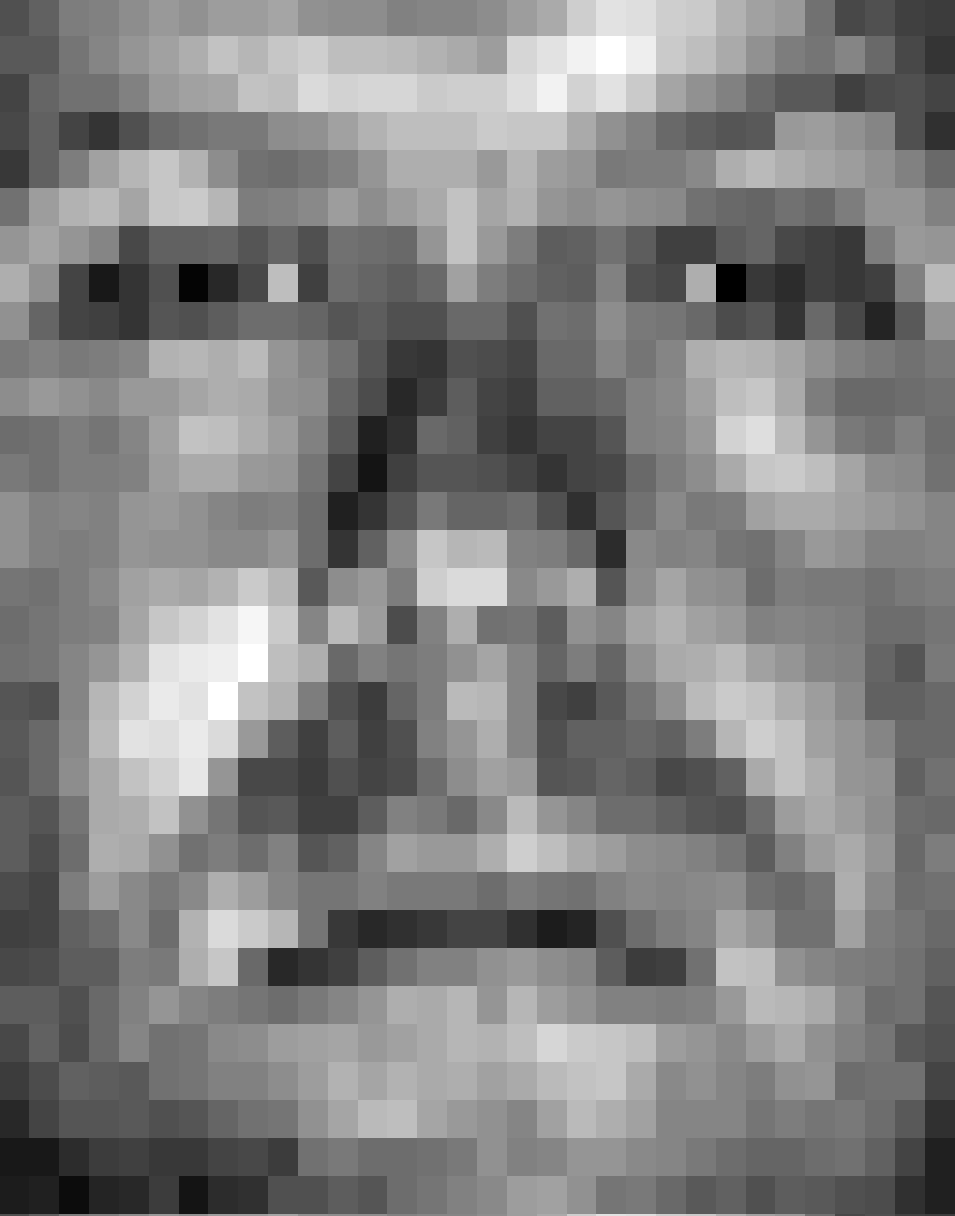} &~
		\includegraphics[width=.105\textwidth]{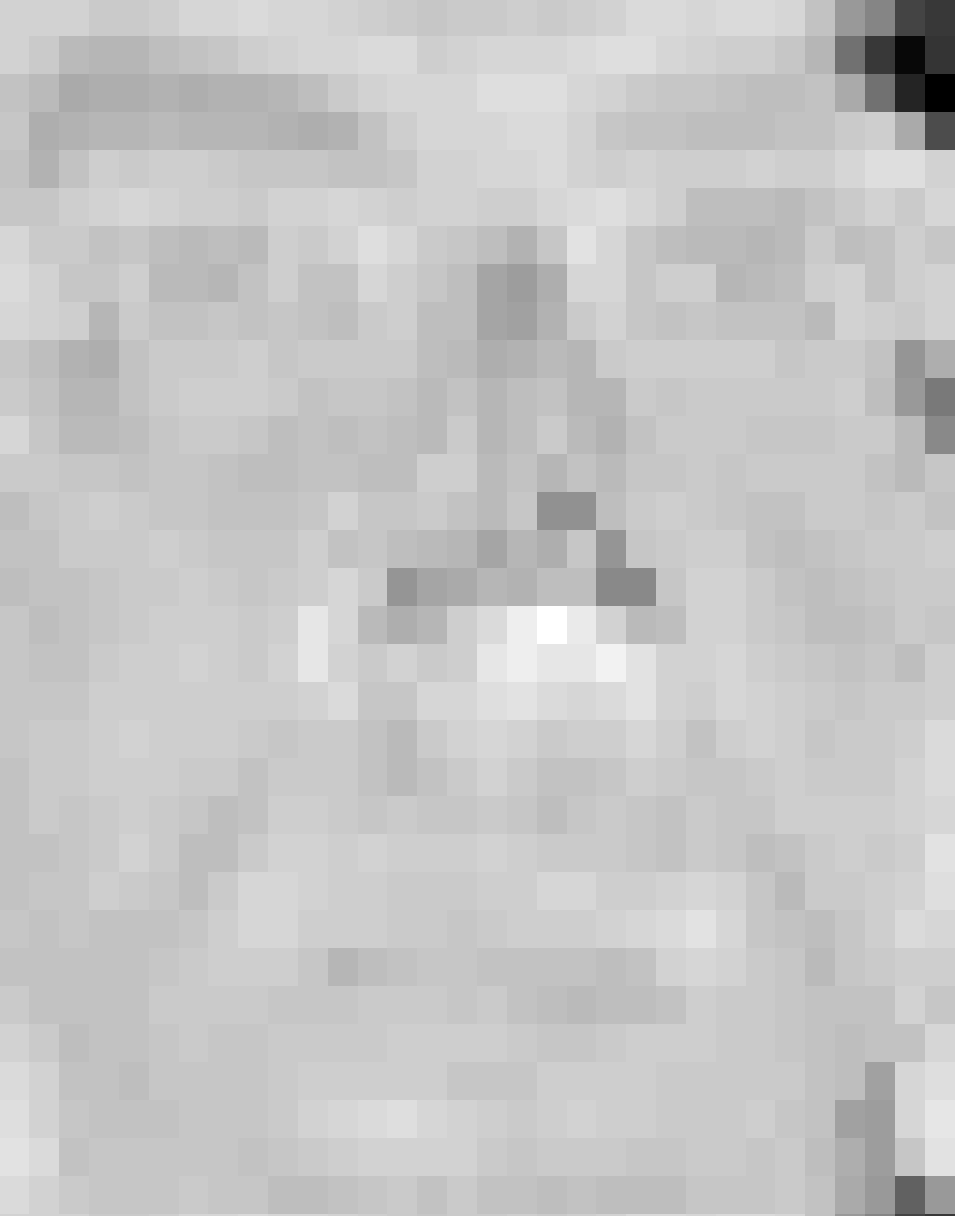} 
		\includegraphics[width=.105\textwidth]{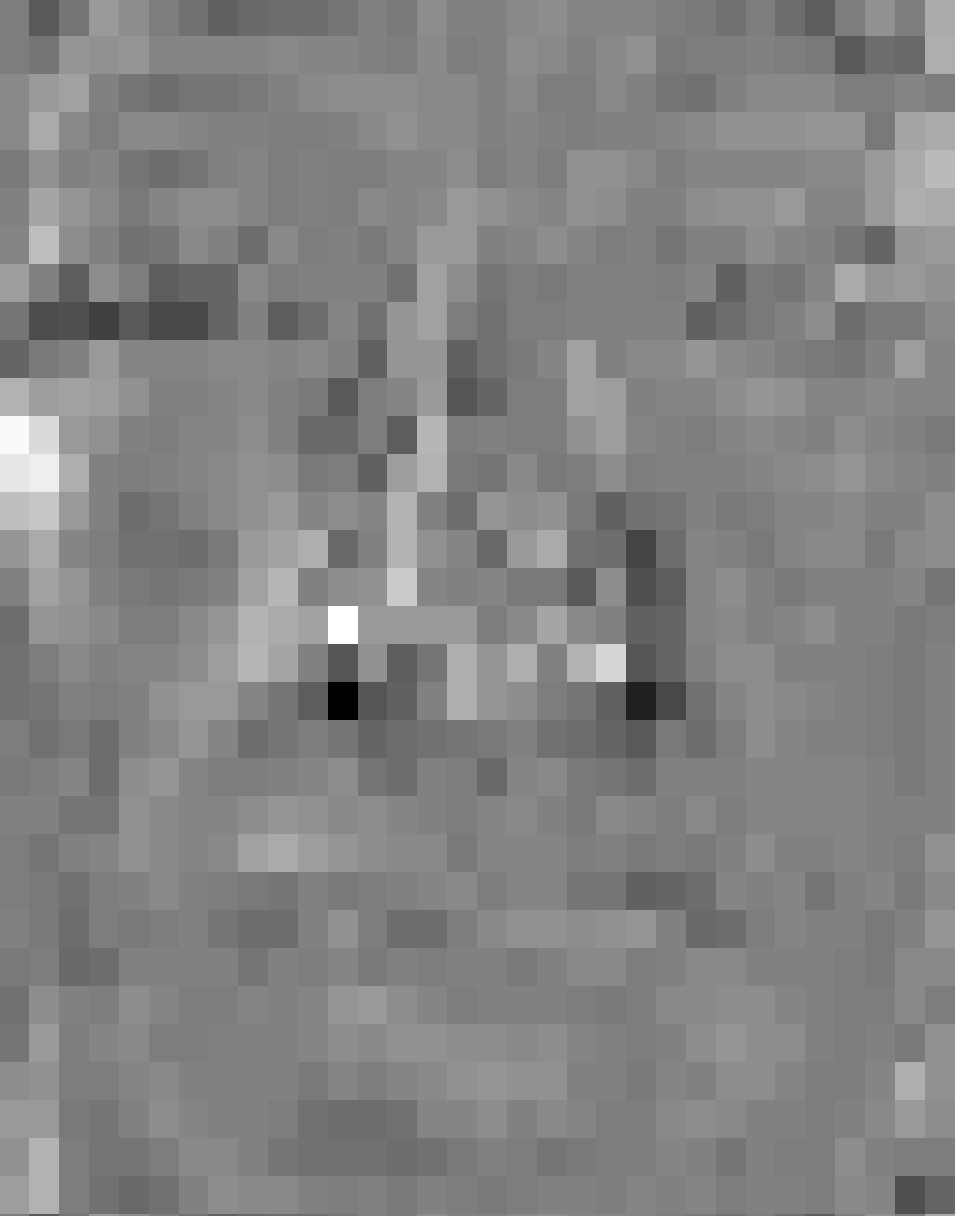} & ~
		\includegraphics[width=.105\textwidth]{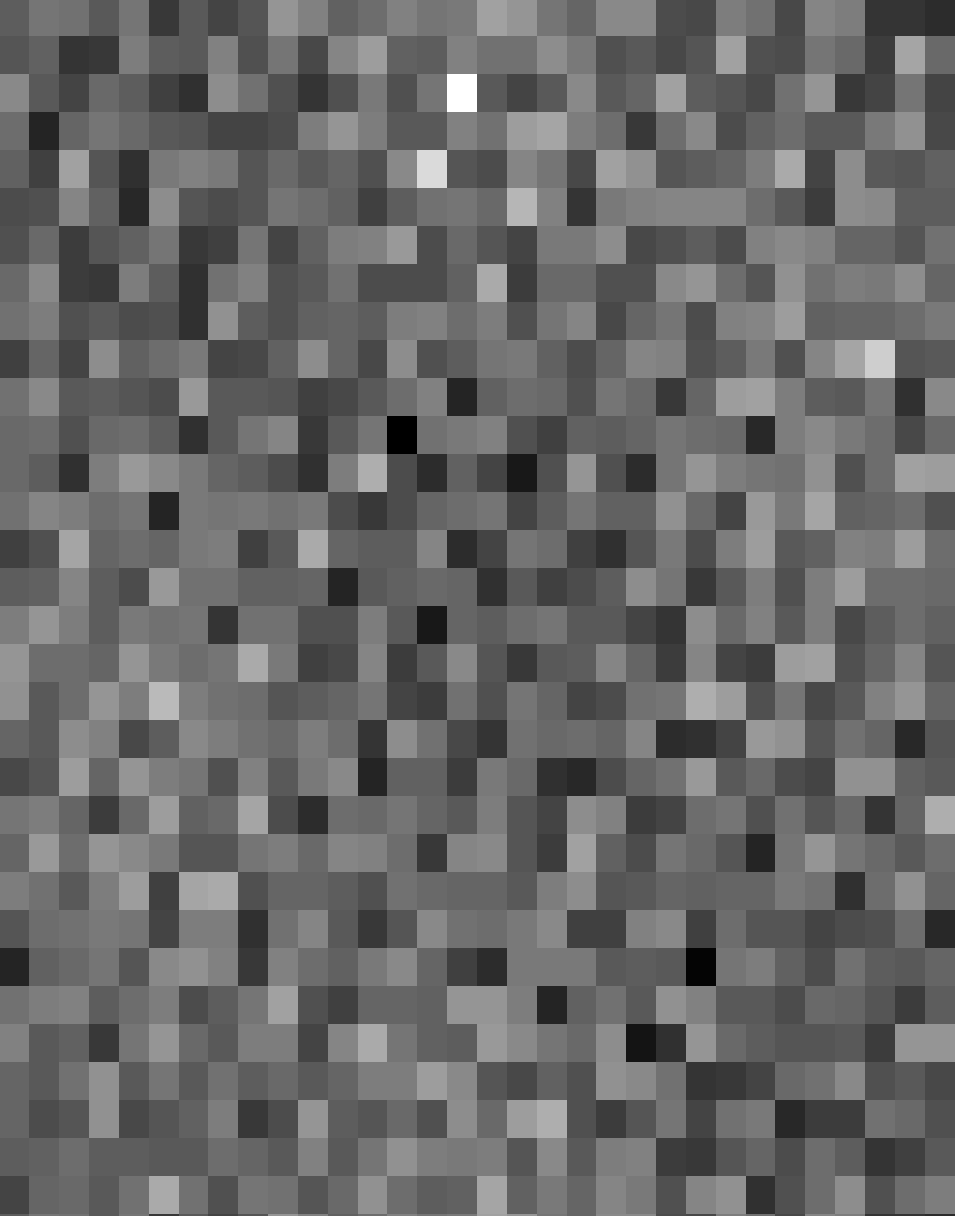}
		\includegraphics[width=.105\textwidth]{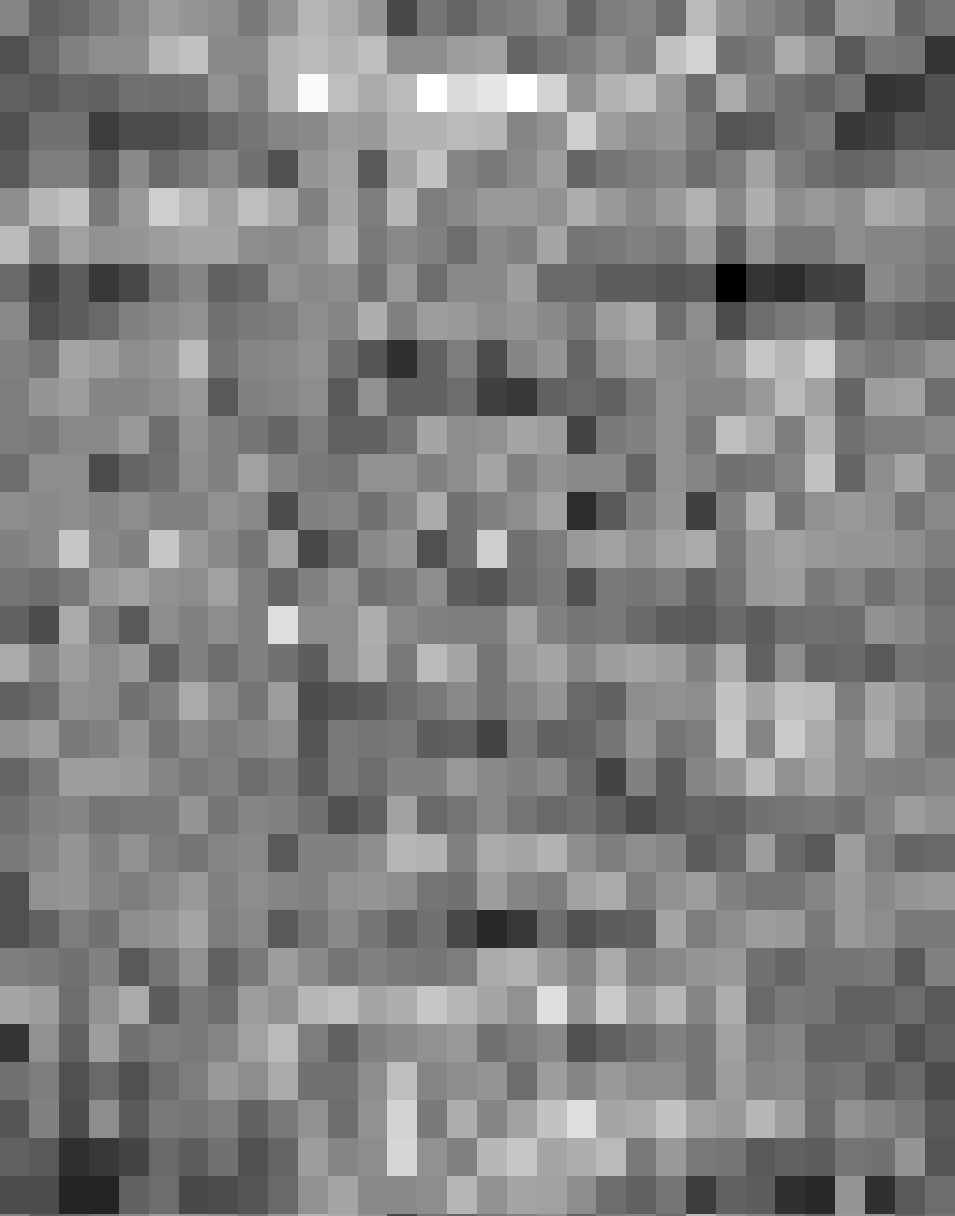}
		\includegraphics[width=.105\textwidth]{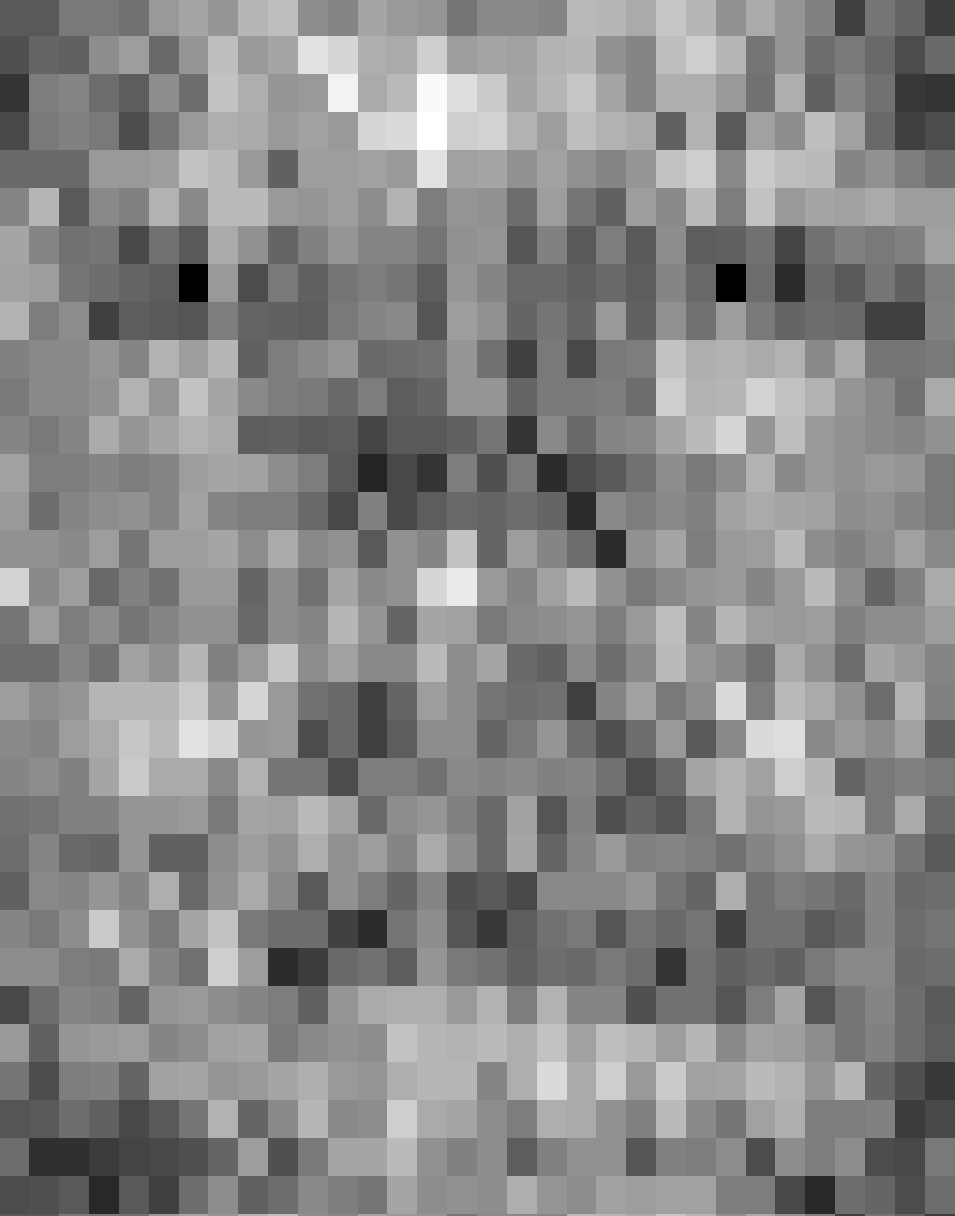}
		\includegraphics[width=.105\textwidth]{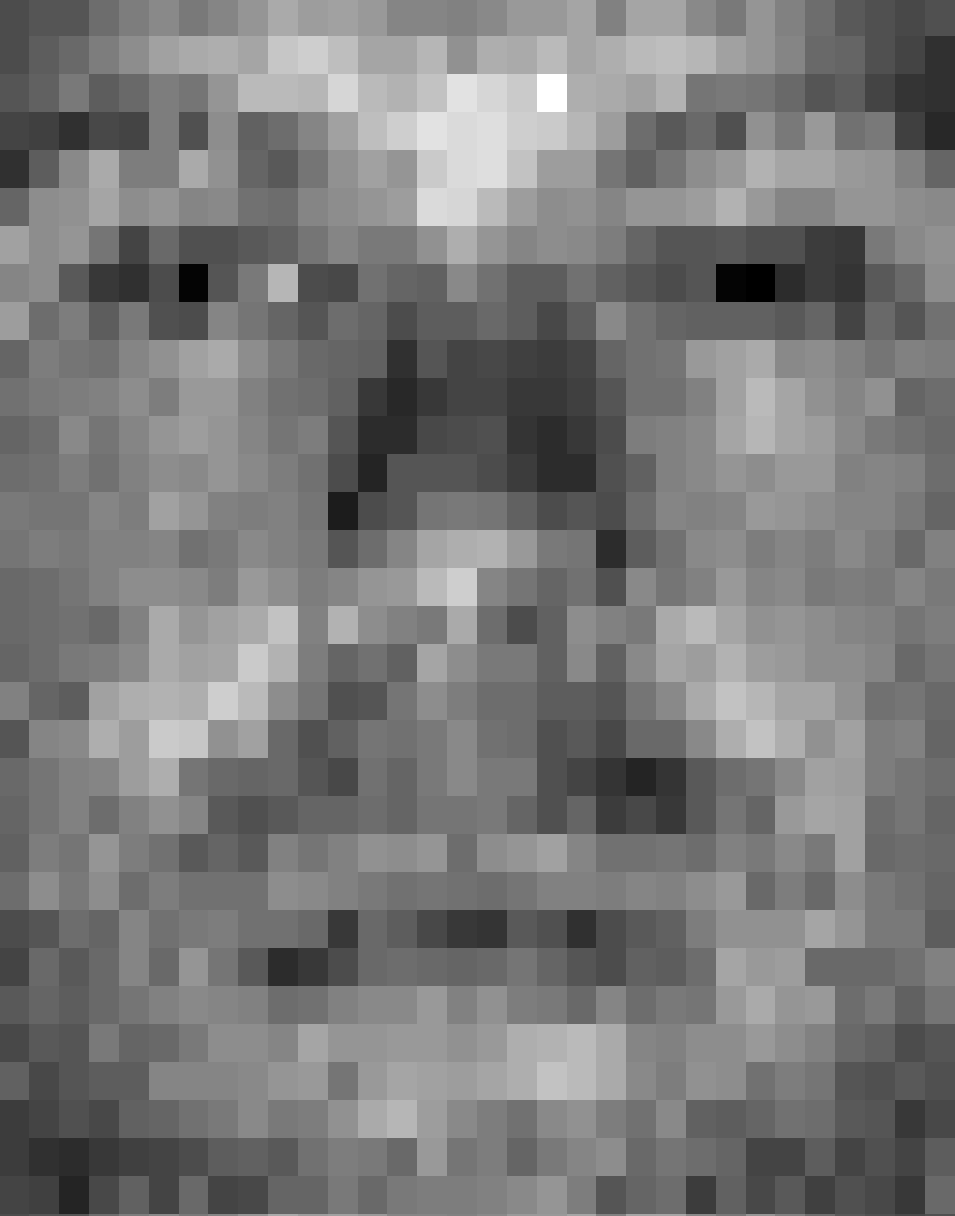} \\ \hline
	\end{tabular}
	\vspace{.05in}
	\caption{Comparison of reconstruction and sensor selection methods. QR sensors emphasize facial features such as the eyes, nose and mouth, and hence achieve adequate reconstruction with as few as 166 sensors (16\% of all pixels). Least squares ($\ell_2$) reconstruction with 50 QR sensors (5\% of pixels) surpasses the performance of compressed sensing with 300 random pixels. In comparison, compressed sensing requires 600 sampled pixels for comparable recovery, while $\ell_2$ reconstruction with random sensors is consistently poor. \label{fig:yalebrecon}}
\end{figure*}

Image processing and computer vision commonly involve high-resolution data with dimension determined by the number of pixels. Cameras and recording devices capture massive images with rapidly increasing pixel and temporal resolution. However, most pixel information in an image can be discarded for subject identification and automated decision-making tasks. 

The extended Yale B face database~\cite{YaleB2001ieee,YaleB2005ieee} is a canonical dataset used for facial recognition, and it is an ideal test bed for recovering low-rank structure from high-dimensional pixel space. 
The data consists of 64 aligned facial images each of 38 stationary individuals in different lighting conditions. 
We validate our sensor (pixel) selection strategy by recovering missing pixel data in a validation image using POD modes or {\em eigenfaces} trained on 32 randomly chosen images of each individual.

Normalized singular values are shown in Fig.~\ref{fig:yalebspectrum}, and the optimal singular value truncation threshold~\cite{Gavish2014ieeetit} occurs at $r=166$, indicating the intrinsic rank of the training dataset. 
Indeed, selected eigenfaces are also shown to reveal no meaningful facial structure beyond eigenface $166$. 
QR pixel selection is performed on the first 50 and first 166 eigenfaces, and selected pixels shown in Fig.~\ref{fig:yalebrecon} cluster around important facial features -- eyes, nose and mouth.
%
%
\begin{figure}
	\begin{overpic}[width=.5\textwidth]{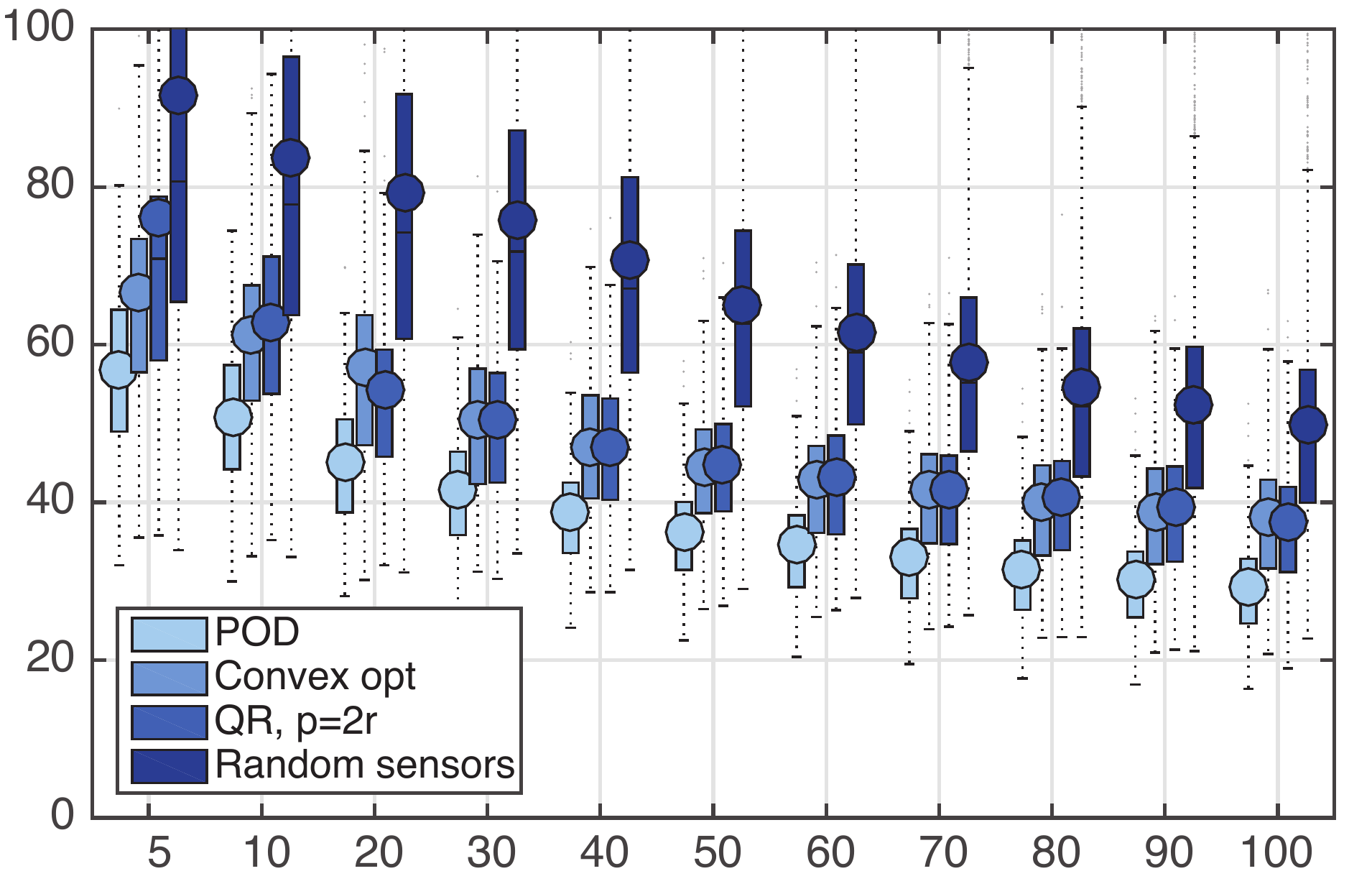}
		\put(40,-3){$r$ modes}
		\put(-2.5,20){\begin{sideways}Percent error\end{sideways}}
	\end{overpic}
	\vspace{.2in}
\caption{Recovery using our oversampled QR pivoting method for $p=2r>r$ sensors. QR pivoting sensors perform comparably to sensors obtained via the convex optimization (Convex opt) method~\cite{Joshi2009ieee}, and both are close to the optimal proper orthogonal decomposition approximation using the full state. In contrast to the convex method, QR sensors are obtained at reduced computational cost using the QR factorization of $\bPsi_r\bPsi_r^T$. Both methods outperform an equal number of randomly selected sensors. \label{fig:yale_conv}}
\end{figure}
\begin{figure*}[t]
	\centering
	\begin{overpic}[width=0.8\textwidth]{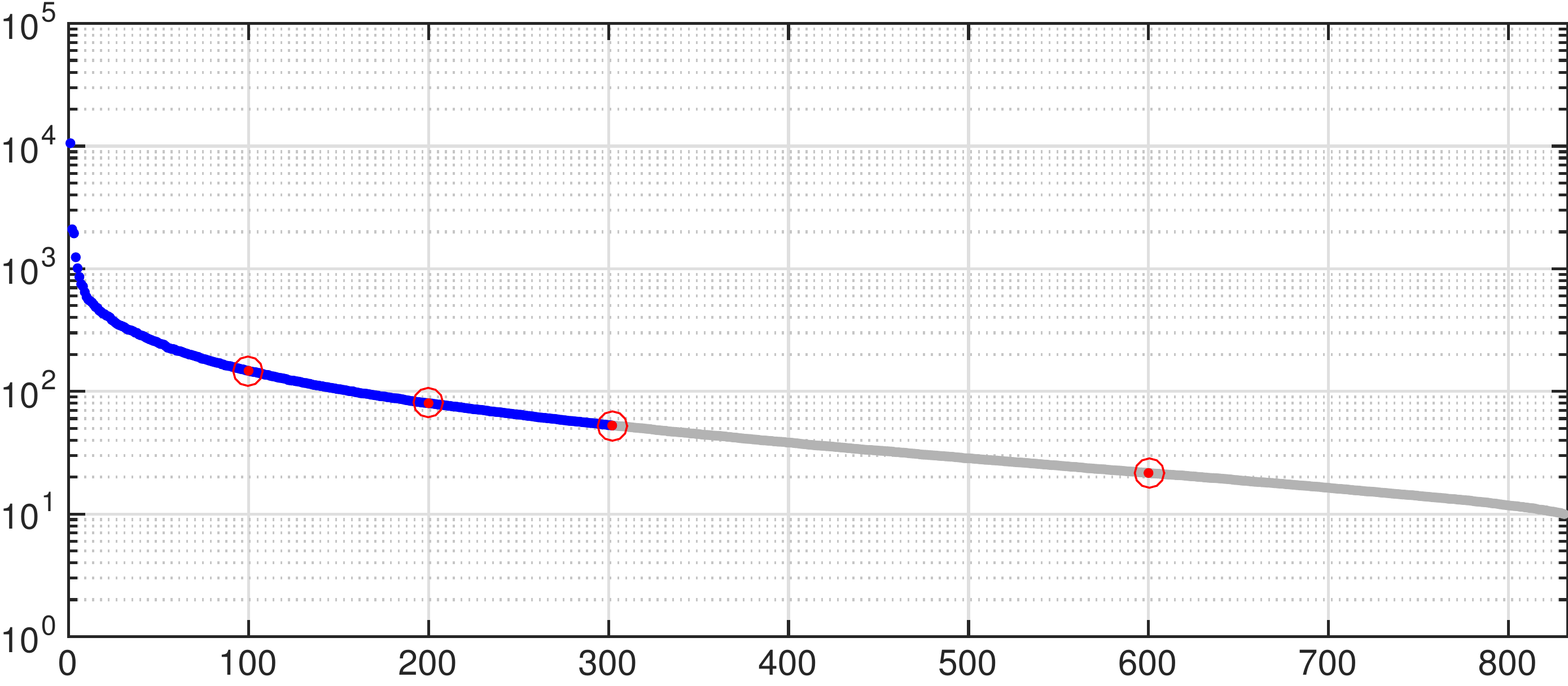}
		\setlength{\fboxsep}{0pt}%
		\put(5,28){\fbox{\includegraphics[width=0.125\textwidth]{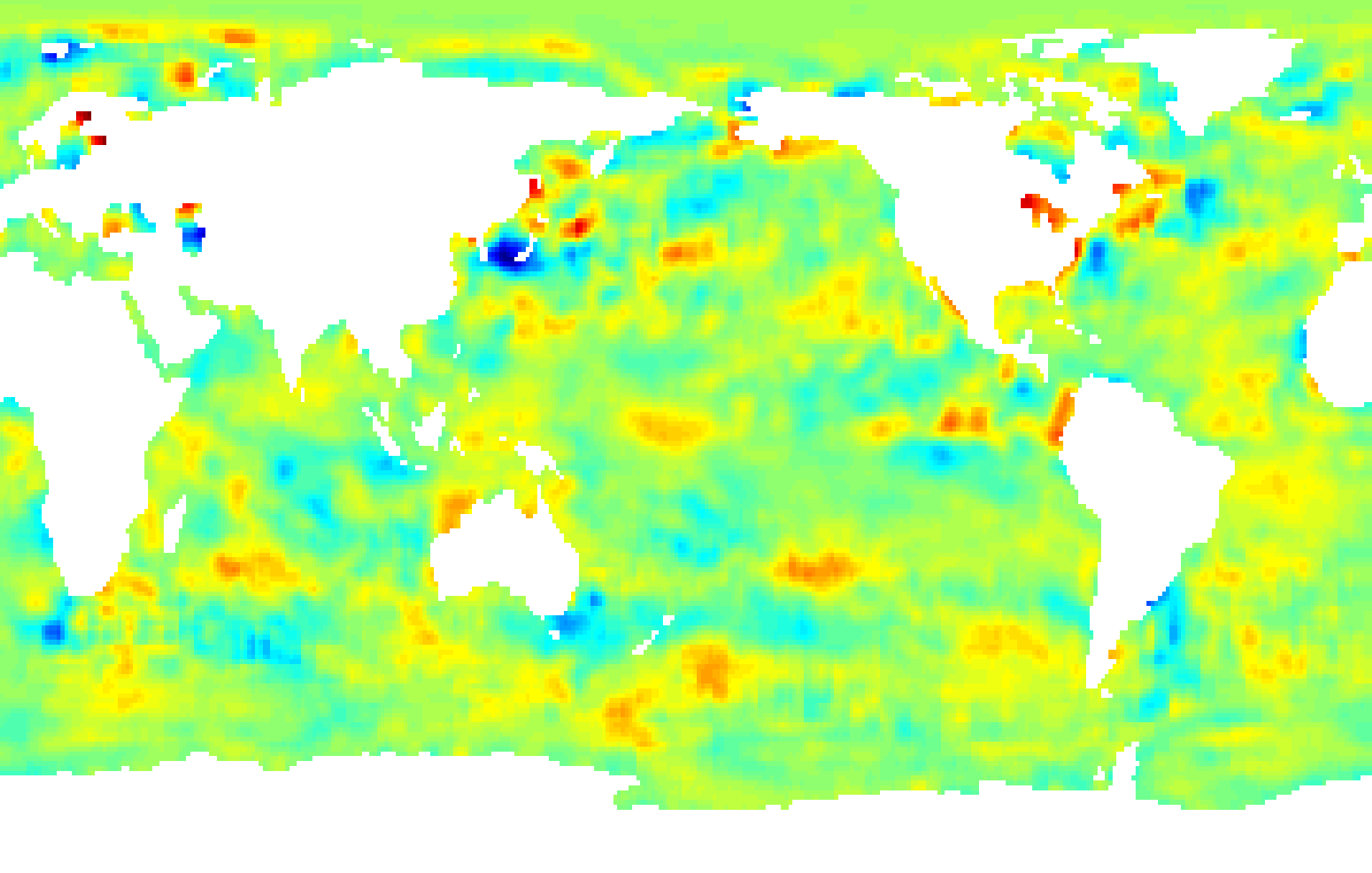}}}
		\put(7,39){Mode 100}
		\put(20,23){\fbox{\includegraphics[width=0.125\textwidth]{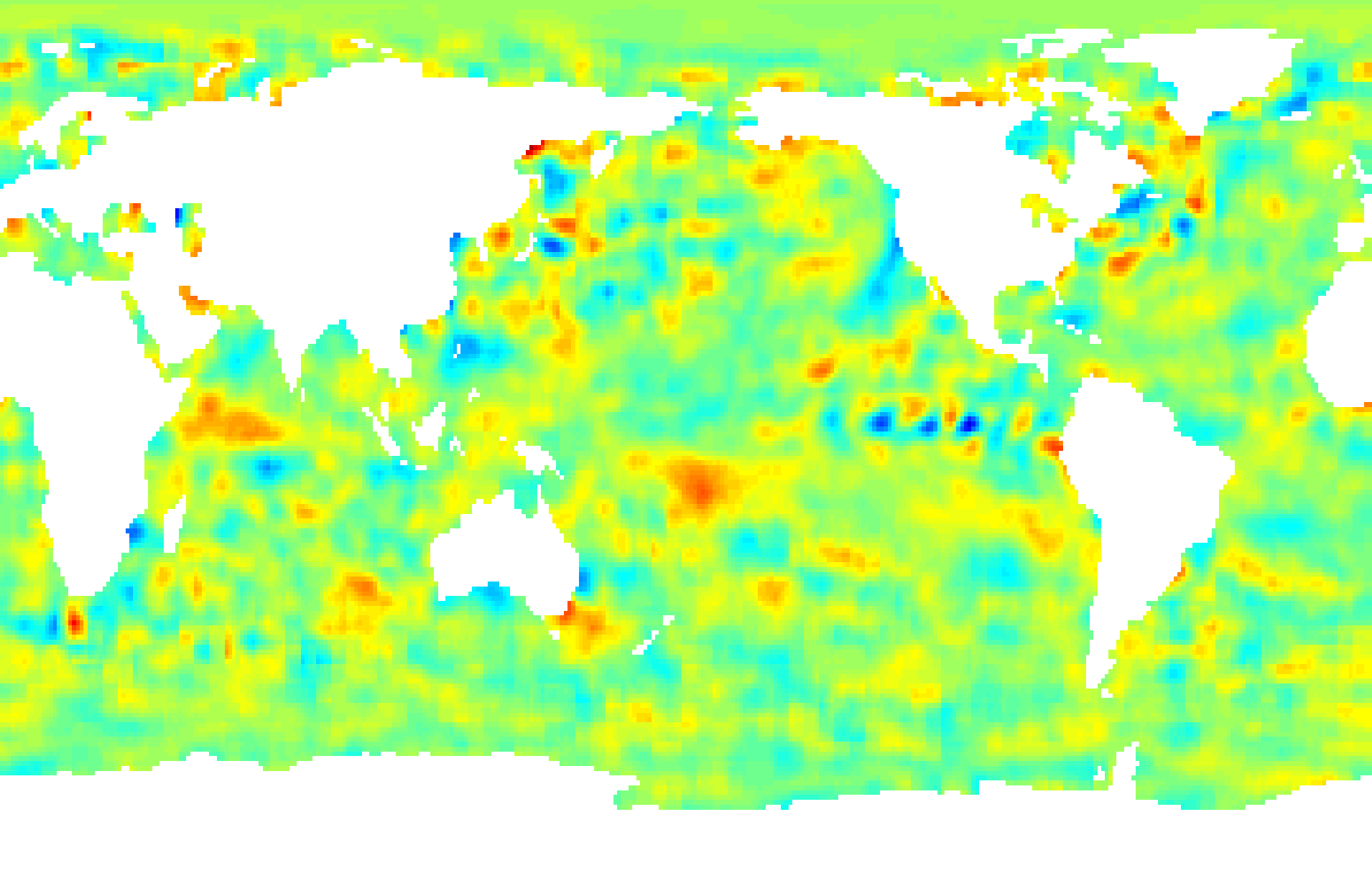}}}
		\put(22,34){Mode 200}
		\put(35,18){\fbox{\includegraphics[width=0.125\textwidth]{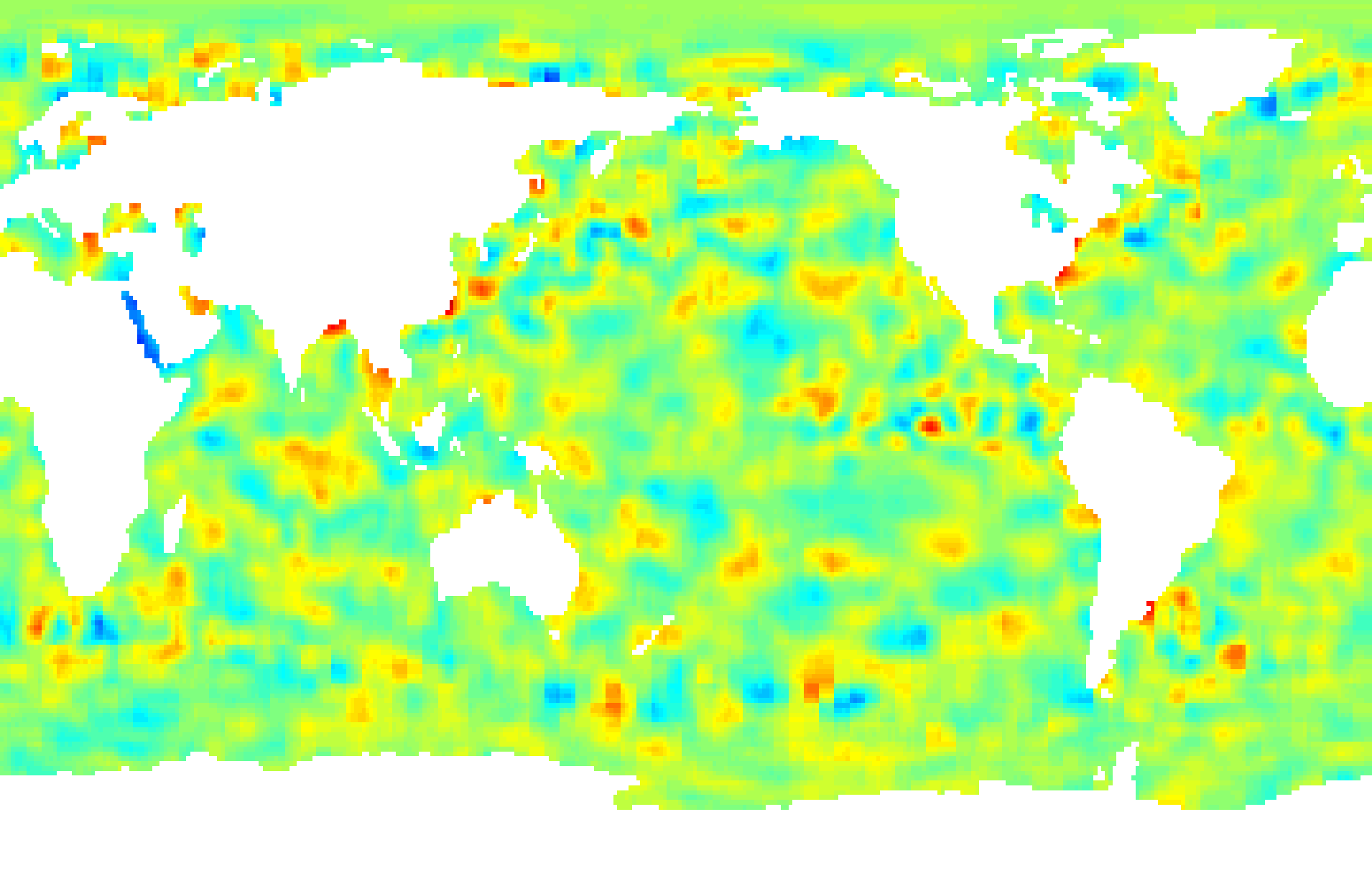}}}
		\put(37,29){Mode 302}
		\put(70,16){\fbox{\includegraphics[width=0.125\textwidth]{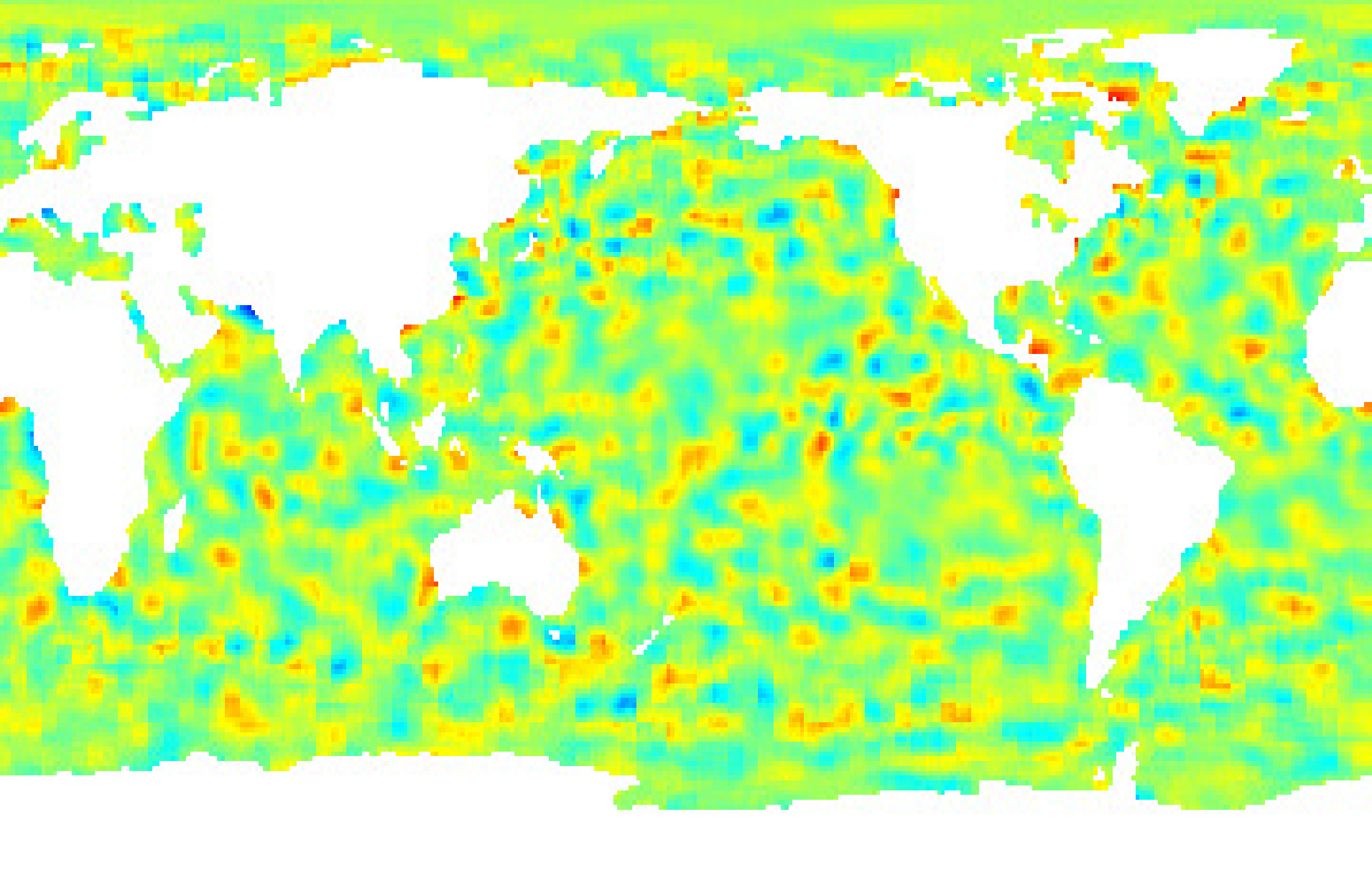}}}
		\put(73,27){Mode 600}	
		\put(50,-1.5){$r$}	
		\put(-2,22){\large$\sigma_r$}
	\end{overpic}
	\caption{Singular values and selected proper orthogonal decomposition (POD) modes for sea surface temperature data. The optimal rank truncation threshold occurs at $r=302$. POD mode 302 contains energetic localized convective phenomena (El Ni\~no) but is largely uninformative on a global scale. Thus, overfitting may occur as more modes are included (see Fig.~\ref{fig:enso_recon}). \label{fig:enso_spect}}
\end{figure*}

Image reconstructions in Fig.~\ref{fig:yalebrecon} are estimated from the same number of selected pixels as the number of modes used for reconstruction. For instance, the 50 eigenface reconstruction is uniquely constructed from 50 selected pixels out of 1024 total -- 5\% of available pixels. Even at lower pixel selection rates, least squares reconstruction from QR selected pixels is more successful at filling in missing data and recovering the subject's face.

For comparison, reconstruction of the same face from random pixels using compressed sensing is shown in Fig.~\ref{fig:yalebrecon}. Compressed sensing in a universal Fourier basis demonstrates progressively improved global feature recovery. However, more than triple the pixels are required for the same quality of reconstruction as in QR selection. Moreover, the convex $\ell_1$ optimization procedure is extremely expensive compared to the single $\ell_2$ regression on subsampled eigenfaces. Therefore data-driven feature selection and structured measurement selection are of significant computational and predictive benefit, and occur at the small training expense of one SVD and pivoted QR operation.

\begin{figure*}[t]
	\setlength{\fboxsep}{0pt}%
	\centering
	\begin{tabular}{| cccc|}
		\hline & {\bf 100 sensors} & {\bf 200 sensors} & {\bf \qquad\quad 302 sensors} \\
		\begin{sideways}{\qquad\quad\bf QR }\end{sideways}& \fbox{\includegraphics[width=.25\textwidth]{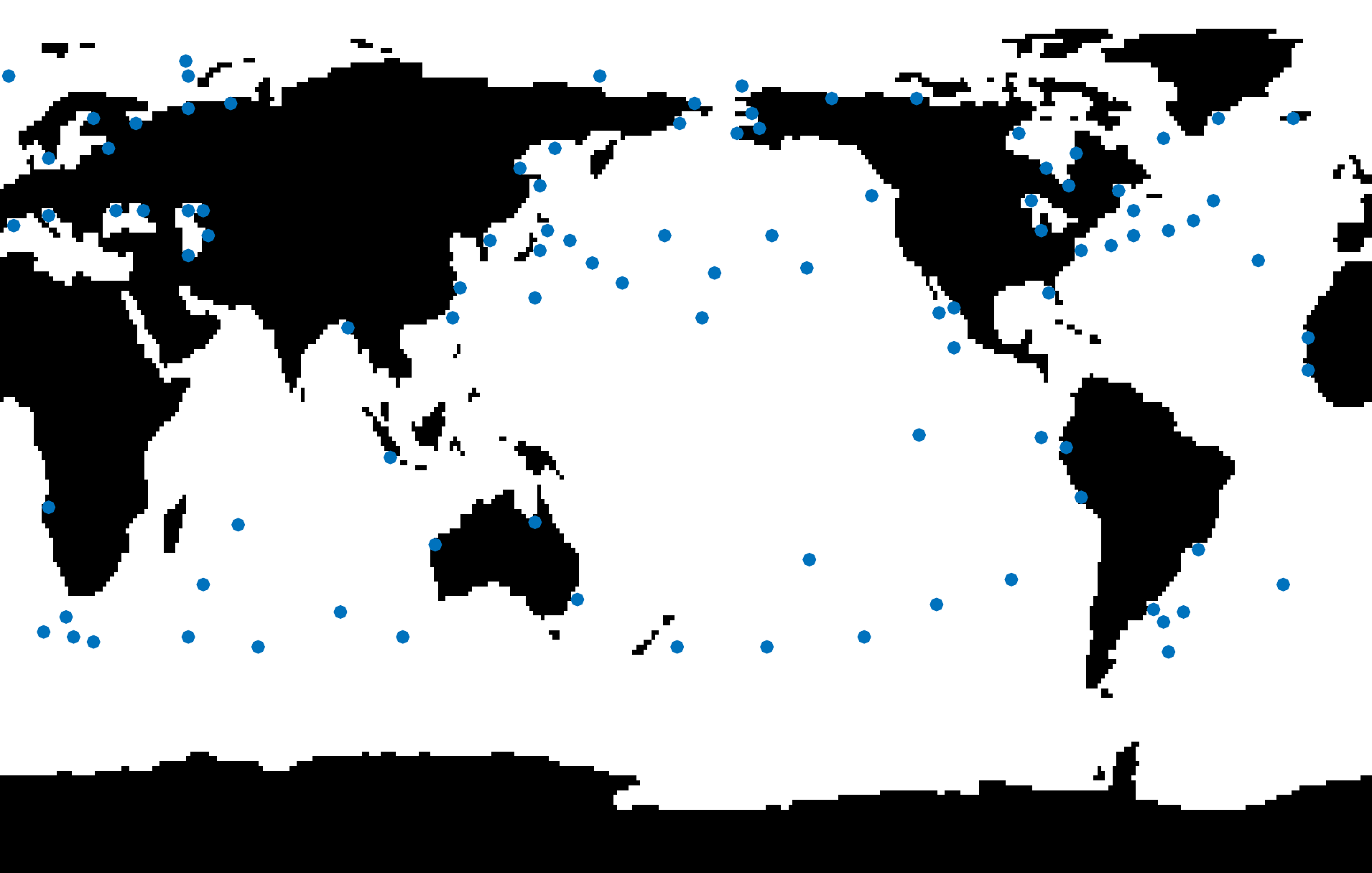}} & \fbox{\includegraphics[width=.25\textwidth]{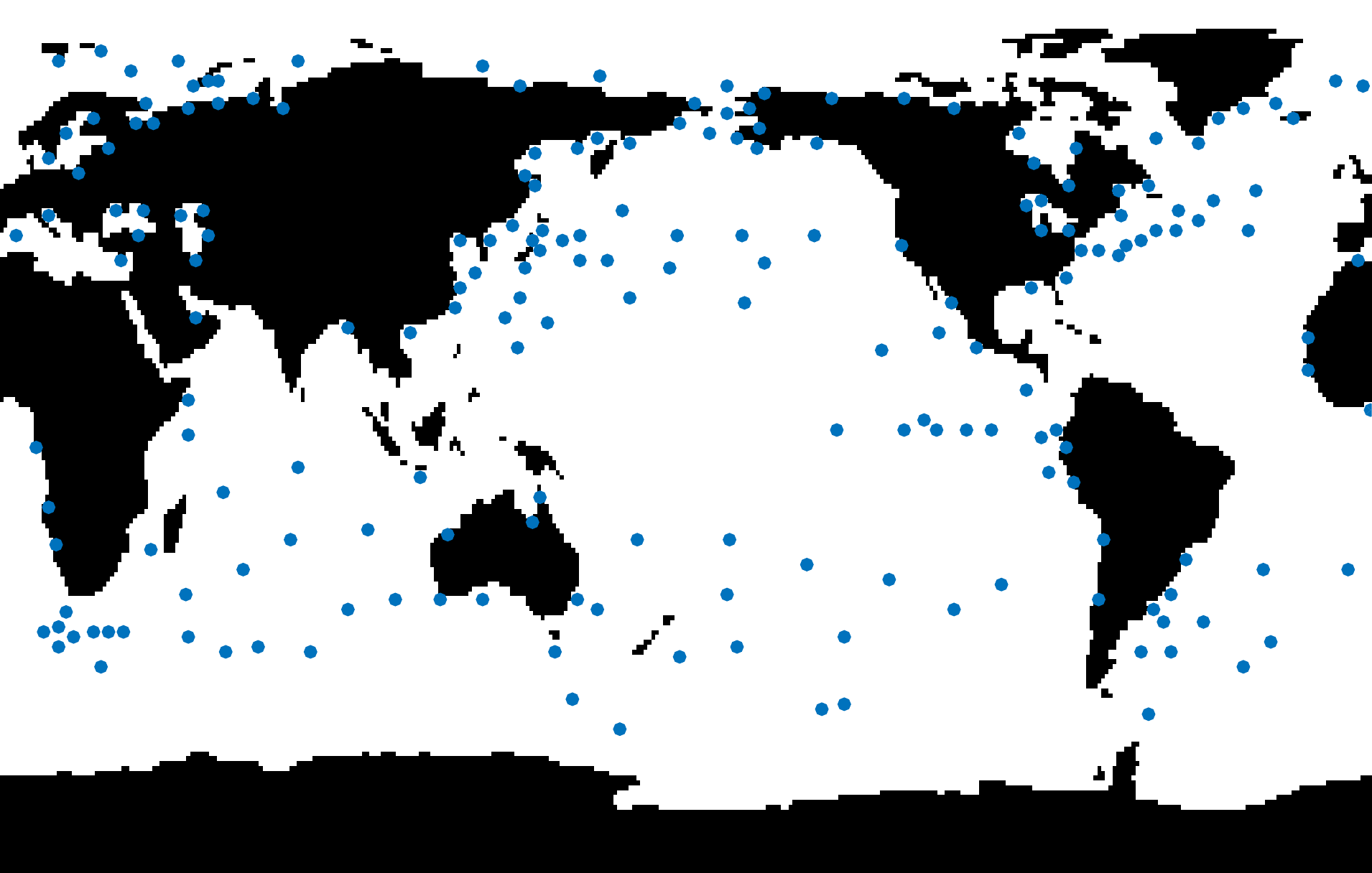}} & \fbox{\includegraphics[width=.25\textwidth]{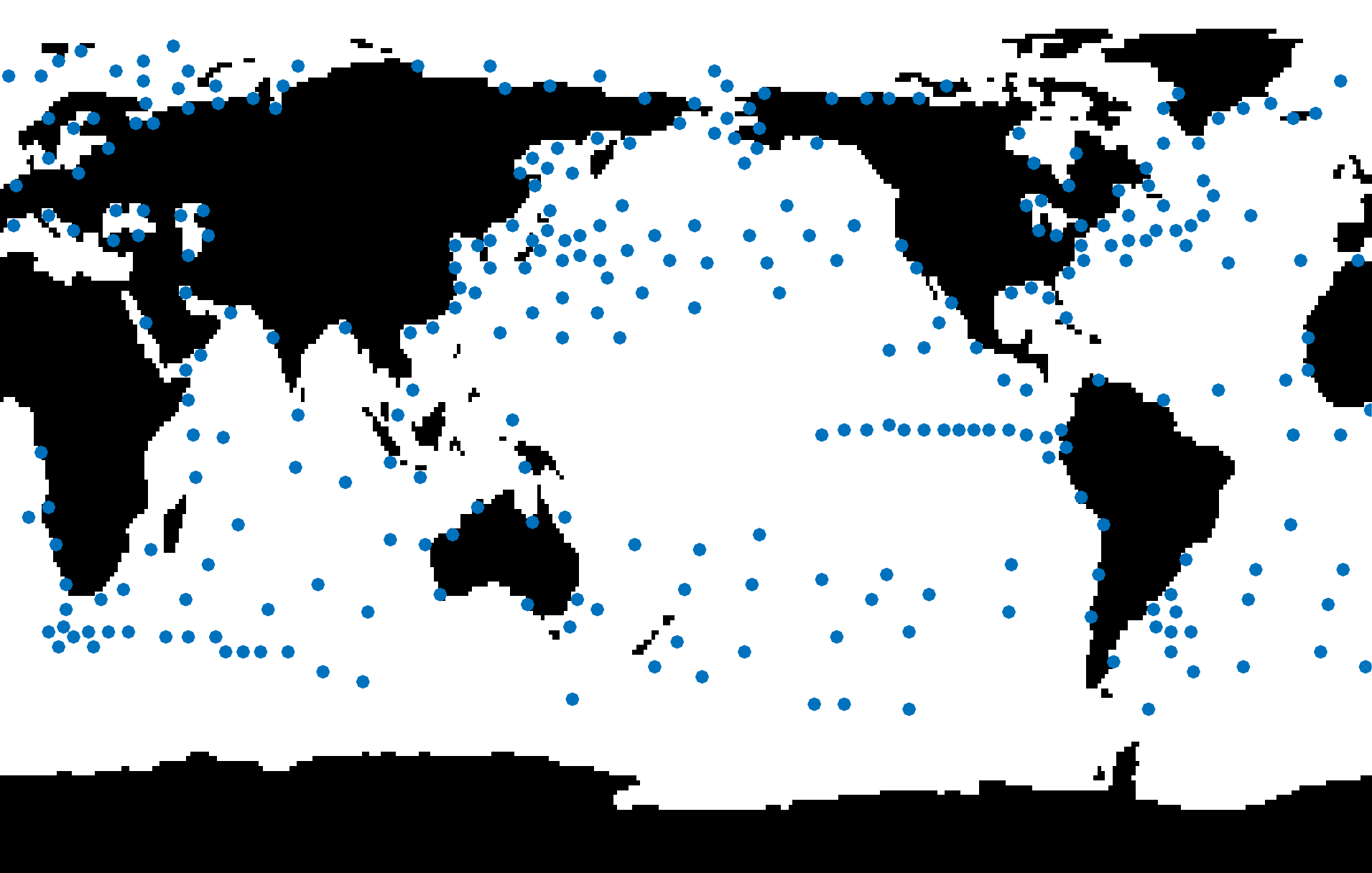}} \\
		& & & \\
		\begin{sideways}{\qquad\bf Random }\end{sideways} & \fbox{\includegraphics[width=.25\textwidth]{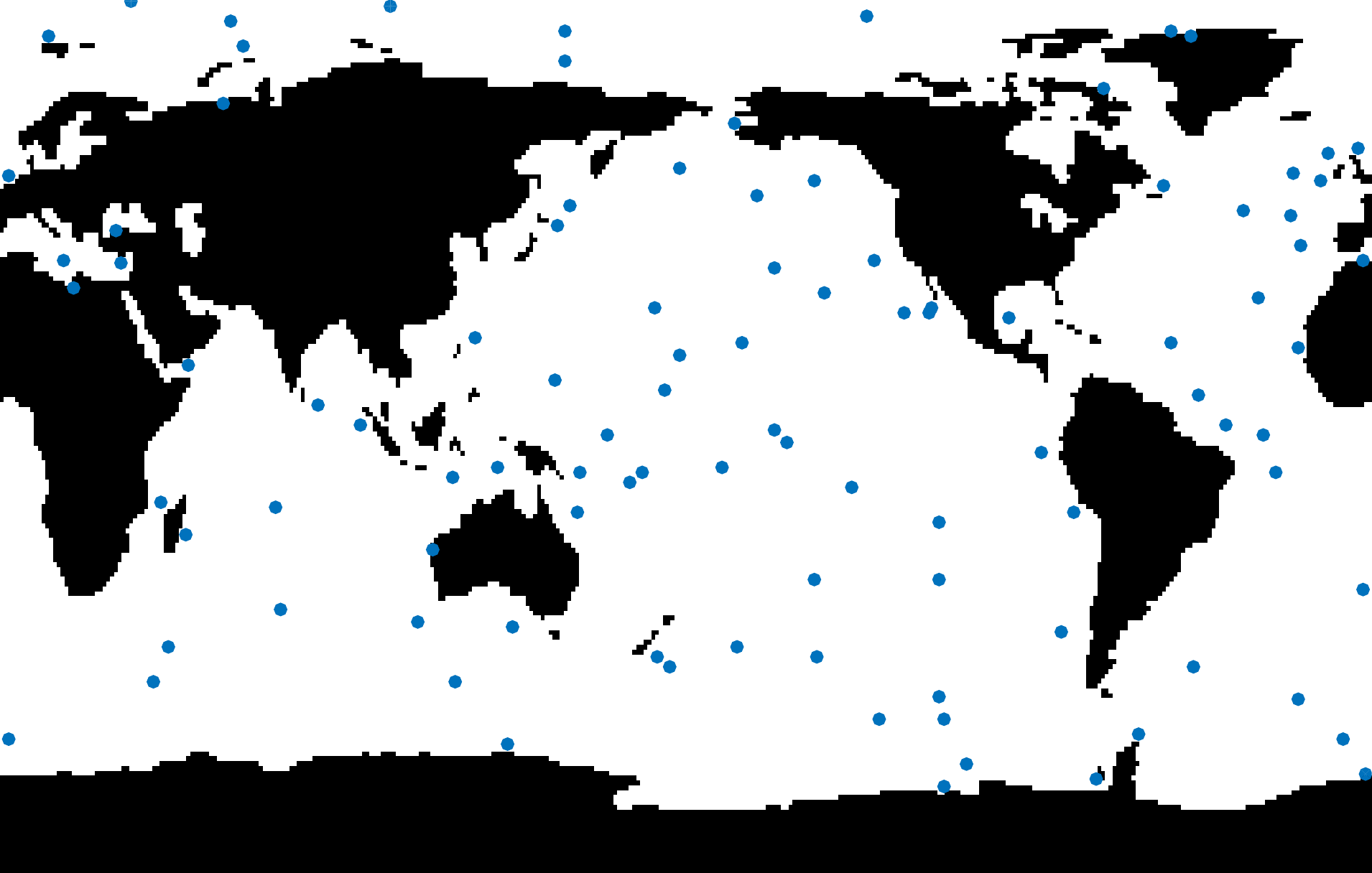}} & \fbox{\includegraphics[width=.25\textwidth]{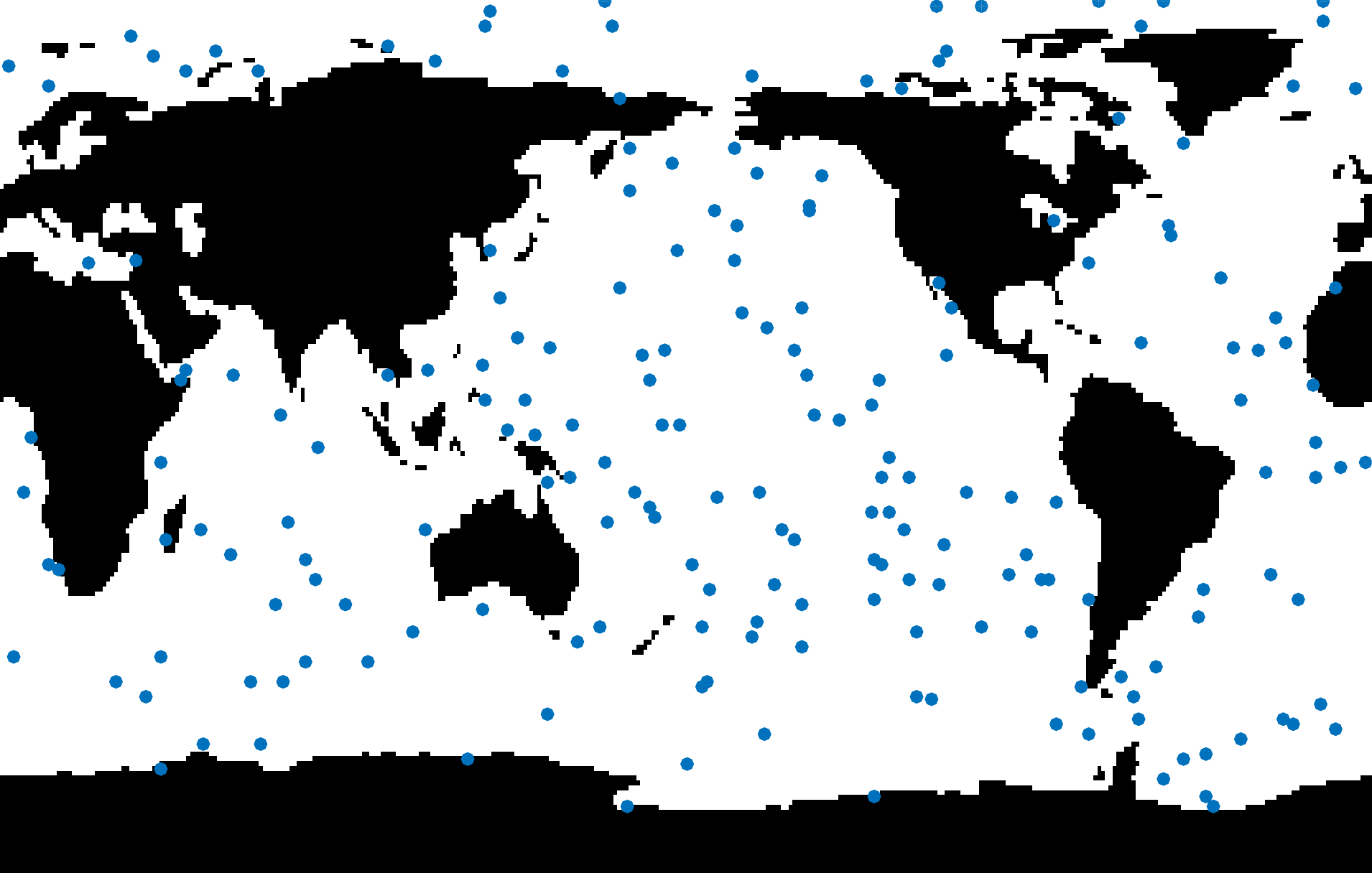}} & \fbox{\includegraphics[width=.25\textwidth]{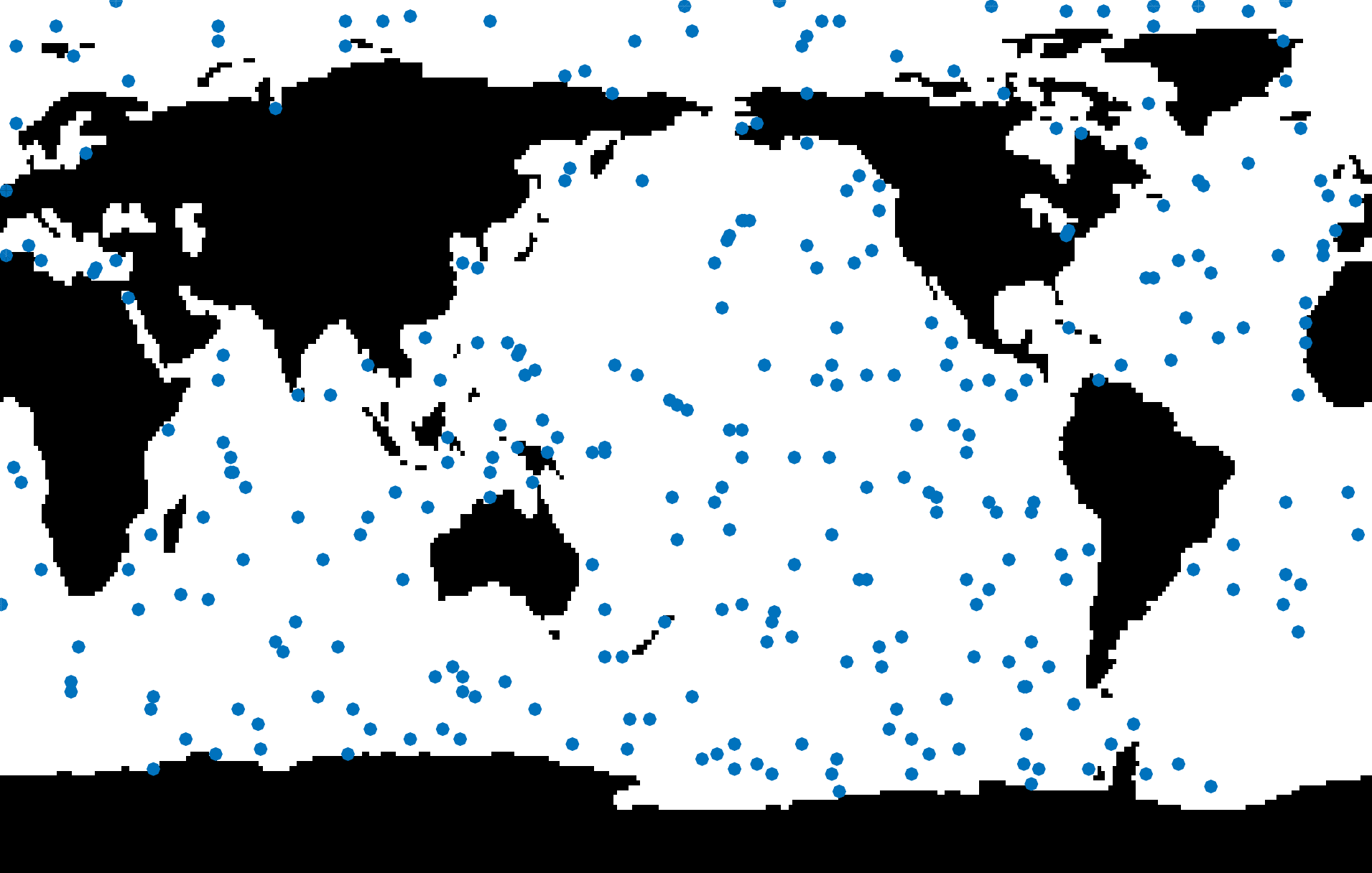}} \\ \hline
	\end{tabular}
	\caption{QR selected sensors used for reconstruction. QR sensors are informative about ocean dynamics, for example capturing convective phenomena such as the El Ni\~no Southern Oscillation off coastal Peru. \label{fig:enso_sensors}}  
\end{figure*}

The convergence of reconstruction with sensors using QR pivoting is shown in Fig.~\ref{fig:yale_conv}. More sensors than modes are used in reconstruction for this example. The expected error dropoff is observed with increasing number of modes and sensors, although the dropoff is slower than for the cylinder flow (Fig.~\ref{fig:ibpm_conv}) due to slower decay of singular values.


\subsection{Sea surface temperature (SST)}

\begin{figure*}[t]
	\setlength{\fboxsep}{0pt}%
	\centering
	\begin{tabular}{|rccc|}
		\hline  & {\bf 100 mode approximation} & {\bf 200 mode approximation} & {\bf 302 mode approximation} \\
		\begin{sideways}{\qquad\bf POD }\end{sideways}& \fbox{\includegraphics[width=.25\textwidth]{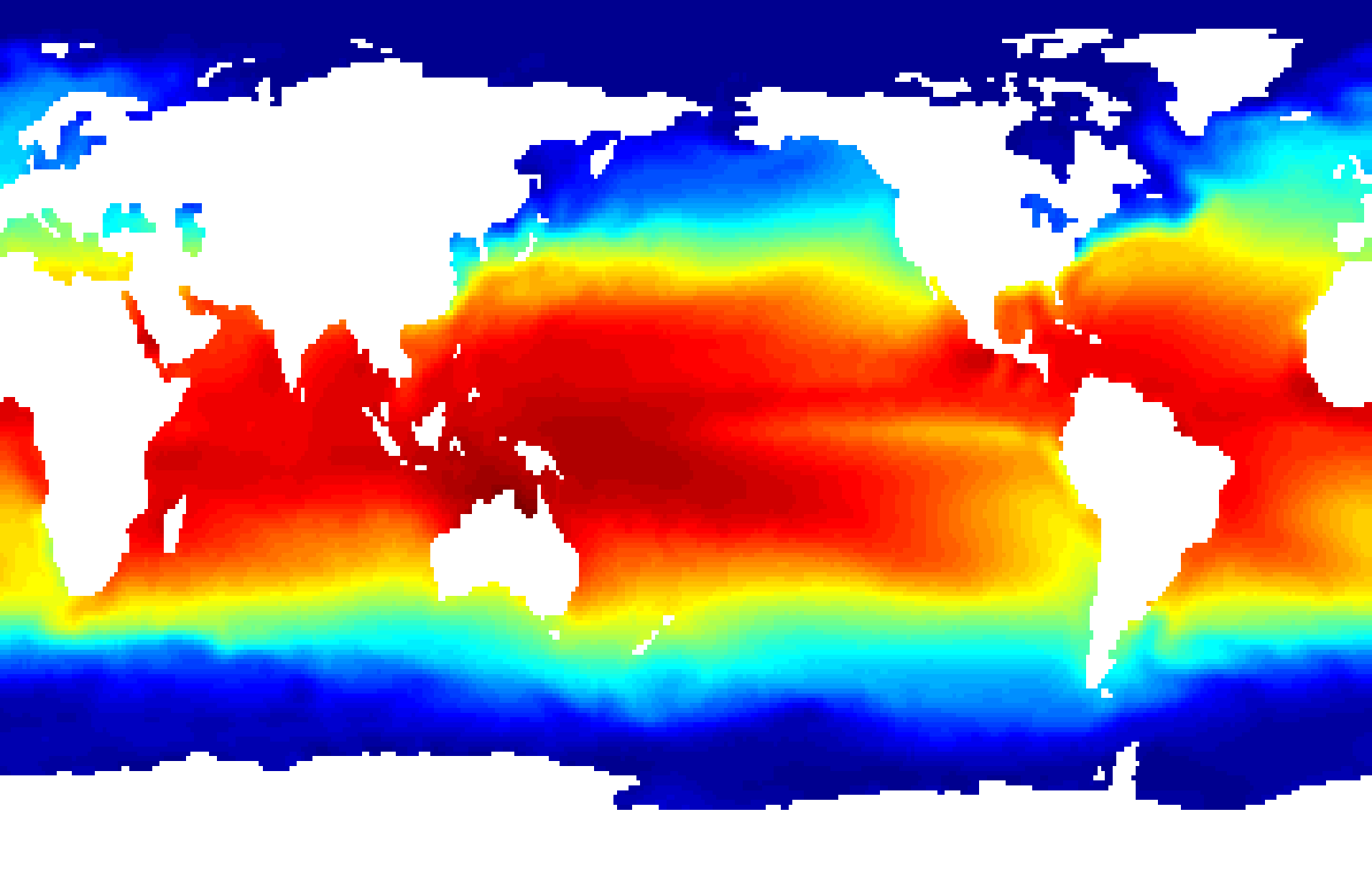}} & \fbox{\includegraphics[width=.25\textwidth]{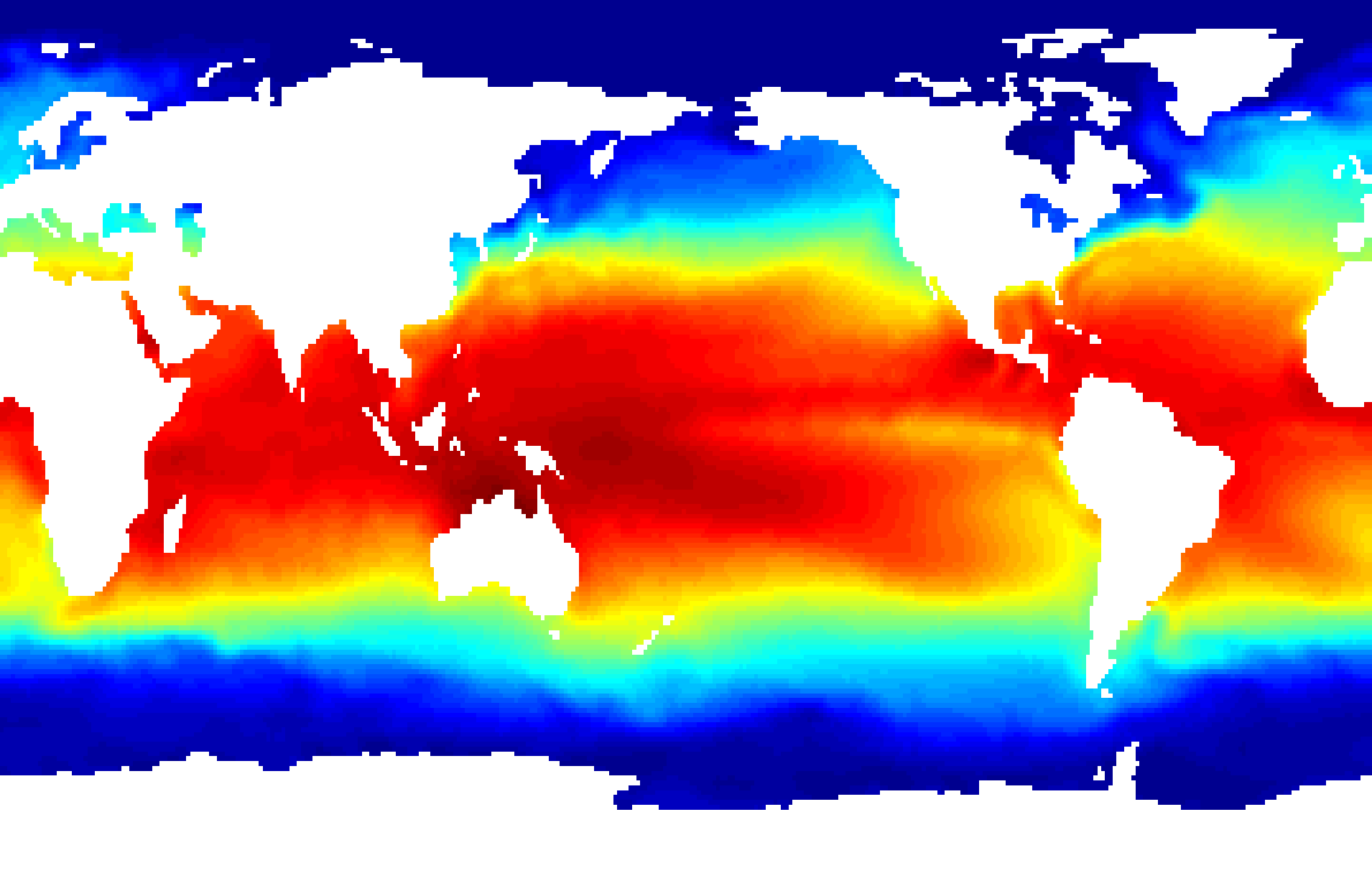}} & \fbox{\includegraphics[width=.25\textwidth]{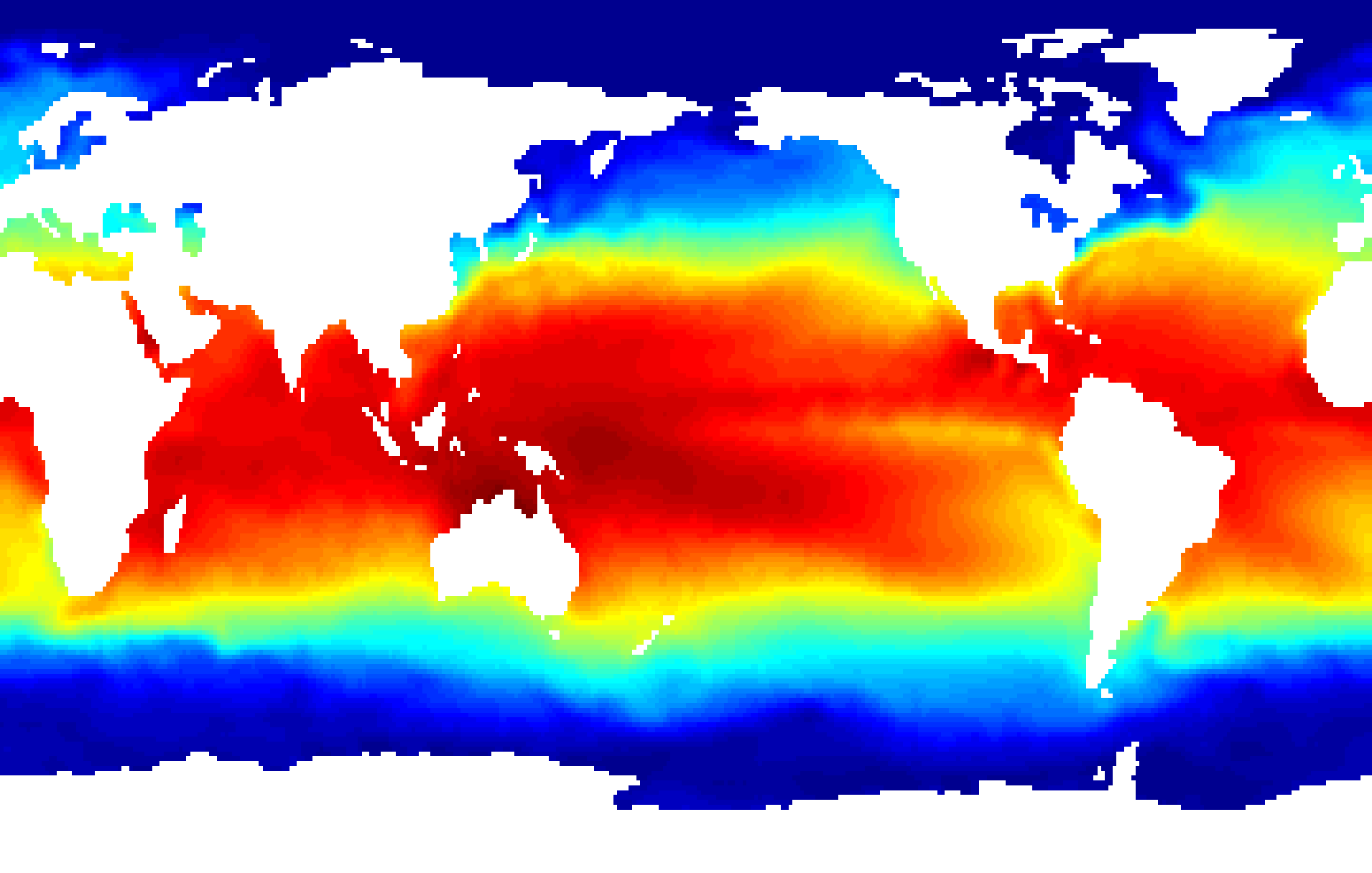}} \\ \hline \hline
		  & {\bf 100 sensors} &  {\bf 200 sensors} & {\bf  \qquad\quad 302 sensors} \\
		\begin{sideways}{\qquad\bf $\ell_2$ QR}\end{sideways} & \fbox{\includegraphics[width=.25\textwidth]{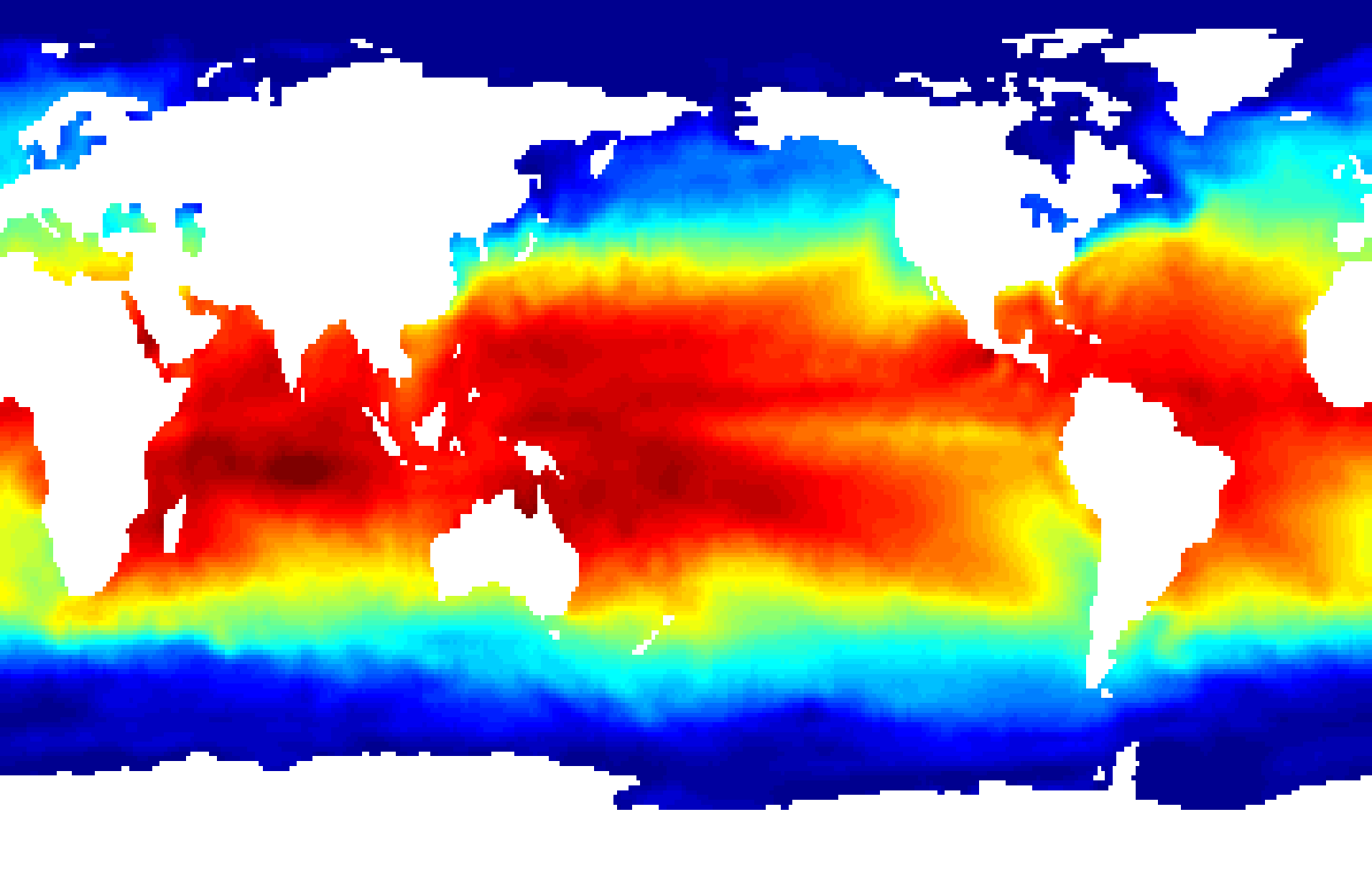}} & \fbox{\includegraphics[width=.25\textwidth]{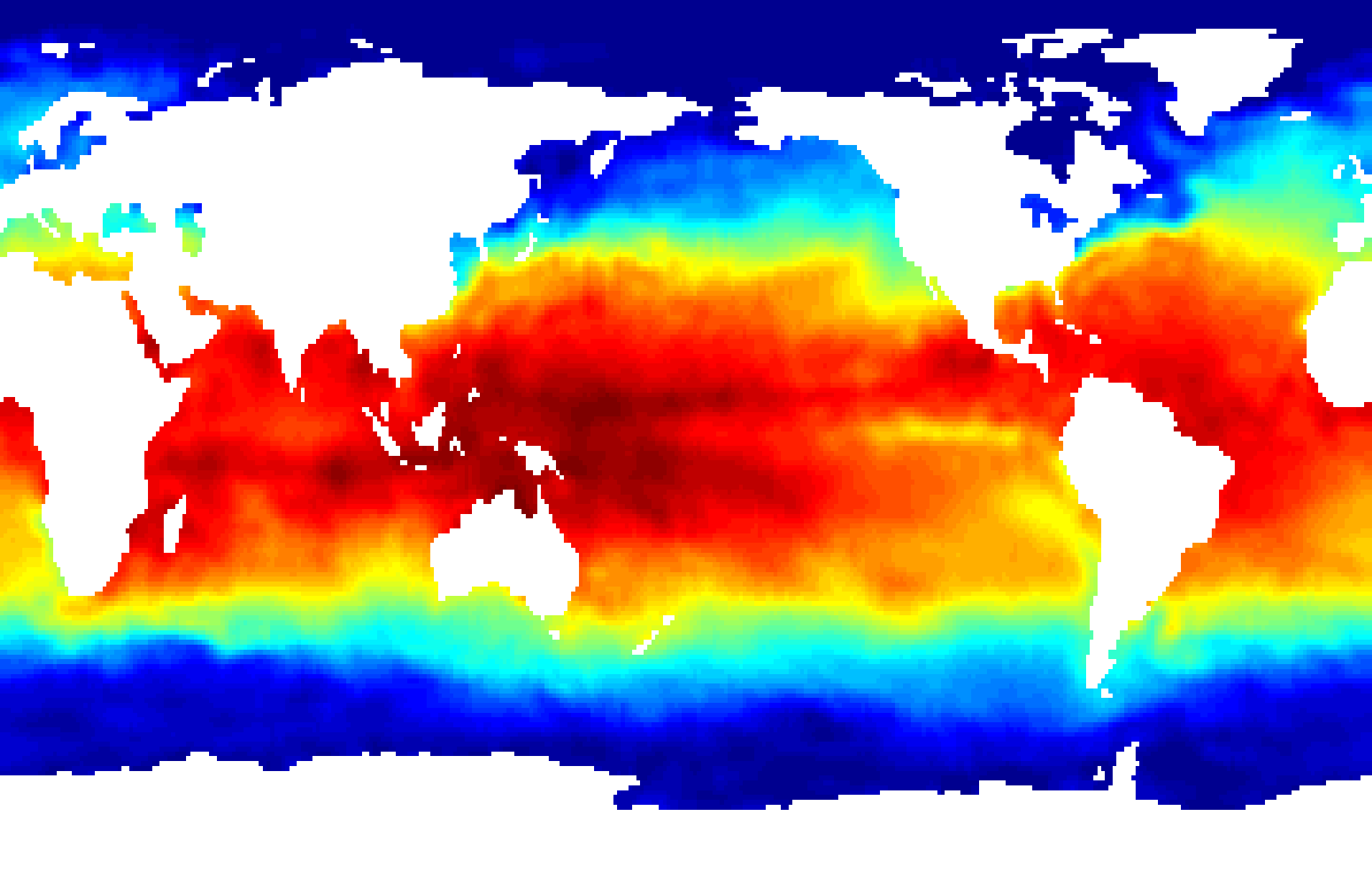}} & \fbox{\includegraphics[width=.25\textwidth]{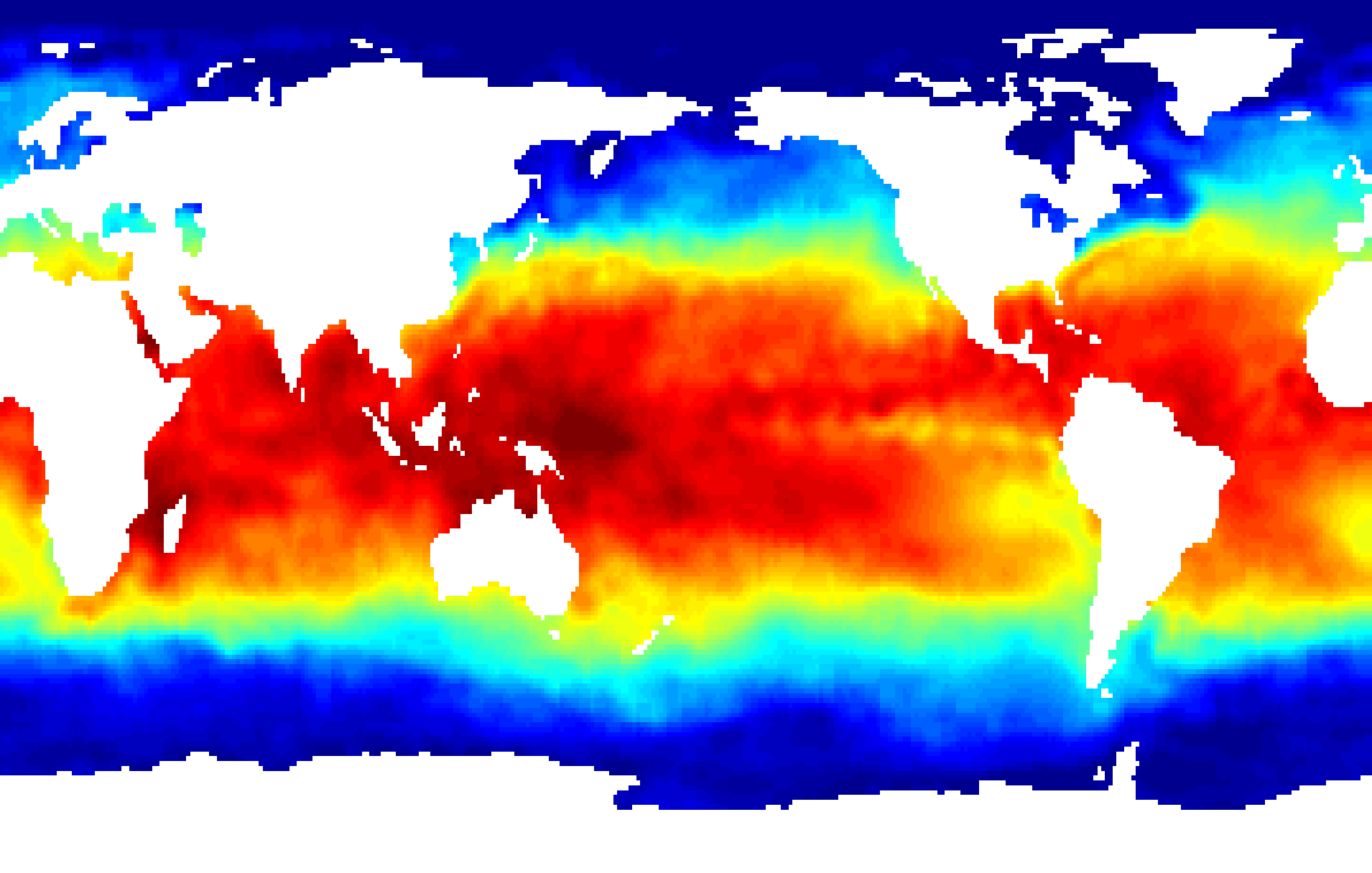}} \\
		\begin{sideways}{\quad\bf $\ell_2$ Random}\end{sideways} & \fbox{\includegraphics[width=.25\textwidth]{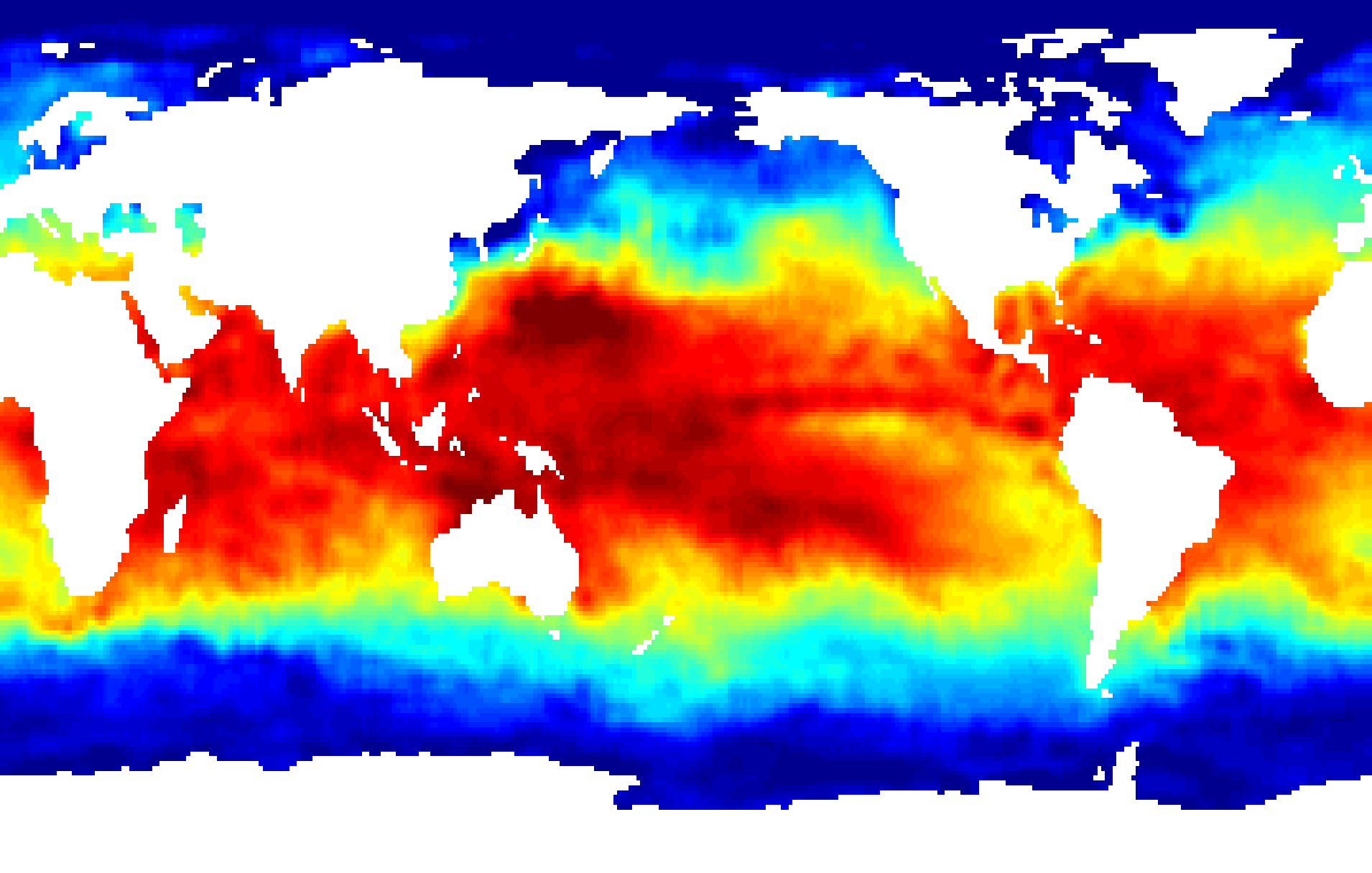}} & \fbox{\includegraphics[width=.25\textwidth]{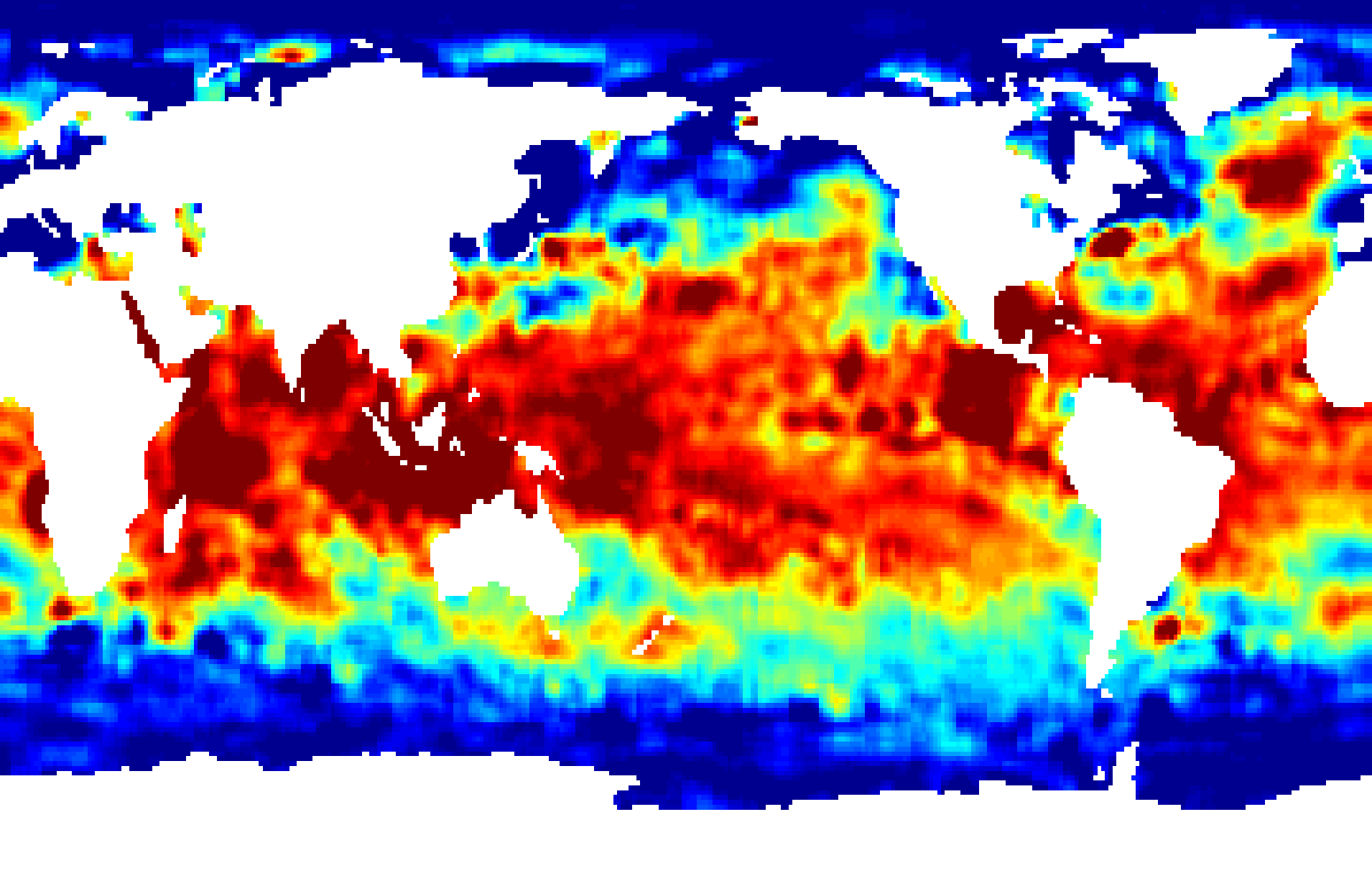}} & \fbox{\includegraphics[width=.25\textwidth]{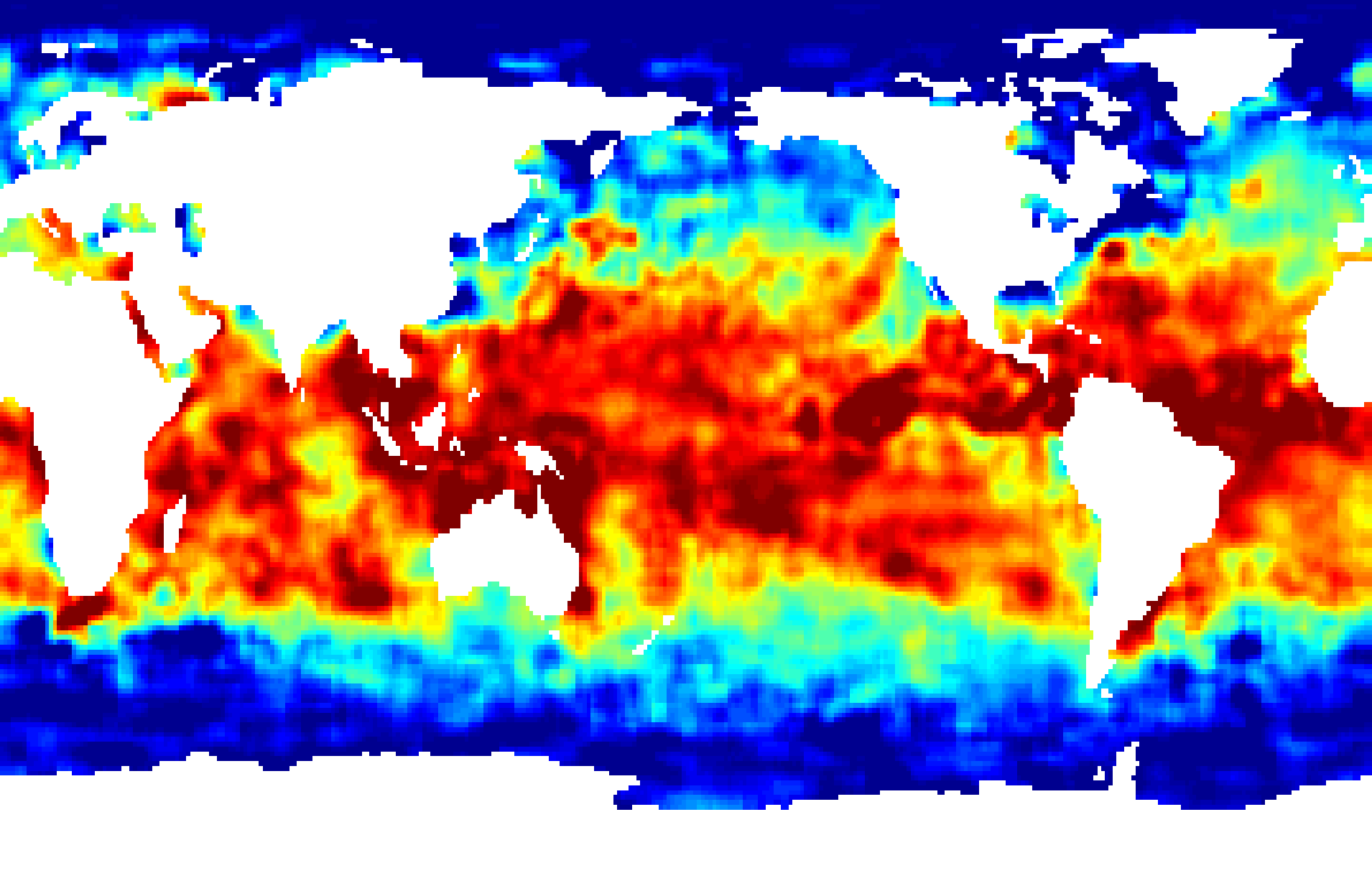}} \\
		\begin{sideways}{\bf Compressed}\end{sideways}\begin{sideways}{\bf Sensing}\end{sideways} & \fbox{\includegraphics[width=.25\textwidth]{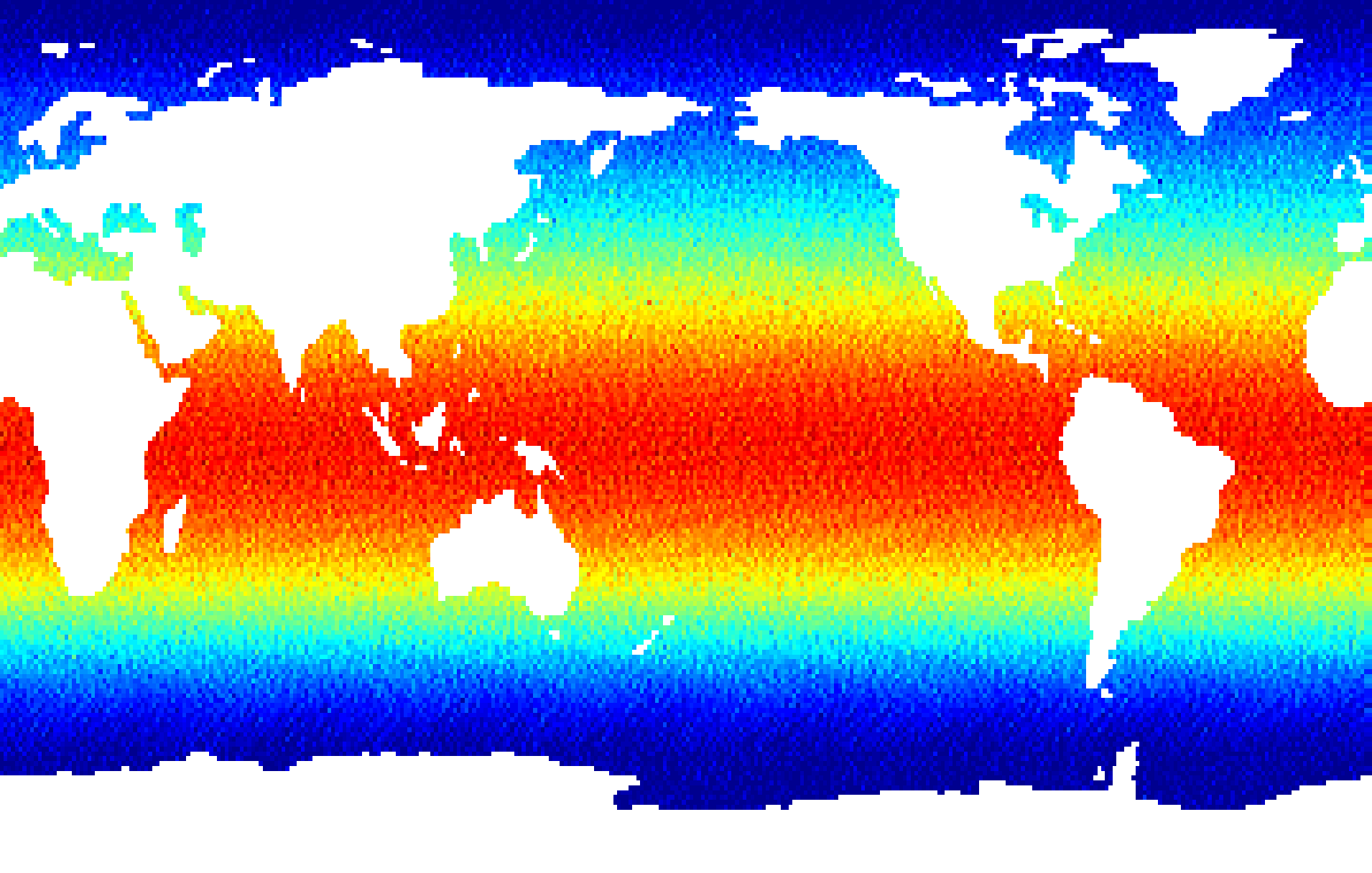}} & \fbox{\includegraphics[width=.25\textwidth]{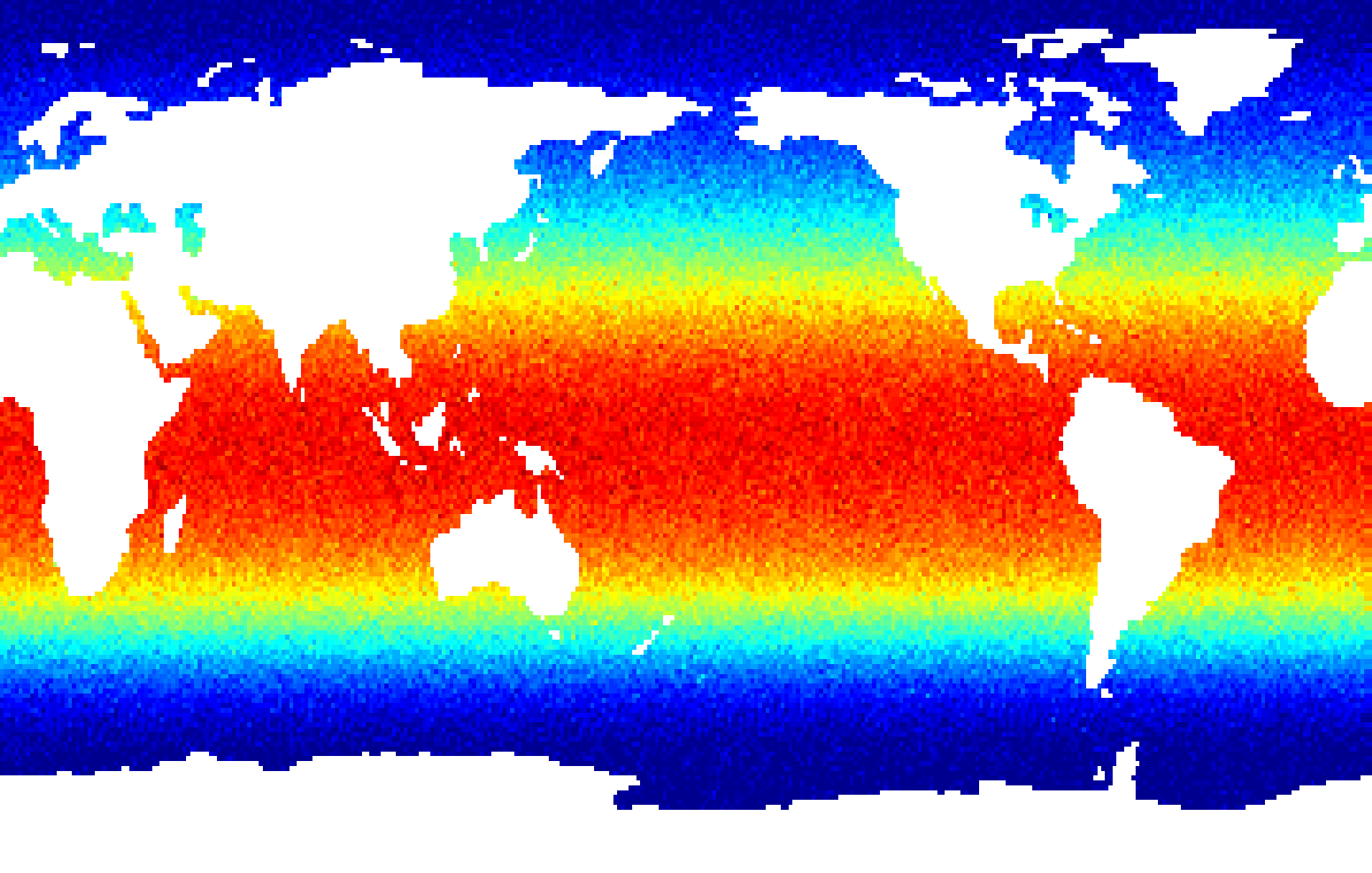}} & \fbox{\includegraphics[width=.25\textwidth]{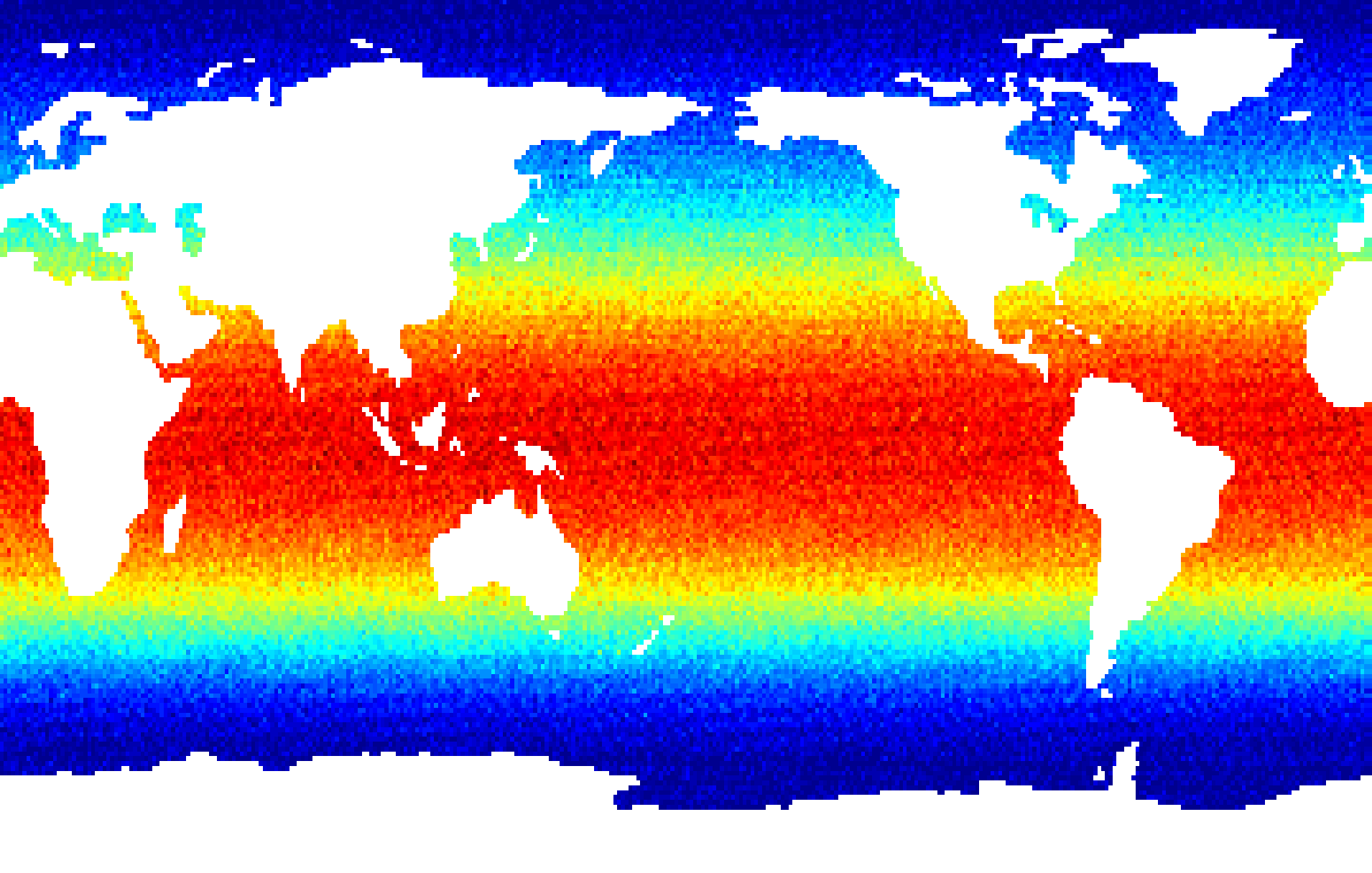}} \\ \hline
	\end{tabular}
	\caption{Comparison of methods for reconstruction of a single snapshot from sensors. QR selected sensors filter uninformative features and achieve better reconstruction. In comparison, random sensors achieve poor reconstruction with least squares ($\ell_2$) and compressed sensing. As the number of proper orthogonal decomposition (POD) modes is increased, modal approximation with the full state (top row) only gets marginally better, which indicates additional features contribute low-energy uninformative features. Hence $\ell_2$ reconstruction expresses low-energy POD modes and suffers from overfitting, with both random and to a lesser extent, QR sensors. \label{fig:enso_recon}}
\end{figure*}

Next we consider the {NOAA\_OISST\_V2} global ocean surface temperature dataset spanning the duration 1990--2016. The data is publicly available online~\cite{noaa_oisst_v2}.  
Unlike eigenfaces, this dataset is a time series, for which a snapshot is recorded every week. Sensor selection must then track energetic temporal signatures. 
Sensors and features are trained on the first 16 years (832 snapshots), and a test snapshot is selected from the excluded validation set. The singular values are shown in Fig.~\ref{fig:enso_spect}.

Like the eigenfaces, localized convective phenomena have energetic contributions to otherwise globally uninformative {\em eigenssts}. This is best seen in the POD snapshot projections, in which the 100 eigensst projection already sufficiently recovers dynamics, while increasing the number of eigenssts in the projection further refines convective phenomena. These lower-energy modes containing convective effects contribute to some degree of overfitting in $\ell_2$ reconstruction (Fig.~\ref{fig:enso_sensors}). The most interesting of these is the El Ni\~no southern oscillation (ENSO) feature that is clearly identified from QR selected sensors. El Ni\~no is defined as any temperature increase of a specified threshold lasting greater than six months in this highlighted region of the South Pacific. It has been implicated in global weather patterns and climate change. 

{\em Remark}: Modal separation of intermittent phenomena such as the El Ni\~no is difficult from a time-invariant POD analysis. Separation of isolated, low-energy temporal events cannot be done from a variance-characterizing decomposition such as the POD -- reordering the snapshots will yield the same dominant modes. On the other hand, tensor decompositions and temporal-frequency analyses such as multiresolution dynamic mode decomposition have succeeded at identifying El Ni\~no where POD has failed. Sensor selection using non-normal modes arising from such decompositions remains an open problem and the focus of ongoing work.

\section{Discussion}\label{Sec:Discussion}
The efficient sensing of complex systems is an important challenge across the physical, biological, social, and engineering sciences, with significant implications for nearly all downstream tasks.  
In this work, we have demonstrated the practical implementation of several sparse sensing algorithms on a number of relevant real-world examples.  
As discussed throughout, there is no all-purpose strategy for the sparse sensing of a high-dimensional system.  
Instead, the choice depends on key factors such as the amount of training data available, the scope and focus of the desired estimation task, cost constraints on the sensors themselves, and the required latency of computations on sensor data.  
Thus, we partition the sparse sensing algorithms into two fundamental categories:  1) optimized sensing in a data-driven tailored basis, and 2) random sensing in a universal basis.

A critical comparison of the two approaches highlights a number of relative strengths and weaknesses.  
The first strategy results in a highly optimized set of sensors that are suitable for tightly scoped reconstruction problems where sufficient training data is available.  
The second strategy requires more sensors for accurate reconstruction but also makes fewer assumptions about the underlying signal, making it more general.  
We emphasize that optimized sensing in a tailored basis typically provides more accurate signal reconstruction than random measurements, facilitating a reduction in the number of sensors by about a factor of two.  
Further, sensor selection and signal reconstruction in the tailored basis is computationally efficient and simple to implement, while compressed sensing generally requires a costly iterative algorithm.  
In problems where sensors are expensive, or when low-latency decisions are required, the reduction in the number of sensors and the associated speed-up of optimized sensing can be significant.  
Thus, when the reconstruction task is well-scoped and a sufficient volume of training data is available, we advocate principled sensor selection rather than compressed sensing.  
In addition, pivoted QR sensors may be used in conjunction with other tailored bases (polynomials, radial basis functions) when signal structure is known. Since these are not data-driven basis functions, QR optimized samples can generalize to different dynamical regimes or flow geometries.

\subsection{Potential applied impact}
Many fields in science and engineering rely on sensing and imaging.  
Moreover, any application involving feedback control for stabilization, performance enhancement, or disturbance rejection relies critically on the choice of sensors.  
We may roughly categorize these sensor-critical problems into two broad categories:  1) problems where sensors are expensive and few (ocean sampling, disease monitoring, espionage, etc.), and 2) problems where sensors are cheap and abundant (cameras, high-performance computation, etc.). 
 
In the first category, where sensors come at a high cost, the need for optimized sparse sensors is clear.  However, it is not always obvious how to collect the training data required to optimize these sensors.  In some applications, high-fidelity simulations may provide insight into coherent structures, whereas in other cases a large-scale survey may be required.  It has recently been shown that it may be possible to optimize sensors based on heavily subsampled data, as long as coherent structures are non-localized~\cite{Brunton2016siap}.  

In the second category, where sensors are readily available, it may still be advantageous to identify key sensors for fast control decisions.  
For example, in mobile applications, such as vision-based control of a quad-rotor or underwater monitoring of an energy harvesting site using an autonomous underwater vehicle, computational and battery resources may be limited.  
Restricting high-dimensional measurements to a small subset of key pixels speeds up computation and reduces power consumption.  
Similar performance enhancements are already exploited in high-performance computing, where expensive function evaluations are avoided by sampling at key interpolation points~\cite{Chaturantabut2012siamjna,Chaturantabut2010siamjsc}.  
Finally, it may also be the case that if measurements are corrupted by noise, reconstruction may improve if uninformative sensors are selectively ignored.  

\section*{Reproducible Research}
A Matlab code supplement is available~\cite{manohar2017code} for reproducing results in this manuscript, including:
\begin{enumerate}
	\item Datasets in Matlab file formats, or links to data that are publicly available online;
	\item Matlab scripts to recreate figures of results.
\end{enumerate}

\footnotesize{
\bibliographystyle{ieeetr}
\bibliography{references}
}

\newpage
\section*{Acknowledgments}
\normalsize
The authors thank Eurika Kaiser, Joshua Proctor, Serkan Gugercin, Bernd Noack, Tom Hogan, Joel Tropp, and Aleksandr Aravkin for valuable discussions.  
SLB and JNK acknowledge support from the Defense Advanced Research Projects Agency (DARPA HR0011-16-C-0016).  
SLB and KM acknowledge support from the Boeing Corporation (SSOW-BRT-W0714-0004).  
BWB and SLB acknowledge support from the Air Force Research Labs (FA8651-16-1-0003).
\revision[]{SLB acknowledges support from the Air Force Office of Scientific Research (AFOSR FA9550-16-1-0650). JNK acknowledges support from the Air Force Office of Scientific Research (AFOSR FA9550-15-1-0385).}{}


\sidebars 

\end{document}